%% file: main.tex
\documentclass[onefignum,onetabnum]{siamart250211}

\input{ex_shared}

\input{pgmacros}

\ifpdf
\hypersetup{
  pdftitle={Fast Direct Solvers},
  pdfauthor={Per-Gunnar Martinsson and Michael O'Neil}
}
\fi

\newcommand{\ignore}[1]{}

\begin{document}

\maketitle

\begin{abstract}
  This survey describes a class of methods known as \emph{fast direct solvers}.
  These algorithms address the problem of solving a system of linear equations
  $\mtx{A}\vct{x}=\vct{b}$ arising from the discretization of either an elliptic
  PDE or of an associated integral equation.  The matrix $\mtx{A}$ will be
  sparse when the PDE is discretized directly, and dense when an integral
  equation formulation is used.  In either case, industry practice for large
  scale problems has for decades been to use iterative solvers such as
  multigrid, GMRES, or conjugate gradients.  A direct solver, in contrast,
  builds an approximation to the inverse of $\mtx{A}$, or alternatively, an
  easily invertible factorization (e.g.~LU or Cholesky).  A major development in
  numerical analysis in the last couple of decades has been the emergence of
  algorithms for constructing such factorizations or performing such inversions
  in linear or close to linear time.  Such methods must necessarily exploit that
  the matrix~$\mtx{A}^{-1}$ is \emph{data-sparse}, typically in the sense that
  it can be tessellated into blocks that have low numerical rank.  This survey
  provides a unifying context to both sparse and dense fast direct solvers,
  introduces key concepts with a minimum of notational overhead, and provides
  guidance to help a user determine the best method to use for a given
  application.
\end{abstract}

\begin{keywords}
  Fast solvers; 
  elliptic PDEs; 
  fast multipole methods; 
  $\mathcal{H}$-matrix methods;
  integral equation methods;
  boundary integral equations;
  multigrid;
  preconditioners; 
  rank structured matrices;
  Krylov methods;
  hierarchically off-diagonal low rank matrices;
  hierarchically block separable matrices; 
  hierarchically semi-separable matrices; 
  data-sparsity; 
  randomized algorithm; 
  sparse direct solver.
\end{keywords}

\begin{AMS}
65-02; 65F05; 65F55; 65N22; 65N38; 65N80.
\end{AMS}

\tableofcontents

\section{Introduction}
\label{sec:intro}
\input{01-intro/intro03}


\section{A one-dimensional boundary value problem}
\label{sec:onedimbvp}
\input{02-onedimbvp/onedimbvp02}


\section{A one-dimensional boundary value problem via an integral equation}
\label{sec:onedimbvpIE}
\input{03-onedimbvpIE/onedimbvpIE03}


\section{Hierarchically Off-Diagonal Low Rank matrices}
\label{sec:HODLR}
\input{04-hodlr/HODLRnew04}


\section{LU factorization with nested dissection ordering}
\label{sec:FDSPDE}
\input{05-nesteddiss/nesteddiss06}




\section{Integral equations for  boundary value problems}
\label{sec:BIE2D}

\input{07-BIE2D/BIE2D02}


\section{Rank structured matrices with nested bases}
\label{sec:BIE_nested}
\input{08-nested/nested_new03}


\section{Rank structured matrices with strong admissibility}
\label{sec:BIE_strong}
\input{09-strong/strong_new02}


\section{Functions of matrices, operator algebra, and preconditioning}
\label{sec:operator_algebra}
\input{operator_algebra3}


\section{The cause of rank deficiencies: Dissipation of information}
\label{sec:dissipation}
\input{10-dissipation/dissipation_new03}


\section{How to obtain the compressed representation}
\label{sec:compression}
\input{11-compression/compression02.tex}


\section{Rank structured matrices in data science and computational statistics}
\label{sec:statistics}

\input{12-statistics/statistics01.tex}


\section{A roadmap of rank structured formats}

\label{sec:taxonomy}
\input{13-taxonomy/taxonomy03.tex}


\section{Conclusions and future outlook}
\label{sec:future}
\input{14-conclusion/conclusion02.tex}

\bibliography{references}
\bibliographystyle{siamplain}
\end{document}

%% file: ex_shared.tex
\usepackage{lipsum}
\usepackage{amsfonts}
\usepackage{graphicx}
\usepackage{epstopdf}
\usepackage{algorithmic}
\ifpdf
  \DeclareGraphicsExtensions{.eps,.pdf,.png,.jpg}
\else
  \DeclareGraphicsExtensions{.eps}
\fi

\newsiamremark{remark}{Remark}
\newsiamremark{hypothesis}{Hypothesis}
\crefname{hypothesis}{Hypothesis}{Hypotheses}
\newsiamthm{claim}{Claim}

\headers{Fast direct solvers}{P.G.~Martinsson and M.~O'Neil}

\title{Fast direct solvers\thanks{The authors wish to acknowledge the valuable
    contributions of members of the fast direct solver community who provided
    feedback on preliminary versions of this manuscript. Particular thanks are
    due to David Persson of the Courant Institute at NYU for his exceptionally
    meticulous and comprehensive review. \funding{Martinsson was supported in
      part by the Office of Naval Research (N00014-18-1-2354), the National
      Science Foundation (DMS-2313434 and DMS-1952735), and the Department of
      Energy ASCR (DE-SC0022251). O'Neil was supported in part by the Office of
      Naval Research (N00014-18-1-2307 and N00014-21-1-2383.)}  } }

\author{Per-Gunnar Martinsson\thanks{Department of Mathematics, University of Texas at Austin,
  (\email{pgm@oden.utexas.edu}, \url{https://users.oden.utexas.edu/~pgm/}).}
\and Michael O'Neil\thanks{Courant Institute, New York University,
  (\email{oneil@cims.nyu.edu}, \url{https://cims.nyu.edu/~oneil}).}
  }

\usepackage{amsopn}


%% file: pgmacros.tex
\usepackage[T1]{fontenc}

\usepackage{color}
\usepackage{graphicx}
\usepackage{bm}
\usepackage{tabto}

\usepackage{pgfplots}

\usepackage[font=footnotesize]{caption}
\usepackage[font=footnotesize]{subcaption}
\usepackage{adjustbox}

\usepackage{kbordermatrix}

\input xy
\xyoption{all}

\usepackage[sort]{cite}

\newcommand{\mtx}[1]{\bm{\mathsf{#1}}}

\newcommand{\pvct}[1]{\bm{#1}}
\newcommand{\vct}[1]{\bm{\mathsf{#1}}}

\newcommand{\mtwo}[4]{\left[\begin{array}{cc} #1 & #2 \\ #3 & #4 \end{array}\right]}
\newcommand{\vtwo}[2]{\left[\begin{array}{cc} #1 \\ #2 \end{array}\right]}

\newcommand{\pcc}{\pvct{c}}

\newcommand{\pnn}{\pvct{n}}

\newcommand{\pxx}{\pvct{x}}
\newcommand{\pww}{\pvct{w}}
\newcommand{\pyy}{\pvct{y}}

\newcommand{\mA}{\mtx{A}}

\newcommand{\mD}{\mtx{D}}
\newcommand{\mE}{\mtx{E}}
\newcommand{\mF}{\mtx{F}}
\newcommand{\mG}{\mtx{G}}

\newcommand{\mU}{\mtx{U}}
\newcommand{\mV}{\mtx{V}}

\newcommand{\mhD}{\hat{\mtx{D}}}

\newcommand{\mhG}{\hat{\mtx{G}}}

\newcommand{\mtA}{\tilde{\mtx{A}}}

\newcommand{\mtD}{\tilde{\mtx{D}}}

\definecolor{MyDarkGreen}{rgb}{0,0.6,0}
\definecolor{MyGray}{rgb}{0.6,0.6,0.6}

\newcommand{\mtb}[1]{\mathcal{U}}
\newcommand{\cO}{\mathcal{O}}
\newcommand{\cA}{\mathcal{A}}

\newcommand{\cV}{\mathcal{V}}
\newcommand{\cD}{\mathcal{D}}
\newcommand{\cI}{\mathcal{I}}

\newcommand{\cY}{\mathcal{Y}}
\newcommand{\cZ}{\mathcal{Z}}

\newcommand{\bbR}{\mathbb R}

%% file: 01-intro/intro03.tex
This survey concerns the problem of numerically computing an approximate
solution to a linear elliptic PDE such as the Laplace, Helmholtz, or Stokes
equation.  Once a mathematical formulation of the problem has been fixed
(i.e.~PDE vs.~integral equation), the first step in a numerical method is to
discretize the continuum equations to obtain a linear system
\begin{equation}
\label{eq:basic}
\mtx{A}\vct{x} = \vct{b},
\end{equation}
where $\mtx{A}$ is an $N\times N$ matrix and the vector~$\vct{b}$ encodes any
boundary conditions and body loads that may have been specified in the
formulation of the problem.  This survey describes direct (as opposed to
iterative) methods for solving the system~(\ref{eq:basic}).

When the mathematical formulation is a partial differential equation, and the
discretization is done through a finite element or finite difference scheme, the
matrix $\mtx{A}$ will be \textit{sparse}, reflecting the fact that it
discretizes a \textit{local} operator. In this context, a direct solver has
traditionally been thought of as a process to compute a Cholesky or LU
factorization of the coefficient matrix that is exact (modulo rounding errors).
The factors computed will be less sparse than the original matrix, however,
which means that the process will not scale linearly with~$N$. With a well
chosen elimination order based on an analysis of the sparsity structure (usually
viewed as a graph of connectivity), complexity is
typically~$\mathcal O(N^{3/2})$ or~$\mathcal O(N^2)$ for problems in two or
three dimensions, respectively~\cite{george_1973, hoffman_1973,
  1989_directbook_duff, 2006_davis_directsolverbook}.  When~$N$ grows large,
such methods quickly become unfeasible and standard practice is to switch to an
iterative solver such as
multigrid~\cite{2016_hackbusch_iterative_book,mccormick_multigrid}, or
(preconditioned) GMRES~\cite{1986_saad_gmres}, which attain $\mathcal O(N)$
complexity in many environments.

On the other hand, when the mathematical formulation is an integral equation,
the continuum equations are~\textit{global} and basically every discretization
scheme will result in a \textit{dense} coefficient matrix~$\mtx{A}$.  A direct
solver in this context typically consists of either computing an LU
factorization of $\mtx{A}$, or just computing the inverse~$\mtx{A}^{-1}$
directly, at a cost of~$\mathcal O(N^3)$ in either case.  This cubic scaling
severely limits the maximum~$N$ that can be handled. Fortunately, for many integral equations arising from classical mathematical physics, there exist
fast algorithms such as the Fast Multipole
Method~\cite{rokhlin_1990,rokhlin1997,1999_ChengGreengardRokhlin} that are
capable of applying the matrix~$\mtx{A}$ to a vector in
complexity~$\mathcal O(N)$.  (To be precise, the operation is carried out to
some specified precision~$\varepsilon$, and the complexity scales
as~$\mathcal O\bigl(\log(1/\varepsilon)^j\,N\bigr)$ for some small integer $j$.)
Such methods have given rise to iterative solvers with overall linear complexity
in $N$ whenever convergence of the iterative method is fast.

The purpose of this text is to describe a class of algorithms that can compute
an approximate invertible factorization of the coefficient matrix~$\mtx{A}$ in
linear or close to linear time. These algorithms exploit the fact that while the
inverse is basically always dense, it is typically \textit{data-sparse} in the
sense that it contains only~$\mathcal O(N)$ pieces of information to any given
numerical precision. The particular form of data-sparsity we focus on here is
\emph{rank structure}, whereby the matrix can be tessellated
into~$\mathcal O(N)$ blocks in such a way that each block has low numerical
rank. (In what follows, we will often refer to \textit{low-rank} factorizations
-- in almost all cases we will mean \textit{approximate low-rank}
factorizations, in which case we say that a matrix has low \textit{numerical
  rank}.)

Traditionally, computing an approximate numerical solution to a PDE has been
viewed as consisting of two separate tasks: (1) A discretization step that
relies on ideas from analysis and PDE theory, and then (2) a solver step that
relies on methods from linear algebra or computer science.  In the direct
solvers that we describe, knowledge of the behavior of elliptic PDEs infuses
both steps. Specifically, \emph{a priori} knowledge of which submatrices can be
expected to have low numerical rank is gained directly from the mathematical
understanding of the behavior of solutions to elliptic PDEs and the elliptic
differential operators themselves. Using this a priori analytical knowledge, the
solver step can further be split into three pieces: (2a) a \textit{compression}
step whereby the matrix~$\mtx{A}$ is factorized into a data-sparse form (usually
exploiting approximate low-rank structures). (2b) An \emph{inversion} (or
\emph{factorization}) stage whereby the data structure constructed in the first
step is used to directly invert the operator. (2c) A \emph{solve} stage where
the computed inverse is applied to given data in the form of boundary conditions
or body loads to produce a solution.

In general, the direct solvers that we describe come with an additional layer of
algorithmic and mathematical complexity compared to classical methods, and they
also tend to be more memory intensive than iterative methods.
Despite these drawbacks, there are important environments where
direct solvers are worth the extra cost:
\begin{itemize}
\item 
Direct solvers provide a path towards solving PDEs that are intractable to iterative methods, in particular for problems with oscillatory solutions.
\item 
Dramatic acceleration often occurs when solving multiple equations with identical coefficient matrices, as the cost of building the approximate inverse is amortized across many solves.
\item Direct solvers are well-suited to many modern communication-constrained
  and heterogeneous computing architectures. On the one hand, they involve
  fewer passes through hierarchical trees than most iterative methods. On the
  other, direct solvers tend to spend most flops in dense linear algebra
  operations, which are very efficient on both CPUs and GPUs.
\item Lastly, since direct solvers give access to the inverse operator
  they can be easily integrated into multiphysics solvers in which
  \emph{operator algebra} forms the basis of many PDE formulations.
  They enable the explicit formation of objects such as time evolution
  operators or Dirichlet-to-Neumann maps, and can often be used to
    construct efficient preconditioners based on the underlying physics of the
    problem.
\end{itemize}
We now turn to briefly introducing some of the themes that will then be pursued
in detail throughout the main part of the survey.

\subsection{Solution formulas for linear elliptic PDEs}

A direct solver can be viewed as a technique for numerically building an
approximation to the solution operator for a PDE. To describe some of the relevant
properties of such operators, let us start by considering a particularly simple
model problem, namely the free-space Poisson equation in three dimensions,

\begin{equation}
\label{eq:laplace}
-\Delta\,u(\pvct{x}) = g(\pvct{x}),   \qquad  \pvct{x} \in \mathbb{R}^{3},
\end{equation}
where $g$ is a given body load, and where
$\Delta u = \sum_{j=1}^{3} \partial^{2}u/\partial x_{j}^{2}$, as
usual. When~(\ref{eq:laplace}) is coupled with appropriate decay conditions
at infinity, the unique solution is 
\begin{equation}
  \label{eq:laplacesoln}
    u(\pvct{x}) = \cV[ g](\pvct{x}) =
\int_{\mathbb{R}^{3}}\phi(\pvct{x} - \pvct{y})\,g(\pvct{y})\,d\pvct{y},
\qquad\pvct{x} \in \mathbb{R}^{3},
\end{equation}
where $\phi$ is the fundamental solution of Laplace's equation:
\begin{equation}
\label{eq:fundlaplace}
\phi(\pvct{x}) = - \frac{1}{4\pi|\pvct{x}|}.
\end{equation}
The operator~$\cV$ is a \emph{volume potential operator}.
Observe that the forward operator in~(\ref{eq:laplace}) is local;
however, the solution operator in~(\ref{eq:laplacesoln}) is global.

When (\ref{eq:laplacesoln}) is evaluated numerically, after
discretization one typically
ends up needing to evaluate a sum of the form
\begin{equation}
\label{eq:nbody}
u_{i} = \sum_{j\neq i}^{N} w_{ij}\,\phi(\pvct{x}_{i} - \pvct{x}_{j})
\,g(\pvct{x}_{j}),
\end{equation}
where $\{\pxx_{i}\}_{i=1}^{N}$ are some discretization nodes and
$\{w_{ij}\}_{i,j=1}^{N}$ are the associated quadrature
weights\footnote{Determining quadrature weights that lead to
high-order accuracy turns out to be a knottier problem than one might have
hoped, since the kernel in (\ref{eq:laplacesoln}) is singular as
$|\pxx - \pyy|\rightarrow 0$. However, the problem is well-understood
and does not fundamentally contribute to the difficulty or asymptotic
complexity of evaluating the sum~(\ref{eq:laplacesoln}) since almost
every scheme leads to weights that are approximately separable in the
sense that $w_{ij} = \tilde{w}_{i}\hat{w}_{j}$ for the vast majority
of elements.}.  While the direct evaluation of~\eqref{eq:nbody} for
each~$i$ would naively require~$\mathcal O(N^2)$ operations, algorithms such as fast
multipole methods (FMMs) enable evaluation of the sum in~$\mathcal
O(N)$ complexity, resulting in a fast solver for (\ref{eq:laplace}).

For more general elliptic PDEs, solution operators analogous
to~(\ref{eq:laplacesoln}) formally exist in a mathematical sense, but can rarely
be expressed explicitly. For instance, consider an inhomogenous elliptic
boundary value problem
\begin{equation}
  \label{eq:bvp}
  \begin{aligned}
    \mathcal A [u](\pvct{x}) =&\ g(\pvct{x}),   &\quad  &\pvct{x} \in \Omega,\\
    \mathcal B [u](\pvct{x}) =&\ f(\pvct{x}),   &    &\pvct{x} \in \Gamma,
  \end{aligned}
\end{equation}
where $\Omega$ is a domain with boundary $\Gamma$,~$\mathcal A$ is a linear
elliptic differential operator (with possibly variable coefficients),
and~$\mathcal B$ specifies the boundary conditions on~$u$
along~$\Gamma$.  The solution formula for~(\ref{eq:bvp}) would generally take
the form
\begin{equation}
  \label{eq:solutionop}
  \begin{aligned}
    u(\pxx) &= \mathcal G[g](\pxx) + \mathcal F[f](\pxx) \\
&= \int_{\Omega}G(\pxx,\pyy)\,g(\pyy)\,d\pyy +
  \int_{\Gamma}F(\pxx,\pyy)\,f(\pyy)\,dS(\pyy),\qquad \pxx \in \Omega,
  \end{aligned}
\end{equation}
where $G$ and $F$ are two kernel functions that depend on $\mathcal A$,
$\mathcal B$, and~$\Omega$.  The solution operator in~(\ref{eq:solutionop}) is
typically mathematically well-behaved; the challenge is that~$G$ and~$F$ are 
known analytically only in the simplest possible environments (e.g.~Laplace or
Helmholtz equations on a sphere or a half space).

The purpose of a \textit{fast direct solver} (FDS) is to numerically build an
approximation to a solution operator such as~(\ref{eq:solutionop}) using
mathematical knowledge of the differential operator~$\mathcal A$ and certain
regularity properties of the solution~$u$.  This problem is in many ways more
mathematically tractable than the traditional approach of first discretizing the
PDE and then solving the resulting linear system.  The reason is that the
solution operator in~(\ref{eq:solutionop}) is a much \emph{better} behaved
mathematical object than is the original differential operator~$\mathcal A$.  In
particular, (1) it is a smoothing operator (even though kernel functions are typically 
singular at the diagonal), (2) it is often bounded with respect to natural norms (e.g.~the
$L^2$ norm on~$\Omega$ or~$\Gamma$), and (3) it is generally stable in the sense
that small perturbations in the inputs~$f$ and~$g$ lead to small changes in the
output.  The opposite could be said of the original differential
operator~$\mathcal A$: (1) it is an anti-smoothing operator, (2) it is
unbounded, and (3) small perturbations in its input can lead to large changes in
the output.

\subsection{Fast algorithms for applying integral operators}
\label{sec:fastIE}

That a global operator such as~(\ref{eq:solutionop}) can rapidly be applied to
vector has, for special cases, been a well established fact for several decades.
One of the most celebrated algorithms developed in the past 40 years, the
\textit{Fast Multipole Method} (FMM)~\cite{rokhlin1985,rokhlin1987},
specifically addresses the problem of efficiently evaluating a sum such
as~(\ref{eq:nbody}). This directly enables the fast evaluation of a solution
formula such as~(\ref{eq:laplacesoln}), and also opens a path to linear
complexity solvers for integral equation formulations of many boundary value
problems via iterative methods.

\begin{figure}[!t]
  \begin{center}
    \setlength{\unitlength}{1mm}
\begin{picture}(120,40)
  \put(080,00){\includegraphics[height=40mm]{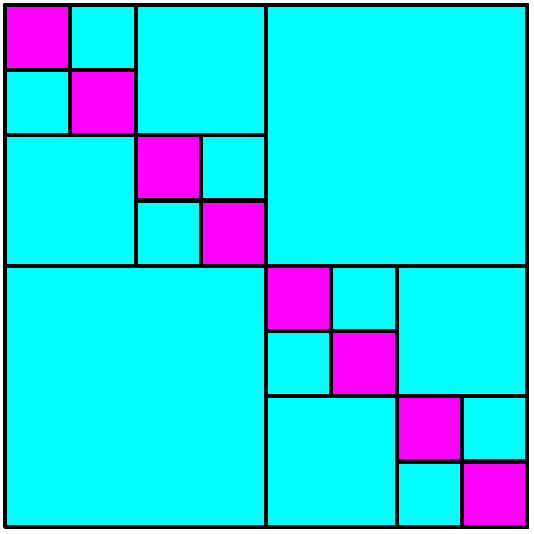}}
  \put(-5,20){\begin{minipage}{81mm}
      \textit{\small{The figure shows an
          example of a \textit{rank-structured} matrix.  Each
          off-diagonal block (cyan) has low numerical rank, and each
          diagonal block (magenta) is small.  Matrices like these arise in
          solvers for elliptic PDEs, for problems in data mining and
          statistics involving smooth kernel functions, and many more.
          The tessellation pattern shown is just one example among many
          possible ones.}}

    \vspace{1.5mm}

    \textit{\small{For matrices of this type, efficient algorithms exist for
        matrix-vector multiplication, matrix inversion, LU decomposition, etc.}}
\end{minipage}}
\end{picture}
\end{center}
\caption{A representative rank-structured matrix, such as an
  ``$\mathcal{H}$-matrix'' or ``HODLR matrix''.}
\label{fig:tessellation}
\end{figure}

The ideas underlying the FMM can be generalized to a much broader class of
linear operators, including integral 
operators such as that
in~(\ref{eq:solutionop}).  In many cases, the underlying reason these ideas work
is that the dense matrix arising upon the discretization
of~(\ref{eq:solutionop}) can be tessellated into $\mathcal O(N)$ blocks in such
a way that each block is of low (numerical)
rank~\cite{2019_martinsson_fast_direct_solvers}.  A representative tessellation
pattern is illustrated in Figure~\ref{fig:tessellation}.  Such a \emph{rank
  structured} matrix can be applied to vectors in linear or close to linear
time, and it is also sometimes possible to compute its inverse, or a triangular
factorization in linear or close to linear time.  These questions can be
approached mathematically from a number of different angles, including, e.g.,
wavelets~\cite{BCR,ABCR} and Calder\'on-Zygmund theory~\cite{meyer1997wavelets}.
(We have so far ignored precisely how to \emph{obtain} the low-rank
approximations of such tessellations. We will return to this aspect of the
algorithms later on.)

Some initial foundational work on rank structured matrices was done by Hackbusch
and co-workers using the so called $\mathcal{H}$- and $\mathcal{H}^{2}$-matrix
formats~\cite{hackbusch, 2008_bebendorf_book, 2002_hackbusch_H2,
  2004_borm_hackbusch, 2010_borm_book}.  Recently, much attention has been
devoted to closely related formats such as \textit{Hierarchically Off-Diagonal
  Low Rank} matrices~\cite{2013_darve_FDS} as well as the \textit{Hierarchically
  Semi Separable (HSS)} or \textit{Hierarchically Block Separable (HBS)}
formats~\cite{2005_martinsson_fastdirect,
  2012_martinsson_FDS_survey,2009_xia_superfast, 2013_xia_randomized_FDS,
  2010_gu_xia_HSS}, among others.

\subsection{Dissipation of information}
\label{sec:intro_diss}

The rank deficiencies that we mentioned in Section~\ref{sec:fastIE} result from
fundamental properties of the integral equations of mathematical physics.  To be
precise, they are the linear algebraic manifestation of the fact that while the
solution operator for a coercive elliptic PDE is global, the amount of
\textit{information} that gets transported from one region of the domain to
another (via the kernel of the solution operator) is limited, to a high degree
of numerical approximation.

For a particularly simple illustration of dissipation of information, let us
consider a linear operator~$\cV_\sigma$ that maps a distribution of electric
charges in a \emph{source} domain $\Omega_{\sigma} \subset \mathbb R^2$ to the
resulting electric potential in a \emph{target} domain $\Omega_{\tau} \subset
\mathbb R^2$, as
illustrated in Figure~\ref{fig:sources_targets} later on in the manuscript.  In
other words, a source density $q$ is mapped to a potential $u$ through
\begin{equation}
\label{eq:poteval}
u(\pxx) =
      \cV_\sigma[q](\pxx) = \int_{\Omega_{\sigma}}\phi(\pxx - \pyy) \, q(\pyy)
         \,d\pyy,\qquad\pxx \in \Omega_{\tau},
\end{equation}
where~$\phi(\pxx) = -(2\pi)^{-1}\log|\pxx|$ is the 2D analog of the 3D fundamental solution in~(\ref{eq:fundlaplace}). As we will
see in Section \ref{sec:dissip_fast_summation}, so long as~$\Omega_\sigma$
and~$\Omega_\tau$ are adequately separated, the kernel
in~(\ref{eq:poteval}), admits what is known as a \textit{multipole expansion}
which expresses the kernel as a sum of terms in which the variables are
separated,
$$
\phi(\pxx - \pyy) = \sum_{p=0}^{\infty}b_{p}(\pxx) \, c_{p}(\pyy),\qquad
\pxx \in \Omega_{\tau},\qquad
\pyy \in \Omega_{\sigma}.
$$ The sum converges exponentially fast when~$\Omega_\sigma$ and
$\Omega_\tau$ are well-separated. 
This means that a high accuracy approximation to~$\phi$ can be attained by
truncating the sum to the first few terms.  Suppose that we keep the first $P$
terms in the sum; then in terms of linear algebra, a matrix obtained by
discretizing (\ref{eq:poteval}) can to high accuracy be approximated by a matrix
of rank $P$.  In Section~\ref{sec:dissip_fast_summation}, we provide more
details to this discussion and also show how the situation changes for other
interaction potentials, such as those associated with the Helmholtz equation.

In mechanics, structural engineers apply this idea of diffusion of information
using the concept of the \textit{Saint-Venant principle}~\cite{toupin1965saint}.
In a simplistic version, it may be stated as follows: \textit{The difference
  between the effects of two different, but statically equivalent loads, becomes
  very small at sufficiently large distances from the load.}  Mathematically,
this notion can be made precise through an asymptotic analysis of the
fundamental solutions to the equations of linear elasticity, in a manner similar
to a multipole expansion for the Laplace equation.

Dissipation of information is also a bedrock principle in a range of different
imaging problems (e.g. radar, seismic probing, etc.), where it arises as the
statement that if you remotely observe a scattered field, then information about
any features that are smaller than the wavelength of the irradiating field can
rarely be retrieved.

In this vein, a key objective of this survey is to elucidate how dissipation of
information is a foundational idea that underlies many popular fast solvers.  It
appears explicitly in classical fast algorithms such as the FMM, and implicitly
in much of the more recent literature on fast direct solvers.  We will describe
what mathematical analysis can tell us about the situation, and also point to
empirical evidence that the exponential convergence of low rank approximations
is astonishingly persistent. In fact, we will in Section \ref{sec:dissipation}
show striking similarities in the
spectral properties of dense matrices arising not only in classical integral
operators of mathematical physics, but also ones that arise in sparse direct
solvers, in domain decomposition methods, and many more.

\subsection{Scaling with dimension --- a serious but manageable problem}
\label{sec:dimensionissue}
While algorithms for efficiently applying global integral operators have become
an essential part of the toolbox in scientific computing, they have a
significant shortcoming in that their practical efficiency scales poorly with
dimension. For instance, our archetypical example of the classical FMM for
evaluating interactions between $N$ electrically charged particles can easily
attain very high practical speed for problems in one and two dimensions.  In
three dimensions, high practical speed can also be attained, but this requires
more sophisticated machinery. In fact, the scaling constant suppressed in
the~$\mathcal O(N)$ notation scales worse than $6^d$ -- this means that
algorithms of this type are practical
in dimensions higher than three only  in special cases.

The poor scaling with dimension is crucial for fast direct solvers as
well. Rather quickly, researchers discovered a bevy of algorithms with very high
practical speed for problems in one and two dimensions.  In contrast, practical
algorithms for problems in three dimensions have only slowly started to emerge
over the past several years.  What makes the three dimensional case tractable at
all is that we have at our disposal two powerful techniques that reduce the
\textit{effective} dimensionality of the dense operators that need to be
inverted from three to two dimensions:
\begin{enumerate}
\item For boundary value problems with standard elliptic operators (e.g.~Laplace
  or Helmholtz) having constant coefficients, the problem can often be
  reformulated mathematically as a boundary integral equation. The resulting
  dense operator is defined only on a 2D surface rather than the full 3D volume

\item For problems that involve variable coefficients, it is standard practice
  to discretize the PDE directly (since analytic knowledge of the relevant
  Green's function rarely exists), leading to a sparse linear system.  When
  factorizing the coefficient matrix using nested dissection or a multifrontal
  method, large dense Schur complement matrices arise.  However, it turns out
  that these Schur complements are defined on \emph{mesh separators} in the
  nested dissection ordering of the mesh; more importantly, these separators
  turn out to be, in essence, two-dimensional substructures within the original
  three-dimensional mesh.  (Figures~\ref{fig:nd2} and~\ref{fig:nd3} may clarify
  the situation.)  Furthermore, and somewhat surprisingly, these Schur
  complements behave very much like discretized integral operators acting on 2D
  boundaries.
\end{enumerate}

\subsection{Scope and purpose}

The purpose of the survey is to provide a roadmap of the still-maturing field of
fast direct solvers, and to weave together seemingly disparate strands of
research. In particular, we seek to highlight the common intellectual ground
between research on FDS for integral equations, FDS for sparse systems arising
from the FEM/FD discretization of elliptic PDEs, and recently developed fast
direct algorithms for kernel matrices in machine learning and
computational statistics.

Our focus is on high level mathematical concepts such as dissipation of
information, different approaches to matrix factorizations, and different ways
to numerically represent Calder\'on-Zygmund operators.  Extensive pointers to the
existing literature are provided for the reader who wishes
to learn more details and to put the methods into practical use.  The text
includes a discussion of both coercive elliptic equations with non-oscillatory
kernels such as Laplace and Stokes, and of equations with oscillatory
kernels such as the Helmholtz and time-harmonic Maxwell equations.

\begin{remark}[High frequency problems]
\label{remark:highfreq}
  In the analysis of direct solvers for problems with oscillatory operators, a
  number of different asymptotic regimes have been considered in the literature.
  The most straight-forward is one where the PDE is kept fixed as the number of
  degrees of freedom~$N$ is increased; in this regime, the methods described in
  this survey often attain close to linear complexity.  A more challenging
  regime is one where the underlying frequency is also increased along with the
  number of degrees of freedom (in order to keep the number of discretization
  points per wavelength fixed).  While promising techniques for attaining
  quasi-linear complexity in this regime have been
  proposed~\cite{2016_michielssen_butterfly}, these techniques remain far less
  well-understood at the time of writing and do not form part of the scope of
  the current survey.  The main roadblock to efficient FDS in this regime is
  that when using similar tessellations of the system matrix as in the
  non-oscillatory case, blocks that were previously low rank cease to be low
  rank anymore; their ranks, in fact, grow with the underlying frequency, and
  are therefore proportional to some power of $N$.
\end{remark}

\subsection{Other surveys and texts}

Several texts and surveys describing fast direct solvers and rank structured
matrices already exist.
In particular, the 2016 survey \cite{2016_kressner_rankstructured_review} addresses similar
questions and also emphasizes high
level concepts, but it approaches the problem from a linear algebraic viewpoint,
rather than from a PDE or applied-analysis viewpoint.

The 2008 and 2010 monographs \cite{2008_bebendorf_book,2010_borm_book} describe
$\mathcal{H}$ and $\mathcal{H}^{2}$ matrices, respectively.
These texts were followed in 2015 by the comprehensive treatment 
\cite{hackbusch2015hierarchical} that covers both the
$\mathcal{H}$ and the $\mathcal{H}^{2}$ formalisms, and also describes
extensions to the tensor case.
The monographs \cite{2008_bebendorf_book,2010_borm_book,hackbusch2015hierarchical} 
are similar to the current text in that they discuss both linear algebraic and PDE
considerations, but they are more technical, and do not reflect
developments in recent years.  The 2009 survey \cite{2009_martinsson_ACTA} is
more limited in scope than the current text, as it restricts attention to
integral equations. 
The 2020 survey \cite{2020_keyes_hierarchical_hierarchical} hones in on
hierarchical algorithms that exploit rank structures broadly. The focus is on
how techniques of this type are well matched to modern hardware, and parallelize
well.  The review~\cite{2016_acta_sparse_direct_survey} summarizes the state of
research on sparse direct solvers.

The 2019 textbook \cite{2019_martinsson_fast_direct_solvers} covers similar
topics to the current text. However, the purpose
of~\cite{2019_martinsson_fast_direct_solvers} was to describe basic techniques
in detail, primarily oriented towards students and end users.  In contrast, our
objective here is to provide an introduction to FDS for applied mathematicians,
with a focus on mathematics and physics.

\subsection{Outline}
The survey is organized as follows.  Sections~\ref{sec:onedimbvp}
and~\ref{sec:onedimbvpIE} introduce key concepts through two simple model
problems involving one dimensional domains.  Section~\ref{sec:HODLR} provides a
first introduction to rank structured matrices by describing perhaps the
simplest hierarchical format (HODLR) that enables near linear complexity
algorithms for matrix-vector multiplication and for matrix inversion.
Section~\ref{sec:FDSPDE} describes how fast direct solvers can be used in
environments where a PDE is discretized directly to result in a sparse
coefficient matrix.  Sections~\ref{sec:BIE2D}~--~\ref{sec:BIE_strong} introduce
some basic ideas around integral equation formulations of classical boundary
value problems, and describe how rank structured matrix algebra can be used to
directly solve the resulting linear systems with linear or close to linear
complexity.  In Section~\ref{sec:operator_algebra}, we give an overview of how
fast direct solvers enable the explicit evaluation of functions of matrices,
which unlocks new operator preconditioning strategies and methods for solving
coupled multiphysics problems.  Section~\ref{sec:dissipation} provides numerical
results that demonstrate the remarkable persistence of numerical rank
deficiencies in dense matrices that arise from a range of different contexts in
scientific simulations.  Section~\ref{sec:compression} describes some popular
techniques for obtaining data-sparse representations of rank structured
matrices, and subsequently Section~\ref{sec:statistics} broadens the scope and
describes how many of the tools that were introduced earlier in the survey can
be used to accelerate key tasks in data science and computational statistics.
Section~\ref{sec:taxonomy} summarizes the material introduced and provides more
detailed guidance for what techniques are best suited for particular problems,
including pointers to existing software packages.

A reader seeking a quicker and more casual read can pick up the core ideas by
focusing on
Sections~\ref{sec:onedimbvp},~\ref{sec:onedimbvpIE},~\ref{sec:HODLR},~\ref{sec:FDSPDE},
and~\ref{sec:dissipation}. To learn more about connections to problems in
computational statistics, the reader could peruse Section~\ref{sec:statistics}.


%% file: 02-onedimbvp/onedimbvp02.tex
In Section \ref{sec:intro}, we claimed that when we approximate a solution
operator to an elliptic PDE the resulting matrix is dense, but data-sparse.  We
will now illustrate what these terms mean in the particularly simple case of an
elliptic boundary value problem in one dimension:
\begin{equation}
\label{eq:twopointbvp}
\begin{aligned}
-u''(x) + m(x)\,u(x) =&\ g(x),\qquad&x \in I,\\
u(a) =&\ f(a),\\
u(b) =&\ f(b),
\end{aligned}
\end{equation}
where $I = (a,b)$ is an interval on the line,  $m$~is a given function on
$I$, and~$f(a)$ and~$f(b)$ specify the Dirichlet boundary conditions.

\begin{figure}
\centering
\begin{tabular}{ccccccc}
$\mtx{A}$ 
&\mbox{}\hspace{5mm}\mbox{}&
$\mtx{A}^{-1}$ 
&\mbox{}\hspace{5mm}\mbox{}&
$\mtx{L}$ 
&\mbox{}\hspace{5mm}\mbox{}&
$\mtx{U}$
\\
\includegraphics[width=14mm]{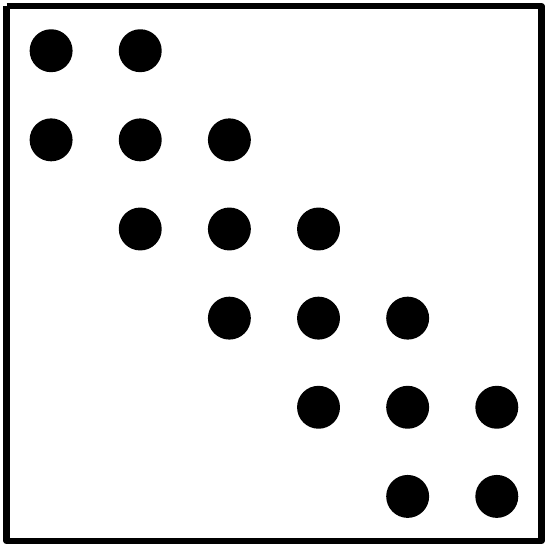} &&
\includegraphics[width=14mm]{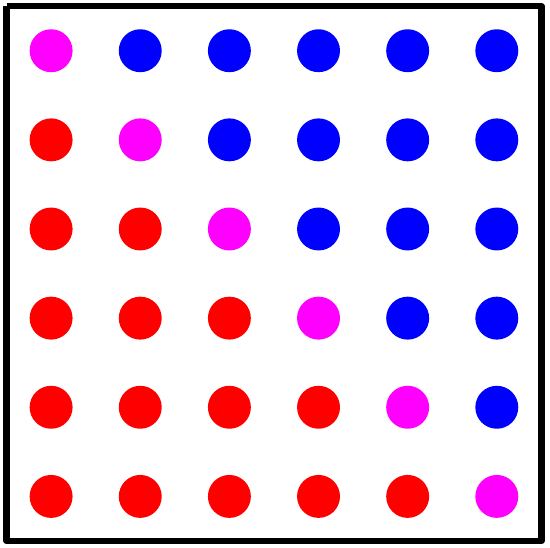} &&
\includegraphics[width=14mm]{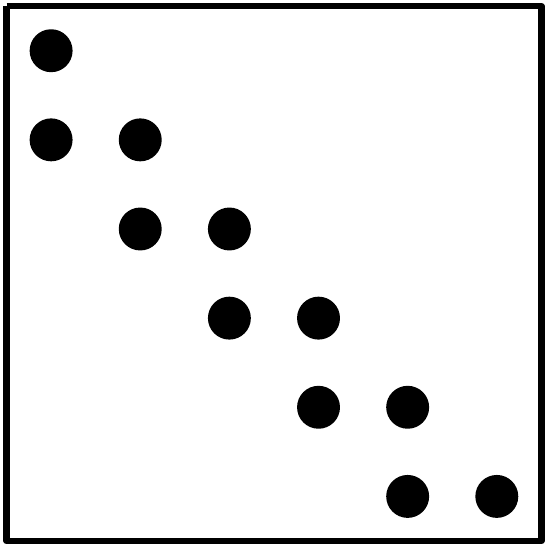} &&
\includegraphics[width=14mm]{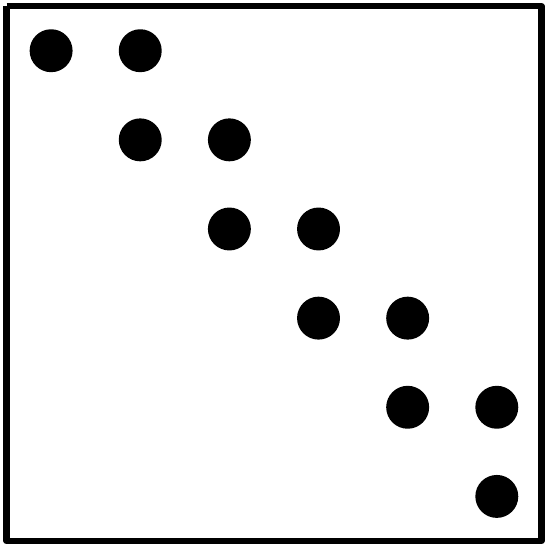}
\\
(a) && (b) && (c) && (d)
\end{tabular}
\caption{
(a) Sparsity pattern of the tridiagonal coefficient matrix $\mtx{A}$ introduced in Section \ref{sec:onedimbvp} for $N=6$.
(b) The inverse $\mtx{A}^{-1}$ is \textit{semi-separable}, which is to say that its upper triangular part (the blue elements) and its lower triangular part (the red elements) are each restrictions of rank-1 matrices. The two rank-one matrices must agree on the diagonal (purple), so a semi-separable matrix is defined by $3N-2$ real numbers, just like a tridiagonal one.
(c,d) The sparsity pattern of the factors $\mtx{L}$ and $\mtx{U}$ in the  factorization $\mtx{A} = \mtx{L}\mtx{U}$.
Since $\mtx{L}(i,i)=1$, there are again $3N-2$ independent numbers.
}
\label{fig:tridiag}
\end{figure}

To keep things simple, we discretize~\eqref{eq:twopointbvp} using
a second-order accurate finite difference method on a set
of uniformly spaced points. With~$N$ a positive integer that determines the
resolution of the mesh, let $\{x_{i}\}_{i=0}^{N+1}$ denote these points, where
$$
x_{i} = a + ih,
\qquad\mbox{and where}\qquad
h = \frac{b-a}{N+1}.
$$
Replacing the second derivative in (\ref{eq:twopointbvp}) by
the second-order accurate approximation
\begin{equation}
  -u''(x_{i}) \approx h^{-2}\bigl(-u(x_{i-1}) + 2u(x_{i}) - u(x_{i+1})\bigr),
\end{equation}
we obtain a linear system
\begin{equation}
\label{eq:finitediff}
h^{-2}\bigl(-u_{i-1} + 2u_{i} - u_{i+1}\bigr) + m(x_{i})\,u_{i} = g(x_{i}),
\qquad i \in \{1,\,2,\,\dots,\,N\}.
\end{equation}
We will denote by $\vct{u} = (u_1 \cdots u_N)^T$ the vector
of approximations to the solution, so that $u_{i} \approx u(x_{i})$.
At the endpoints, we use the specified boundary values to set
\begin{equation}
  \label{eq:endpoints}
u_{0} = f(a)
\qquad\mbox{and}\qquad
u_{N+1} = f(b).
\end{equation}
The linear system~(\ref{eq:finitediff}) can then
conveniently be written in matrix form as
\begin{equation}
\label{eq:tridiagsys}
\mtx{A}\vct{u} = \vct{g},
\end{equation}
where~$\mtx{A}$ is an~$N\times N$ matrix and $\vct{g}$ is a vector in
$\mathbb{R}^{N}$ that encodes both the body load $g$ and the given boundary data
$f$. The matrix $\mtx{A}$ is tridiagonal, with entries defined
by~\eqref{eq:finitediff}, cf.~Figure~\ref{fig:tridiag}(a).

Subsequently, the linear system~(\ref{eq:tridiagsys}) can be solved through
Gaussian elimination using only $\mathcal \cO(N)$ operations since each
elimination step involves only three variables. Another way to phrase this fact
is to say that the LU factorization of a tridiagonal matrix consists of
two bidiagonal matrices, as illustrated in Figures~\ref{fig:tridiag}(c) and (d).

On the other hand, as an alternative to computing an LU factorization,
it is possible in the present case to directly compute the inverse
$$
\mtx{B} = \mtx{A}^{-1}.
$$ The matrix $\mtx{B}$ is of course dense, but it turns out that it
is data sparse and can be explicitly computed, stored, and applied to vectors
using only $\mathcal \cO(N)$ operations. Specifically, $\mtx{B}$ is
\textit{semi-separable}, which is to say that it consists of two rank-one
matrices that are spliced together along the diagonal. To be precise, there
exist vectors $\vct{c},\, \vct{d},\, \vct{e},\, \vct{f} \in \mathbb{R}^{N}$,
for which
\begin{equation}
\label{eq:semisep}
\mtx{B}(i,j) =
\left\{\begin{aligned}
c_{i}d_{j},\qquad i \geq j \qquad\mbox{(on or below the diagonal)},\\
e_{i}f_{j},\qquad i \leq j \qquad\mbox{(on or above the diagonal)}.
\end{aligned}\right.
\end{equation}
We see that the subdiagonal half of $\mtx{B}$ is the restriction of
the rank-1 matrix $\vct{c}\vct{d}^{*}$ and the superdiagonal half is
the restriction of $\vct{e}\vct{f}^{*}$. On the diagonal, the two
halves must agree, so the vectors satisfy the constraints
$$
c_{i}d_{i} = e_{i}f_{i},\qquad i \in \{1,\,2,\,\dots,\,n\}.
$$ Taking into account the ambiguity due to scaling (in the sense that
if $t\neq 0$, then of course \mbox{$(t\vct{c})(t^{-1}\vct{d}^{*}) =
\vct{c}\vct{d}^{*}$}), we see that precisely $3N-2$ real numbers are
required to define a semi-separable matrix.
(Not coincidentally, there are also $3N-2$ entries in a tridiagonal matrix. In
fact, there is a one-to-one correspondence between invertible tridiagonal
matrices and invertible semi-separable matrices \cite{2008_vandebril_semiseparable_book}.)

To summarize, the key finding of this section is that the
matrix $\mtx{A}$ that encodes the discretized differentiation
operator is {sparse}, and that its inverse is dense but
{data-sparse}. Furthermore, they can both be applied to vectors using $\cO(N)$
operations.


%% file: 03-onedimbvpIE/onedimbvpIE03.tex
Let us next formulate a solver for our one dimensional model problem
(\ref{eq:twopointbvp}) that is based on an integral equation formulation.  This
formulation leads to a linear system with a dense coefficient matrix, but we
will see that it has an internal structure that also allows for constructing a
linear complexity solver.
From a mathematical point of view, the reformulation in integral-form leads to a
well-conditioned linear system which greatly improves the numerical stability of
the solver.

\subsection{Integral equation formulation}
In order to derive an integral equation formulation of (\ref{eq:twopointbvp}),
we recall~\cite{coddington1984,lee1997stiff} that in the special case where the
coefficient function $m$ is zero, and we have homogeneous boundary data, we can
write down an analytic solution explicitly.
To be precise, the equation
\begin{equation}
\label{eq:specialtwopointbvp}
\begin{aligned}
-u''(x) =&\ g(x),\qquad&x \in I,\\
u(a) =&\ 0,\\
u(b) =&\ 0,
\end{aligned}
\end{equation}
has the solution
$$
u(x) = \int_{a}^{b} G(x,y) \, g(y) \,dy,\qquad x \in I,
$$
where the \textit{Green's function} $G$ is defined via
\begin{equation}
\label{eq:green_function_def}
G(x,y) =
\begin{aligned}
\frac{(b-x)(y-a)}{b-a},\qquad&\mbox{when}\ x \geq y\qquad\mbox{(on or below the diagonal)},\\
\frac{(x-a)(b-y)}{b-a},\qquad&\mbox{when}\ x \leq y\qquad\mbox{(on or above the diagonal)}.
\end{aligned}
\end{equation}
(The structural similarity between (\ref{eq:green_function_def}) and
(\ref{eq:semisep}) is of course no coincidence.)
Now consider the more general equation
\begin{equation}
\label{eq:homtwopointbvp}
\left\{\begin{aligned}
-u''(y) + m(y) u(y) = &\ g(y),\qquad&y \in I,\\
u(a) =&\ 0,\\
u(b) =&\ 0.
\end{aligned}\right.
\end{equation}
Multiplying the differential equation in~\eqref{eq:homtwopointbvp} by~$G(x,y)$,
and then integrating in $y$ from $a$ to $b$ results in the integral equation
\begin{equation}
  \label{eq:onedimLS}
u(x) + \int_{a}^{b}G(x,y) \, m(y) \, u(y) \,dy = 
\int_{a}^{b} G(x,y) \, g(y) \,dy,
\qquad x \in I.
\end{equation}
We discretize (\ref{eq:onedimLS}) using a simple Nystr\"om method that
relies on the same grid $\{x_{i}\}_{i=0}^{N+1}$ that we introduced in
Section \ref{sec:onedimbvp}. Collocating equation (\ref{eq:onedimLS})
to the interior discretization nodes, and replacing the integral
by its approximation by the trapezoidal rule, we obtain the linear system
\begin{equation}
\label{eq:onedimLSdisc}
u_{i} + \sum_{j=1}^{N}h \, G(x_{i},x_{j}) \, m(x_{j}) \, u_{j} = 
\sum_{j=1}^{N}h \, G(x_{i},x_{j}) \, g(x_{j}),
\qquad i \in \{1,\,2,\,\dots,\,n\}.
\end{equation}
Observe that when we applied the trapezoidal rule, the boundary terms with the
weights~$h/2$ vanish since $u(a) = u(b) = 0$.  To be concise, we
write~\eqref{eq:onedimLSdisc} in matrix form as
\begin{equation}
\label{eq:kim}
\bigl(\mtx{I} + \mtx{G}\mtx{M}\bigr)\vct{u} = \mtx{G}\vct{f},
\end{equation}
where the $N\times N$ matrix $\mtx{G}$ has entries
\begin{equation}
\label{eq:jong-il}
\mtx{G}(i,j) = h\,G(x_{i},x_{j}),
\end{equation}
and where $\mtx{M}$ is the diagonal matrix whose diagonal entries are 
$\{m(x_{i})\}_{i=1}^{N}$.

\begin{remark}[Inhomogeneous boundary data]
\label{remark:inhom}
The technique we described in this section for solving two-point
boundary value problems with zero boundary conditions can easily
be extended to general Dirichlet data.
Specifically, if we seek to solve the more general
problem~(\ref{eq:twopointbvp}), we can look for a solution of the form $u
= v + w$, where
\[
w(x) = f(a)\frac{b-x}{b-a} + f(b)\frac{x-a}{b-x}
\]
is a linear function that satisfies the boundary conditions (the ``homogeneous solution''). Since $-w'' = 0$, the function $v$ (the ``particular solution'') satisfies
\begin{equation}
\label{eq:hombvp2}
-v''(x) + m(x)v(x) = f(x) -
m(x)w(x)
\end{equation}
with zero boundary data. Equation~\eqref{eq:hombvp2} can then be solved using
the techniques described in this section.
\end{remark}

\subsection{Solving the integral equation efficiently}

At first, it may appear that forming the linear system (\ref{eq:kim}) would have
cost $\cO(N^2)$, since all matrices involved are dense, and then the cost to solve
it would be $\cO(N^3)$. However, the coefficient matrix turns out to be
data-sparse, which enables us to solve the equation in linear time. The key
observation is that relation (\ref{eq:green_function_def}) implies that the
subdiagonal and the superdiagonal halves of the matrix
$\mtx{I} + \mtx{G}\mtx{M}$ are each restrictions of matrices of rank~1, so the
matrix as a whole has the structure ``diagonal + semi-separable''.  It is known
that the inverse of such a matrix is also of the form ``diagonal +
semi-separable'', and can be constructed using $\cO(N)$ operations. We will not
describe specialized algorithms for this task, but will instead describe a
more general machinery that can handle the matrix $\mtx{I} + \mtx{G}\mtx{M}$,
and also a range of other dense matrices that arise when discretizing PDEs and
integral equations. Section \ref{sec:HODLR} describes a basic method with
complexity $\cO(N \log^2 N)$, and then Section~\ref{sec:BIE_nested} 
describes a scheme with~$\cO(N)$ complexity.

\subsection{Conditioning}
One motivation for introducing the integral equation formulation of the BVP was
to provide another illustration of the notion of a dense but data-sparse
matrix. A more important reason was to illustrate how the conversion to an
integral equation often provides an exact mathematical technique for vastly
improving the numerical stability of the discretized system.

Considering first the finite difference formulation of Section
\ref{sec:onedimbvp}, we observe that as $N$ grows, the matrix~$\mtx{A}$ should
provide a successively more faithful approximation to the differentiation
operator $u \mapsto -u'' + m\,u$. Since this operator is \textit{unbounded}, 
the largest eigenvalue of $\mtx{A}$ must grow as $N$
grows. In fact, one can easily demonstrate that the condition number
$\kappa(\mtx{A})$ of $\mtx{A}$ satisfies $\kappa(\mtx{A}) \sim N^{2}$ as
$N\rightarrow \infty$.

\begin{figure}[t]
  \centering
  \begin{subfigure}{.45\linewidth}
    \centering
    \includegraphics[width=.95\linewidth]{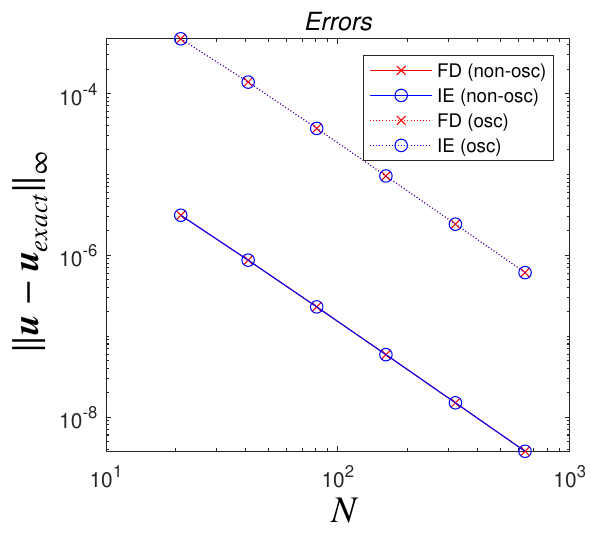}
    \caption{A comparison of errors in solving a two-point boundary value
      problem when using a finite difference versus integral equation formulation.}
  \end{subfigure}
  \begin{subfigure}{.45\linewidth}
    \centering
    \includegraphics[width=.95\linewidth]{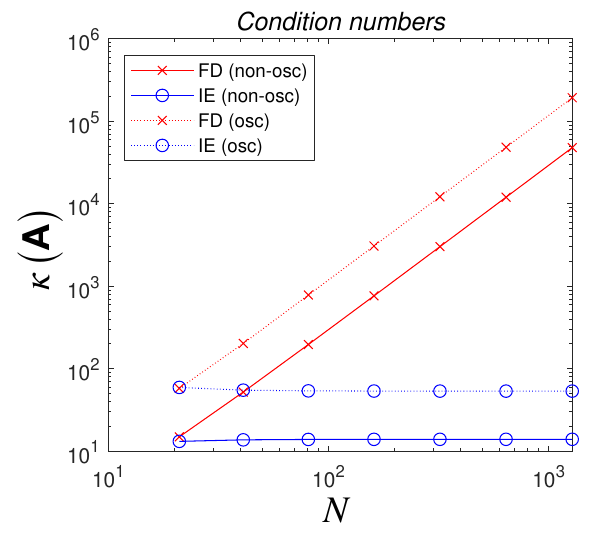}
    \caption{Condition numbers for the matrices associated with
      solving a two-point boundary value
      problem when using a finite difference versus integral equation formulation.}
  \end{subfigure}
  \caption{A numerical example illustrating the behavior of the finite
      difference (FD) method of Section~\ref{sec:onedimbvp} and the integral
      equation (IE) method of Section~\ref{sec:onedimbvpIE}.
      We first considered
      a problem on the domain $[0,1]$ with a non-oscillatory solution
      ``(non-osc)'' where $m(x) = 100(1+x)\cos(x)$ and $g(x) = 1 +
      \cos(1+x)$. We then swapped the sign of $m$ (but kept everything else the
      same) to get a problem with an oscillatory solution ``(osc)''.
      A key point here is
      that the condition number of the integral equation formulation does not
      grow with~$N$. A secondary point is that elliptic problems with
      oscillatory solutions are far more challenging.
      }
    \label{fig:FD_vs_IE}
\end{figure}

In contrast, the coefficient matrix $\mtx{I} + \mtx{G}\mtx{M}$ in~\eqref{eq:kim}
will, as $N$ grows, become a successively more accurate approximation to a
continuum operator of the form ``identity + compact'' -- or, in other words, a
Fredholm operator of the second kind~\cite{1990_riesz}. Such an operator is very
well-behaved mathematically; in particular, upon its discretization, it
leads to a linear system that is typically well-conditioned. A full discussion
of Fredholm theory is beyond the scope of this survey, but a classic treatment
in the integral operator context is~\cite{1990_riesz}. For our purposes, let us
simply point to Figure~\ref{fig:FD_vs_IE} that, through a numerical example,
illustrates how the condition number of $\mtx{I} + \mtx{G}\mtx{M}$ is basically
constant (and
small).

\begin{remark}[A curious connection]
  In the numerical example shown in Figure~\ref{fig:FD_vs_IE}, we plot not only
  the condition numbers of the coefficient matrices but also the errors in the
  computed approximate solution for the two methods. We see that they are
  basically identical, which is the result of the fact that the matrix~$\mtx{G}$
  defined by~(\ref{eq:jong-il}) happens to be the \textit{exact} inverse of the
  one-dimensional finite difference stencil. In other words, if we write the
  matrix $\mtx{A}$ introduced in Section~\ref{sec:onedimbvp} as
  $\mtx{D} + \mtx{M}$, then $\mtx{G} = \mtx{D}^{-1}$. It follows that
$$
\bigl(\mtx{D} + \mtx{M}\bigr)^{-1} = 
\bigl(\mtx{D}(\mtx{I} + \mtx{D}^{-1}\mtx{M})\bigr)^{-1} = 
\bigl(\mtx{I} + \mtx{D}^{-1}\mtx{M}\bigr)^{-1}\mtx{D}^{-1} = 
\bigl(\mtx{I} + \mtx{G}\mtx{M}\bigr)^{-1}\mtx{G}.
$$
In other words, the discretized systems we derived in Sections
\ref{sec:onedimbvp} and \ref{sec:onedimbvpIE} are mathematically equivalent (up
to errors from floating point arithmetic).
\end{remark}


%% file: 04-hodlr/HODLRnew04.tex
In this section, we revisit a class of rank structured matrices that we
introduced in Figure \ref{fig:tessellation}, and is related to the
semi-separable format we saw in Sections~\ref{sec:onedimbvp}
and~\ref{sec:onedimbvpIE}. This \textit{hierarchically off-diagonal low rank}
(HODLR) format allows close to linear complexity algorithms for solving a system
such as the discretized integral equation (\ref{eq:kim}).  The coefficient
matrix in~(\ref{eq:kim}) has the very special structure ``semi-separable +
diagonal''. Specialized fast algorithms exist for this task, but we choose to
work with the HODLR format as it can handle a wide range of problems and
introduces several key ideas that will play an important role when solving
problems in dimensions two and three.

\subsection{Informal definition and basic operations}

To mathematically describe the tessellation implied by
Figure~\ref{fig:tessellation}, we decompose a given $N\times N$ matrix
$\mtx{A}$ as a sum
\begin{equation}
\label{eq:HODLR_decomp}
\mtx{A} = \mtx{A}^{(1)} + \mtx{A}^{(2)} + \cdots + \mtx{A}^{(L)} + \mtx{D},
\end{equation}
where $\mtx{A}^{(1)}$ holds the two largest off-diagonal blocks of~$\mtx{A}$,~$\mtx{A}^{(2)}$ holds to four next largest blocks, and so on until the remaining
(small) diagonal blocks are left in the matrix $\mtx{D}$. In other words, for
the matrix shown in Figure \ref{fig:tessellation}, we write
$$
\begin{array}{ccccccccc}
\mtx{A} &=& \mtx{A}^{(1)} &+& \mtx{A}^{(2)} &+&\mtx{A}^{(3)} &+& \mtx{D}\\
\includegraphics[width=14mm]{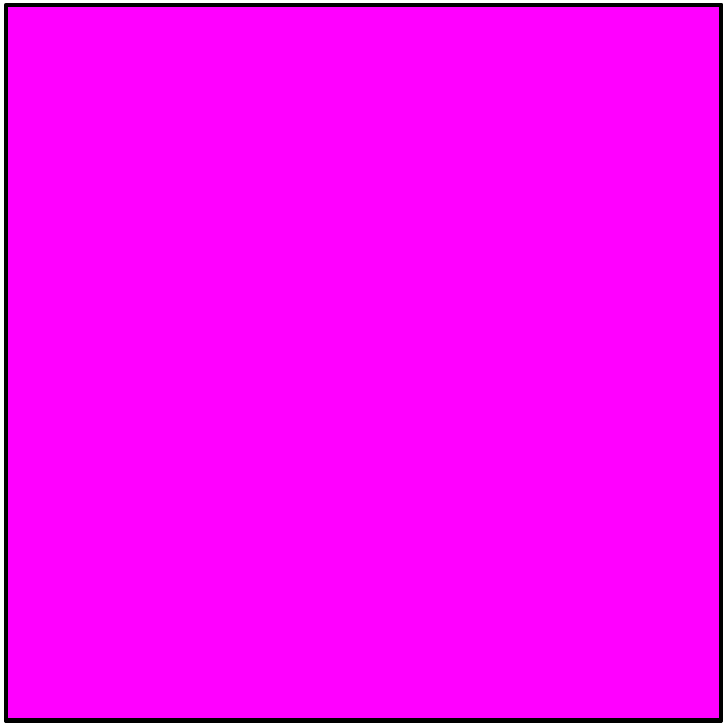}
&&
\includegraphics[width=14mm]{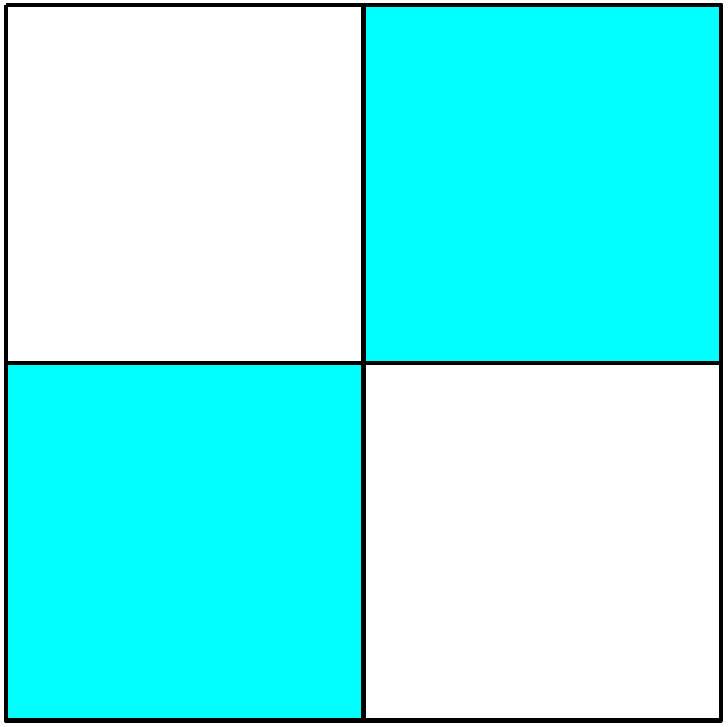}
&&
\includegraphics[width=14mm]{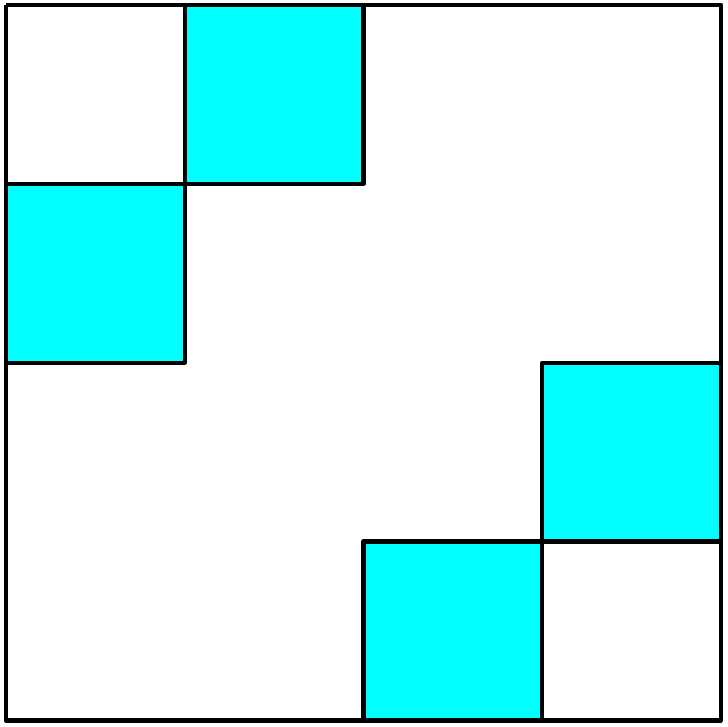}
&&
\includegraphics[width=14mm]{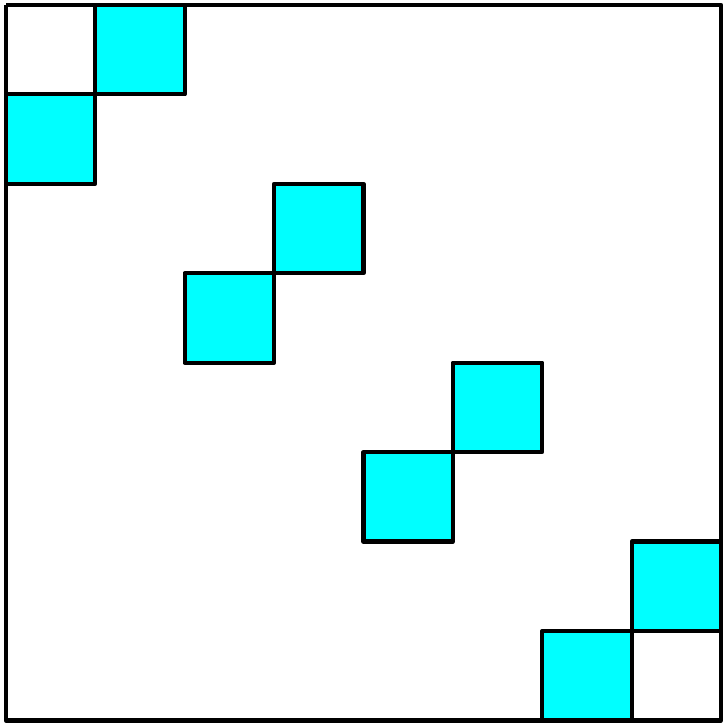}
&&
\includegraphics[width=14mm]{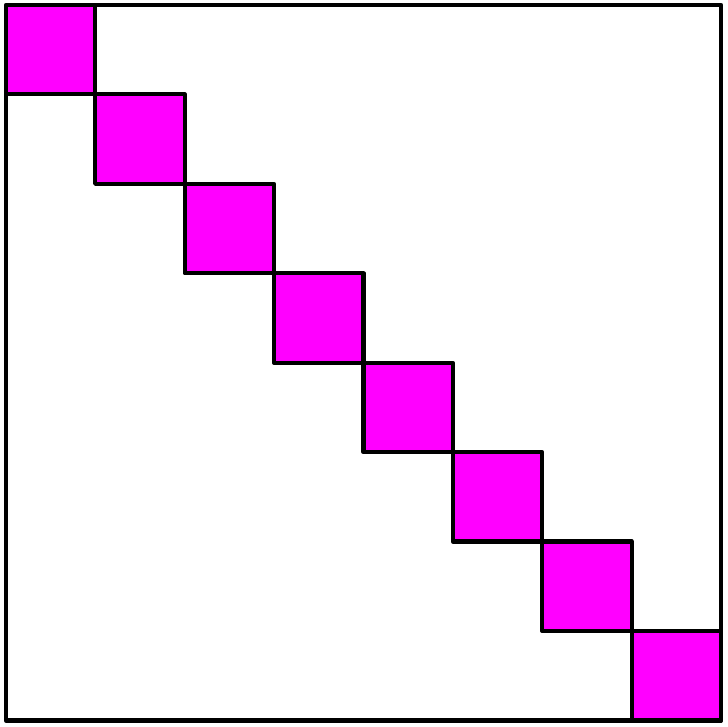}
\end{array} .
$$
We now say that $\mtx{A}$ is a HODLR matrix of rank $k$ if all blocks of the
matrices $\mtx{A}^{(\ell)}$ have rank $k$ or less and if the diagonal blocks of
$\mtx{D}$ are roughly of size $2k\times 2k$ or less. Observe that the last condition
implies that the number of levels $L$ will satisfy
$$
L \approx \log_{2}(N/(2k)).
$$

The storage cost of a HODLR matrix can easily be determined by observing that
each matrix $\mtx{A}^{(\ell)}$ requires $\sim Nk$ storage (since
it holds $2^{\ell}$ blocks, each of size $N\,2^{-\ell} \times N\,2^{-\ell}$ and
of rank at most $k$). Consequently, the cost to store all matrices
$\{\mtx{A}^{(\ell)}\}_{\ell=1}^{L}$ is $\cO(NLk) = \cO(N\log(N) k)$. Turning
next to the matrix $\mtx{D}$, we see that it holds about $N/(2k)$ blocks, each
of size at most $2k \times 2k$, for a total of $\sim Nk$ storage as well. In all, we
see that $\mtx{A}$ requires $\cO(N\log(N) k)$ storage.

Analyzing the cost of executing a matrix-vector multiplication for a HODLR
matrix is entirely analogous to analyzing the storage cost. We simply observe
that each matrix in the sum (\ref{eq:HODLR_decomp}) can be applied to a vector
using $\cO(Nk)$ operations, and since there are $L \sim \log(N)$ terms
in the sum, the cost of the matrix-vector multiplication becomes $\cO(N\log(N)k)$.

\subsection{A recursive inversion algorithm}
\label{sec:HODLRwoodbury}

A conceptually simple algorithm for inverting a HODLR matrix is obtained by
recursively applying a version of the well-known Woodbury formula. To be
precise, we will use the fact that if~$\mtx{A}$ is an~$N\times N$ matrix that
admits a split
$$
\begin{array}{cccccccc}
\mtx{A} &=& \mtx{D} & + & \mtx{W} & \mtx{V}^{*},\\
N\times N && N\times N && N\times p & p\times N
\end{array}
$$
where $\mtx{D}$ is in some sense ``easy to invert'' and $\mtx{W}\mtx{V}^{*}$ is
a low rank term, then
\begin{equation}
\label{eq:woodburymini}
\mtx{A}^{-1} = \mtx{D}^{-1} - \mtx{D}^{-1}\mtx{W}\mtx{S}^{-1}\mtx{V}^{*}\mtx{D}^{-1},
\end{equation}
where
\begin{equation}
\label{eq:formulaS}
\mtx{S} =
\mtx{I} + \mtx{V}^{*}\mtx{D}^{-1}\mtx{W}
\end{equation}
is a $p\times p$ (i.e.~small) matrix that needs to be invertible for
(\ref{eq:woodburymini}) to hold.

To apply the Woodbury formula to a HODLR matrix $\mtx{A}$, let us start by
partitioning it into $2\times 2$ blocks
\begin{equation}
\label{eq:huey}
\mtx{A} = \left[\begin{array}{rr}
\mtx{A}_{\alpha,\alpha} & \mtx{A}_{\alpha,\beta} \\
\mtx{A}_{\beta ,\alpha} & \mtx{A}_{\beta ,\beta}
\end{array}\right].
\end{equation}
Now observe that the diagonal blocks $\mtx{A}_{\alpha,\alpha}$ and
$\mtx{A}_{\beta,\beta}$ are themselves HODLR matrices; furthermore, the
off-diagonal blocks are low rank matrices that can be written as 
\begin{equation}
\label{eq:dewey}
\begin{array}{ccccccc}
\mtx{A}_{\alpha,\beta} &=& \mtx{W}_{\alpha,\beta}&\mtx{V}_{\beta}^{*},\\
\tfrac{N}{2} \times \tfrac{N}{2} && \tfrac{N}{2} \times k & k \times \tfrac{N}{2}
\end{array}
\qquad\mbox{and}\qquad
\begin{array}{ccccccc}
\mtx{A}_{\beta,\alpha} &=& \mtx{W}_{\beta,\alpha}&\mtx{V}_{\alpha}^{*}.\\
\tfrac{N}{2} \times \tfrac{N}{2} && \tfrac{N}{2} \times k & k \times \tfrac{N}{2}
\end{array}
\end{equation}
Combining (\ref{eq:huey}) and (\ref{eq:dewey}), we write
$$
\mtx{A} =
\underbrace{\left[\begin{array}{rr}\mtx{A}_{\alpha,\alpha} & \mtx{0} \\ \mtx{0} & \mtx{A}_{\beta,\beta}\end{array}\right]}_{=:\mtx{D}}
+
\underbrace{\left[\begin{array}{rr}\mtx{0} & \mtx{W}_{\alpha,\beta} \\ \mtx{W}_{\beta,\alpha} & \mtx{0}\end{array}\right]}_{=:\tilde{\mtx{W}}}
\underbrace{\left[\begin{array}{rr}\mtx{V}_{\alpha}^{*} & \mtx{0} \\ \mtx{0} & \mtx{V}^{*}_{\beta}\end{array}\right].}_{=:\mtx{V}^{*}}
$$
The Woodbury formula (\ref{eq:woodburymini}) is now immediately
applicable:~$\mtx{D}^{-1}$ can be computed by recursing
on~$\mtx{A}_{\alpha,\alpha}$ and~$\mtx{A}_{\beta,\beta}$, and~$\mtx{S}$
in~\eqref{eq:formulaS} takes the form
\begin{equation}
\mtx{S} =
\begin{bmatrix}
\mtx{I} & \mtx{V}_{\alpha}^{*}\mtx{A}_{\alpha,\alpha}^{-1}\mtx{W}_{\alpha,\beta} \\
\mtx{V}_{\beta}^{*}\mtx{A}_{\beta,\beta}^{-1}\mtx{W}_{\beta,\alpha} & \mtx{I}
\end{bmatrix}.
\end{equation}
Since $\mtx{S}$ is of size $2k\times 2k$, its inversion is inexpensive.

The algorithm described in this section works and is sometimes
efficient. It has a key drawback, however, which is that the ranks of the HODLR
submatrices involved \textit{grow} during the course of its execution. The issue
is due to the second term in formula (\ref{eq:woodburymini}), which adds a
low-rank component to all the off-diagonal blocks of $\mtx{D}^{-1}$. For
matrices involving a small number of levels, and a low HODLR rank, the rank
growth may be tolerable, but in many situations it slows the algorithm down
considerably.

The issue of rank growth is not due to a flaw in the algorithm, since it turns
out that the inverse of a HODLR matrix of rank $k$ is \textit{not} a HODLR
matrix of rank $k$, but instead a HODLR matrix of rank bounded by $(L-1)k$. To
circumvent this challenge, one must look for an algorithm that represents the
computed inverse using a different format; we will show one such method in
Section \ref{sec:HODLRinv}.

\begin{remark}[Recompression]
  An alternative strategy is to numerically recompress the off-diagonal blocks
  to keep their ranks constant.  This would lead to~$\cO(1)$ errors in the
  general case, but often works well in practice, since many of the matrices
  under consideration have additional structure beyond merely being HODLR.
  Intuitively, if the matrix arises from discretizing a differential equation,
  and we know in advance that its inverse should also be compressible, then it
  is perhaps not surprising that recompression works well.
\end{remark}

\subsection{Exact non-recursive inversion of a HODLR matrix}
\label{sec:HODLRinv}
We will next describe an algorithm for inverting a HODLR matrix that does not
lead to rank increases, and does not involve recursions. To achieve this, it
relies on a \textit{multiplicative} representation of the inverse, rather than
the additive one used Section \ref{sec:HODLRwoodbury}. The price to pay is that
the inverse computed will not be represented in the HODLR format.

To avoid the need to introduce a complicated notational apparatus, we describe
the algorithm for the basic 3-level HODLR matrix shown earlier in
Figure~\ref{fig:tessellation}, with the off-diagonal blocks labeled as shown in
Figure~\ref{fig:tree}. (Subsequently, Section~\ref{sec:trees} explains the
rationale for the indices of the off-diagonal blocks.) We will construct a
representation of $\mtx{A}^{-1}$ that takes the form
\begin{equation}
\label{eq:HOLDRinvfact}
\mtx{A}^{-1} = \mtx{B}_{0}\,\mtx{B}_{1}\,\mtx{B}_{2}\,\mtx{B}_{3},
\end{equation}
where each matrix $\mtx{B}_{\ell}$ is block diagonal and has diagonal blocks
that are rank-$k$ perturbations of the identity matrix. This means that each
$\mtx{B}_{\ell}$ can be applied to a vector in $\cO(Nk)$ flops.

The first step is to form a block diagonal matrix $\mtx{B}_{3}$ whose diagonal
blocks are the inverses of the diagonal blocks of $\mtx{A}$ (so that
$\mtx{B}_{3} = \mtx{D}^{-1}$ with $\mtx{D}$ as in
(\ref{eq:HODLR_decomp})). Applying $\mtx{B}_{3}$ to $\mtx{A}$, we form a matrix
$\mtx{A}_{3} = \mtx{B}_{3}\mtx{A}$ whose diagonal blocks are the identity
matrix:
$$
\begin{array}{cccccccc}
\mtx{A}_{3} &=& \mtx{B}_{3} & \mtx{A} \\
\includegraphics[width=37mm]{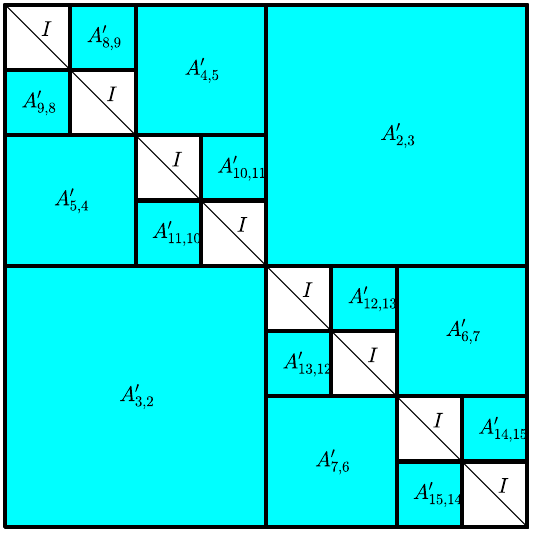}
&&
\includegraphics[width=37mm]{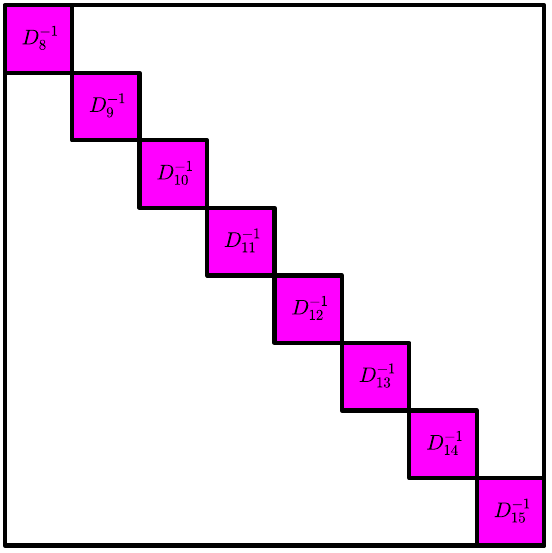}
&
\includegraphics[width=37mm]{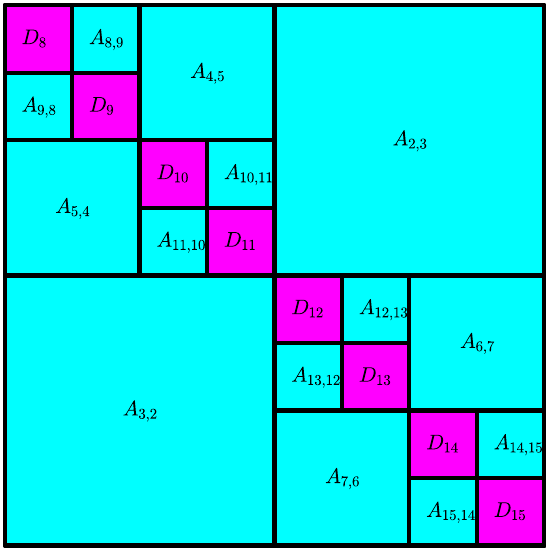}
\end{array}
$$
Observe that while the off-diagonal blocks in $\mtx{A}$ were all modified in
forming $\mtx{A}_{3}$, their \textit{ranks} remain bounded by $k$.

In the second step, let $\mtx{D}_{4},\,\mtx{D}_{5},\,\mtx{D}_{6},\,\mtx{D}_{7}$ denote
the diagonal blocks of $\mtx{A}_{3}$, so that
\begin{multline*}
\mtx{D}_{4} = 
\left[\begin{array}{rr}
\mtx{I} & \mtx{A}_{8,9}'\\
\mtx{A}_{9,8}' & \mtx{I} 
\end{array}\right],
\qquad
\mtx{D}_{5} = 
\left[\begin{array}{rr}
\mtx{I} & \mtx{A}_{10,11}'\\
\mtx{A}_{11,10}' & \mtx{I} 
\end{array}\right],
\\
\mtx{D}_{6} = 
\left[\begin{array}{rr}
\mtx{I} & \mtx{A}_{12,13}'\\
\mtx{A}_{13,12}' & \mtx{I} 
\end{array}\right],
\qquad
\mtx{D}_{7} = 
\left[\begin{array}{rr}
\mtx{I} & \mtx{A}_{14,15}'\\
\mtx{A}_{15,14}' & \mtx{I} 
\end{array}\right].
\end{multline*}
Since each of these matrices is of the form ``identity + low rank'', they can be
inexpensively inverted. We form a matrix $\mtx{B}_{2}$ by putting the resulting
inverses $\mtx{D}_{4}^{-1}$, $\mtx{D}_{5}^{-1}$, $\mtx{D}_{6}^{-1}$, and
$\mtx{D}_{7}^{-1}$ into its diagonal blocks. Then applying $\mtx{B}_{2}$ to
$\mtx{A}_{3}$, we introduce larger identity matrices into the diagonal blocks:
$$
\begin{array}{cccccccc}
\mtx{A}_{2} &=& \mtx{B}_{2} & \mtx{A}_{3}, \\
\includegraphics[width=36mm]{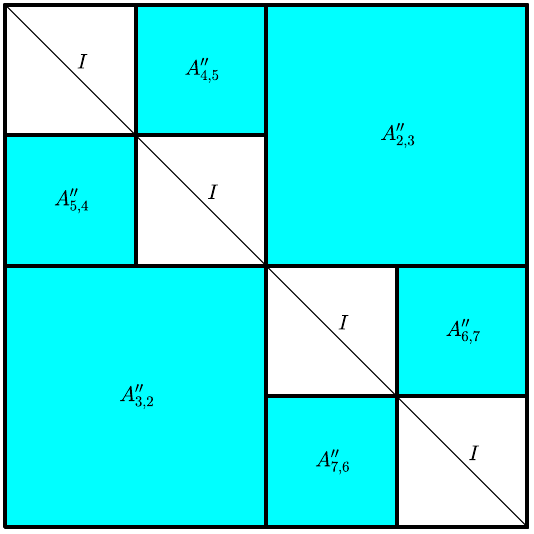}
&&
\includegraphics[width=36mm]{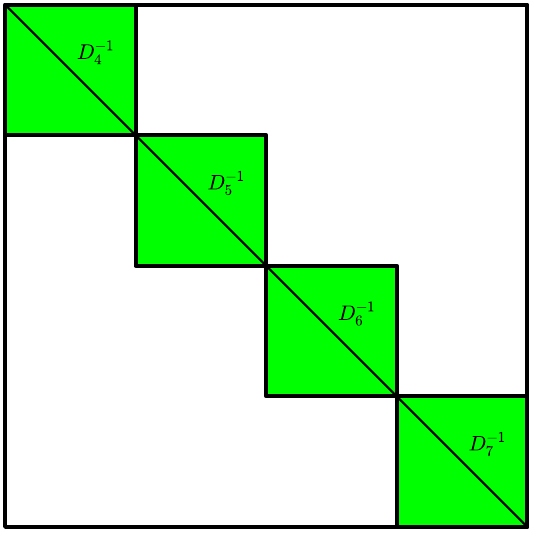}
&
\includegraphics[width=36mm]{04-hodlr/figs/fig_HODLR_inv_A3.pdf}
\end{array}
$$
Here, green blocks indicate ``identity + low rank''.
Observe again that the off-diagonal blocks of $\mtx{A}_{2}$ all have rank at most $k$.

The pattern is that we create successively larger and larger identity matrices along the diagonal, while maintaining the rank structure in the remainder of the matrix. The last two steps follow the same pattern, and lead to the final factorization
$$
\begin{array}{ccccccccccccccccc}
\mtx{I} &=&
\mtx{B}_{0}&\mtx{B}_{1}&\mtx{B}_{2}&\mtx{B}_{3}&\mtx{A},\\
&&
\includegraphics[width=20mm]{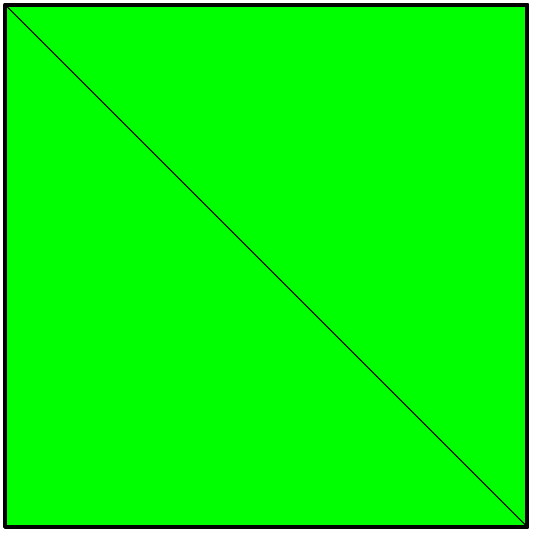}&
\includegraphics[width=20mm]{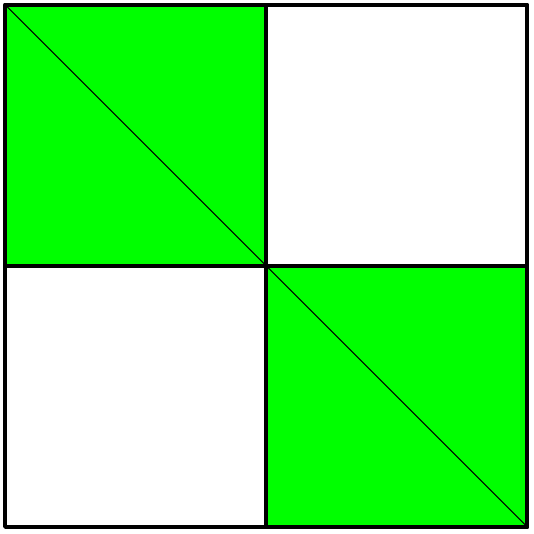}&
\includegraphics[width=20mm]{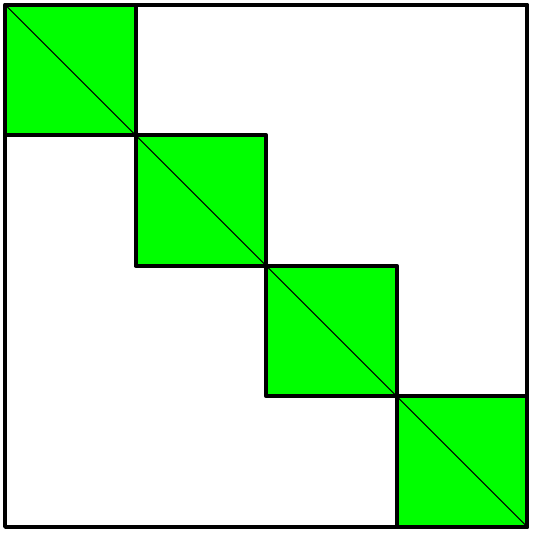}&
\includegraphics[width=20mm]{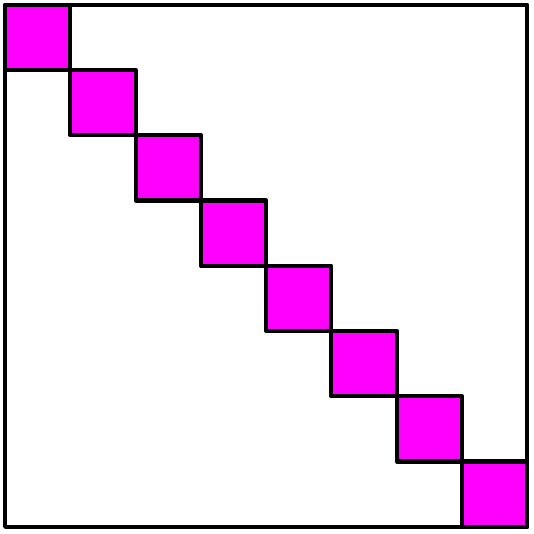}&
\includegraphics[width=20mm]{04-hodlr/figs/fig_HODLR_inv_A_simple.pdf}
\end{array}
$$
which immediately implies that $\mtx{A}^{-1} = \mtx{B}_{0}\mtx{B}_{1}\mtx{B}_{2}\mtx{B}_{3}$, so that (\ref{eq:HOLDRinvfact}) holds.

\subsection{Cluster trees and formal definition of a HODLR matrix}
\label{sec:trees}
Before closing the discussion of HODLR matrices, let us for completeness
describe a more formal way of defining the HODLR property through a hierarchical
partitioning of the index vector $I = [1,2,3,\dots,N]$. This will explain the
labeling of the blocks we used in Section \ref{sec:HODLRinv}, and will set the
ground for a discussion of more general rank structured formats in Sections \ref{sec:BIE_nested}
and \ref{sec:BIE_strong}.

\begin{figure}
\centering

\includegraphics[width=90mm]{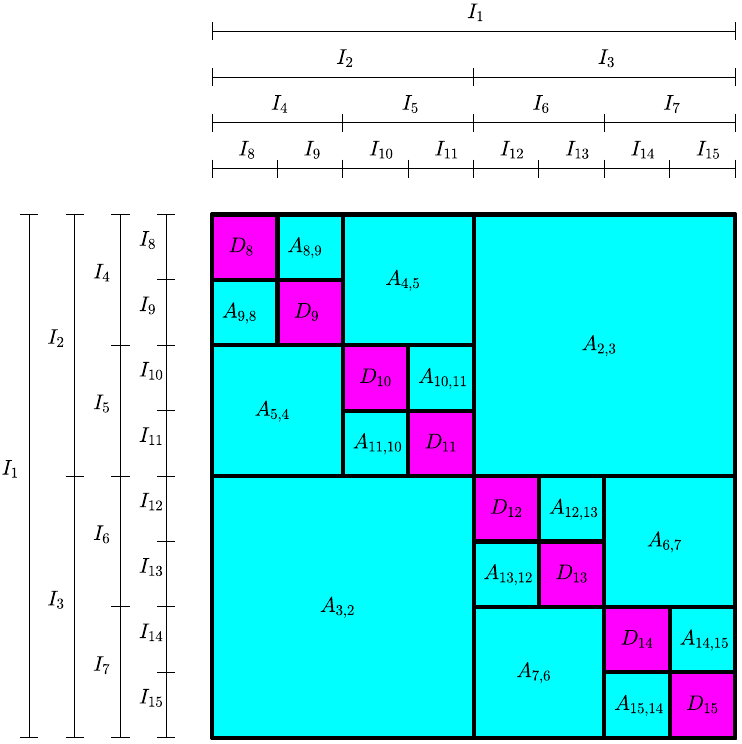}

\caption{Illustration of how a binary tree $\{I_{\tau}\}_{\tau=1}^{15}$ of subindices to the full index vector $I = [1,2,\dots,N]$ induces a hierarchical tessellation of a matrix, shown for the case of a tree of depth $L=3$. The root node is $I_{1} = I$.}
\label{fig:tree}
\end{figure}

The idea is to partition the index vector $I = [1,2,3,\dots,N]$ into a binary
tree $\{I_{\tau}\}_{\tau=1}^{15}$ of successively smaller subsets, cf.~Figure
\ref{fig:tree}.  The root of the tree is given the index $\tau=1$, and we
associate it with the full index vector, so that $I_{1} = I$. At the next finer
level of the tree, we split $I_{1}$ into two parts $I_{2}$ and $I_{3}$ so that
$I_{1} = I_{2} \cup I_{3}$ forms a disjoint partition.  Then continue splitting
the index vectors, so that $I_{2} = I_{4} \cup I_{5}$ and
$I_{3} = I_{6} \cup I_{7}$ and so forth.  For future reference, let us introduce
some basic terminology: We let $\ell$ denote a \textit{level} of the tree, with
$\ell=0$ denoting the root, so that level $\ell$ holds $2^{\ell}$ nodes.  We use
the terms \textit{parent} and \textit{child} in the obvious way, and say that a
pair of nodes $\{\alpha,\beta\}$ forms a \textit{sibling pair} if they have the
same parent.  A \textit{leaf} node is of course a node that has no children.

We now say that the matrix $\mtx{A}$ has (numerical) HODLR rank~$k$ if for every
sibling pair $\{\alpha,\beta\}$, the corresponding off-diagonal
block~$\mtx{A}_{\alpha,\beta} = \mtx{A}(I_{\alpha},I_{\beta})$ and off-diagonal
blow~\mbox{$\mtx{A}_{\beta,\alpha} = \mtx{A}(I_{\beta},I_{\alpha})$} have
(numerical) rank at most $k$.  To complete the definition, for each leaf
node~$\tau$ let~$\mtx{D}_{\tau} = \mtx{A}(I_{\tau},I_{\tau})$ denote the
corresponding diagonal block.

\subsection{Background and context}
The format that we here refer to as HODLR has been described in numerous papers, including 
\cite{2016_darve_block_low_rank_FEM,2009_martinsson_FEM,2011_bremer_rokhlin_hodlr_style,2016_kressner_rankstructured_review}. The inversion technique of Section \ref{sec:HODLRinv} is described in \cite{2013_darve_FDS,2022_martinsson_gpu_hodlr} and may well have been known earlier. Finally, we note that HODLR is a special case of the $\mathcal{H}$-matrix framework \cite{hackbusch,2002_hackbusch_H2,2004_borm_hackbusch} of Hackbusch and co-workers, which we describe further in Sections \ref{sec:BIE_nested} and \ref{sec:BIE_strong}.


%% file: 05-nesteddiss/nesteddiss06.tex
\begin{figure}
    \centering
\begin{tabular}{cccc}
$\mtx{A}$ & 
$\mtx{C} = \texttt{chol}(\mtx{A})$ & 
$\mtx{B} := \mtx{A}(I,I)$ & 
$\mtx{R} = \texttt{chol}(\mtx{B})$
\\
\includegraphics[width=20mm]{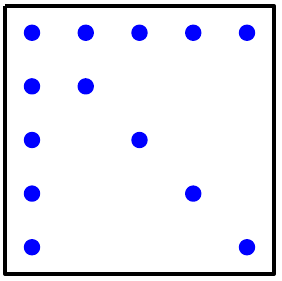} & 
\includegraphics[width=20mm]{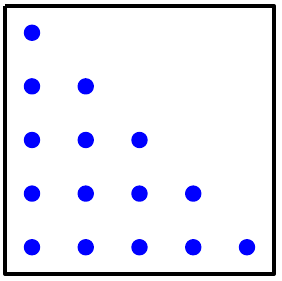} & 
\includegraphics[width=20mm]{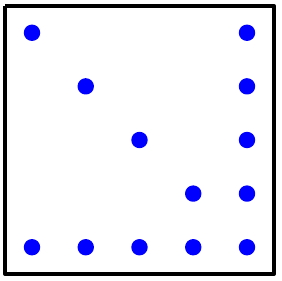} & 
\includegraphics[width=20mm]{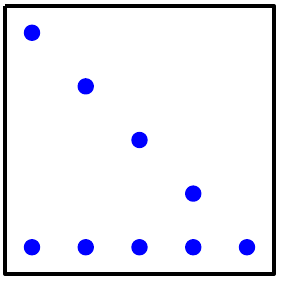} 
\\
(a) & (b) & (c) & (d)
\end{tabular}

    \caption{Illustration of the importance of ordering in reducing fill-in when factorizing sparse matrices. (a) Let $\mtx{A}$ be a positive definite matrix with the sparsity pattern show in the figure. (b) The Cholesky factor $\mtx{C}$ such that $\mtx{A} = \mtx{C}\,\mtx{C}^{*}$. (c) Form a new matrix $\mtx{B}$ by simply permuting the rows and columns of $\mtx{A}$. (d) The Cholesky factor $\mtx{R}$ such that $\mtx{B} = \mtx{R}\,\mtx{R}^{*}$.} 
    \label{fig:basic_cholesky}
\end{figure}

Having spent the last two section describing techniques for working with
structured dense matrices, let us now return to the question of how to construct
an LU factorization of a sparse matrix $\mtx{A}$ arising from the discretization
of an elliptic PDE. We saw in Section \ref{sec:onedimbvp} that for one
dimensional problems, the matrix $\mtx{A}$ is typically banded with a very
narrow band. This means that LU factorization is straightforward, as the
triangular factors inherit the bandedness property from the original matrix. We 
will now consider the two and three dimensional cases, where the
factorization problem is more challenging. It turns out that for a
typical mesh in dimension two or higher, we will quickly lose sparsity as we
attempt to compute an LU factorization of the coefficient matrix. However, the
amount of ``fill-in'' depends strongly on the ordering of the grid
(cf.~Figure \ref{fig:basic_cholesky}), and if a ``good'' ordering is chosen,
then fill-in is often a manageable issue for medium scale 2D problems and small
scale 3D problems.

In this section, we will briefly introduce the key ideas of sparse direct
solvers, and describe how rank structure in the dense matrices that arise can
be exploited to obtain close to linear complexity for problems in both two and
three dimensions.

\subsection{Nested dissection orderings and multifrontal methods}
\label{sec:NDLU}
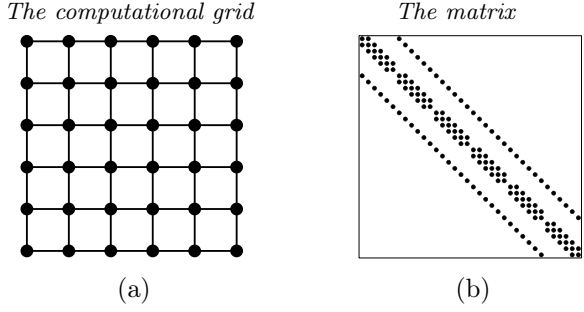
\begin{figure}[t]
\centering
\setlength{\unitlength}{1mm}
\begin{picture}(78,40)
\put(00,05){\includegraphics[width=30mm]{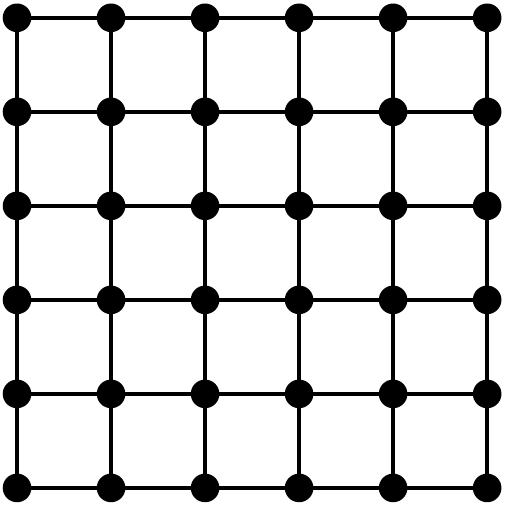}}
\put(45,05){\includegraphics[width=30mm]{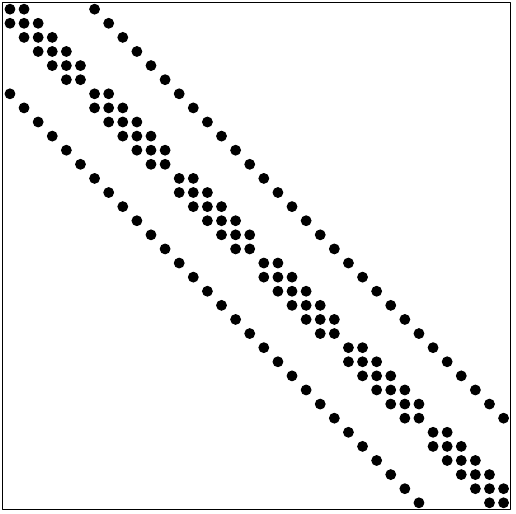}}
\put(13,00){(a)}
\put(-2,37){\small\textit{The computational grid}}
\put(58,00){(b)}
\put(50,37){\small\textit{The matrix}}
\end{picture}
\caption{The model problem involving a simple finite difference discretization
introduced in Section \ref{sec:NDLU}. (a) An $n\times n$ uniform grid, shown for $n=6$.
(b) The sparsity pattern of the $n^{2} \times n^{2}$ matrix $\mtx{A}$ in 
(\ref{eq:skyfall}). The matrix is banded, with bandwidth $n$.}
\label{fig:show_basic_grid}
\end{figure}

Our discussion will be framed around the model problem
\begin{equation}
\label{eq:buck1}
\begin{aligned}
-\Delta u(x) + m(x)\,u(x) =&\ f(x),\qquad&x \in \Omega,\\
u(x) =&\ g(x),\qquad&x \in \Gamma,
       \end{aligned}
\end{equation}
where~$\Omega = (0,1)^{2}$,~$f$ is a given body load, and~$g$ is
given Dirichlet data. (Observe that (\ref{eq:buck1}) is a two dimensional
analog of (\ref{eq:twopointbvp}).) We discretize (\ref{eq:buck1}) using the
standard five-point stencil on a regular finite difference grid of $N = n\times
n$ points, as shown in Figure \ref{fig:show_basic_grid}(a). The resulting linear
system takes the form
\begin{equation}
\label{eq:skyfall}
\mtx{A}\vct{x} = \vct{b},
\end{equation}
where the coefficient matrix $\mtx{A}$ has the familiar sparsity pattern shown
in Figure \ref{fig:show_basic_grid}(b).

Observing that the matrix in Figure  \ref{fig:show_basic_grid}(b) is banded, we
may be inclined to proceed with an LU factorization as we did in Section
\ref{sec:onedimbvp}, and put our faith in the fact that the triangular factors
will be banded just like the original matrix. The complication that arises is
that the bands will get almost completely filled in. Since the width~$b$ of the
band in two dimensions satisfies $b \sim N^{1/2}$ (unlike in one dimension
where $b = \cO(1)$), computing the triangular factors requires $\sim b^2 N \sim
N^2$ flops. For the analogous problem in three dimensions on an $N = n\times n
\times n$ grid, we would find a bandwidth of $b\sim N^{2/3}$ and a catastrophic
overall complexity of $\sim b^2 N \sim N^{7/3}$.

\begin{figure}[t!]
\centering
\begin{center}
\setlength{\unitlength}{1mm}
\begin{picture}(120,52)
\put(010,12){\includegraphics[width=40mm]{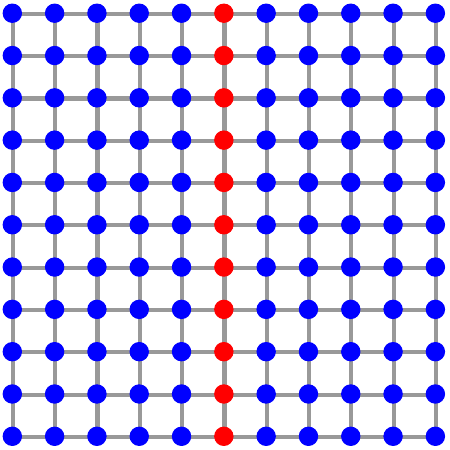}}
\put(029,09){\color{red} $I_{1}$}
\put(005,30){\color{blue}$I_{2}$}
\put(052,30){\color{blue}$I_{3}$}
\put(027,00){(a)}
\put(070,12){\includegraphics[width=40mm]{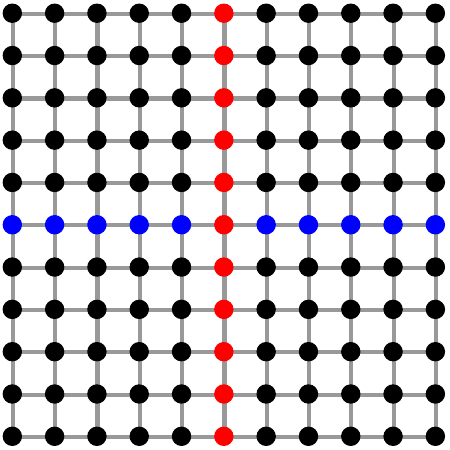}}
\put(089,09){\color{red}$I_{1}$}
\put(065,30){\color{blue}$I_{2}$}
\put(112,30){\color{blue}$I_{3}$}
\put(065,20){\color{black}$I_{4}$}
\put(065,40){\color{black}$I_{5}$}
\put(112,20){\color{black}$I_{6}$}
\put(112,40){\color{black}$I_{7}$}
\put(087,00){(b)}
\end{picture}
\end{center}
    \caption{Dissection of a computational mesh, as described in Section \ref{sec:NDLU}. (a) An $n\times n$ finite difference mesh. We excise a line of nodes in $I_{1}$ (red) to separate the remaining nodes into two unconnected subdomains $I_{2}$ and $I_{3}$ (blue). (b) One additional step of grid partitioning results in four unconnected subdomains (black).}
    \label{fig:nd1}
\end{figure}
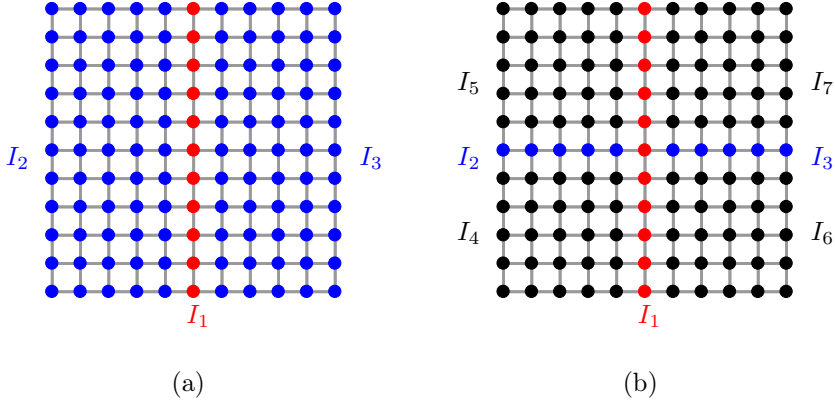

The key to reducing fill-in is to first reorder the nodes in the mesh. To
illustrate, let us consider the popular ``nested dissection'' ordering, which is
based on recursively splitting the mesh into roughly equisized disconnected
pieces by eliminating a small amount of nodes from the mesh. For the present
geometry, we could begin by by cutting out the nodes
marked~$I_{1}$ (red) in Figure~\ref{fig:nd1}(a) to split the mesh into the two
disconnected parts marked $I_{2}$ and $I_{3}$. In other words,
$$
I = [I_{3},\, I_{2},\, I_{1}]
$$
forms an ordering of the degrees of freedom in which the matrix
$\mtx{A}(I,I)$ has zero blocks in the $(2,3)$ and $(3,2)$ position:
$$
\mtx{A}(I,I) =
\left[\begin{array}{rrr}
\mtx{A}_{33} & \mtx{0}      & \mtx{A}_{31} \\
\mtx{0}      & \mtx{A}_{22} & \mtx{A}_{21} \\
\mtx{A}_{13} & \mtx{S}_{12} & \mtx{A}_{11}
\end{array}\right].
$$
This means that if we \textit{temporarily} assume that we have
access to a mechanism for building the triangular factorizations
$$
\mtx{A}_{22} = \mtx{L}_{22}\mtx{U}_{22}
\quad\mbox{and}\quad
\mtx{A}_{33} = \mtx{L}_{33}\mtx{U}_{33},
$$
then we immediately obtain the partial LU factorization
$$
\mtx{A}(I,I)
=
\left[\begin{array}{c|c|c}
\mtx{L}_{33}                  & \mtx{0}                       & \mtx{0} \\ \hline
\mtx{0}                       & \mtx{L}_{22}                  & \mtx{0} \\ \hline
\mtx{A}_{13}\mtx{U}_{33}^{-1} & \mtx{A}_{12}\mtx{U}_{22}^{-1} & \mtx{I}
\end{array}\right]
\left[\begin{array}{c|c|c}
\mtx{I} & \mtx{0} & \mtx{0} \\ \hline
\mtx{0} & \mtx{I} & \mtx{0} \\ \hline
\mtx{0} & \mtx{0} & \mtx{S}_{11}
\end{array}\right]
\left[\begin{array}{c|c|c}
\mtx{U}_{33} & \mtx{0}      & \mtx{L}_{33}^{-1}\mtx{A}_{31} \\ \hline
\mtx{0}      & \mtx{U}_{22} & \mtx{L}_{22}^{-1}\mtx{A}_{21} \\ \hline
\mtx{0}      & \mtx{0}      & \mtx{I}
\end{array}\right],
$$
where
$
\mtx{S}_{11} =
\mtx{A}_{11} - \mtx{A}_{13}\mtx{U}_{33}^{-1}\mtx{L}_{33}^{-1}\mtx{A}_{31} -
               \mtx{A}_{12}\mtx{U}_{22}^{-1}\mtx{L}_{22}^{-1}\mtx{A}_{21}
$
is a Schur complement. The factorization is completed by factoring $\mtx{S}_{11} = \mtx{L}_{11}\mtx{U}_{11}$, so that
\begin{equation}
\label{eq:pamart}
\mtx{A}(I,I) = 
\left[\begin{array}{c|c|c}
\mtx{L}_{33}                  & \mtx{0}                       & \mtx{0} \\ \hline
\mtx{0}                       & \mtx{L}_{22}                  & \mtx{0} \\ \hline
\mtx{A}_{13}\mtx{U}_{33}^{-1} & \mtx{A}_{12}\mtx{U}_{22}^{-1} & \mtx{L}_{11}
\end{array}\right]
\left[\begin{array}{c|c|c}
\mtx{U}_{33} & \mtx{0}      & \mtx{L}_{33}^{-1}\mtx{A}_{31} \\ \hline
\mtx{0}      & \mtx{U}_{22} & \mtx{L}_{22}^{-1}\mtx{A}_{21} \\ \hline
\mtx{0}      & \mtx{0}      & \mtx{U}_{11}
\end{array}\right].
\end{equation}

To summarize, we build the factorization (\ref{eq:pamart}) by executing three steps:
\begin{itemize}
\item Factorize $\mtx{A}_{33}$:
$[\mtx{L}_{33},\mtx{U}_{33}] = \texttt{lu}(\mtx{A}_{33})$ \hfill \textit{\color{blue}size $\sim N/2 \times N/2$}
\item Factorize $\mtx{A}_{22}$:
$[\mtx{L}_{22},\mtx{U}_{22}] = \texttt{lu}(\mtx{A}_{22})$ \hfill \textit{\color{blue}size $\sim N/2 \times N/2$}
\item Factorize $\mtx{S}_{11}$:
$[\mtx{L}_{11},\mtx{U}_{11}] = $\\ \tabto{.5in} $\texttt{lu}(\mtx{A}_{11} - 
                                           \mtx{A}_{12}\mtx{U}_{22}^{-1}\mtx{L}_{22}^{-1}\mtx{A}_{21} -
                                           \mtx{A}_{13}\mtx{U}_{33}^{-1}\mtx{L}_{33}^{-1}\mtx{A}_{31})$
      \hfill \textit{\color{blue} size $\sim \sqrt{N} \times \sqrt{N}$}
    \end{itemize}
The idea is now to execute the factorizations of $\mtx{A}_{22}$ and $\mtx{A}_{33}$
by applying the same idea recursively. We first split them in two halves to get the
partition shown in Figure~\ref{fig:nd1}(b), and then continue down until the
pieces are small enough that brute force dense factorization can be done efficiently.

\begin{figure}[t!]
\begin{center}
\begin{tabular}{ccc}
\textit{Levels 0 (red) and 1 (blue) } &
\textit{Levels 2 (red) and 3 (blue) } &
\textit{Level 4, the leaves} \\
\includegraphics[width=0.25\textwidth]{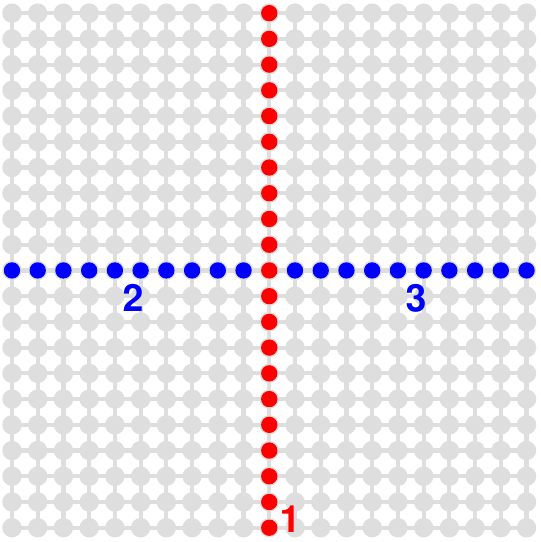} &
\includegraphics[width=0.25\textwidth]{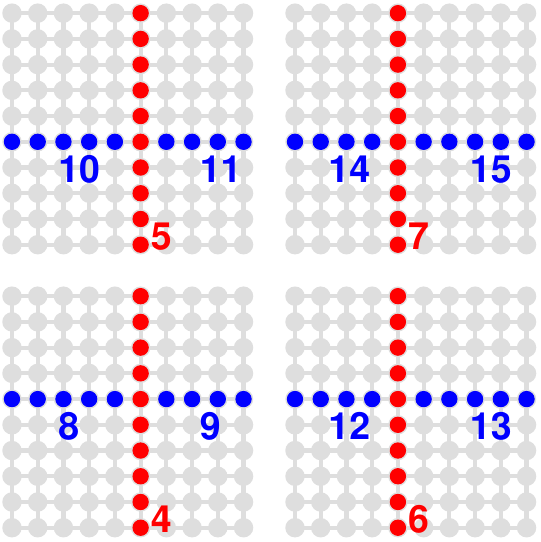} &
\includegraphics[width=0.25\textwidth]{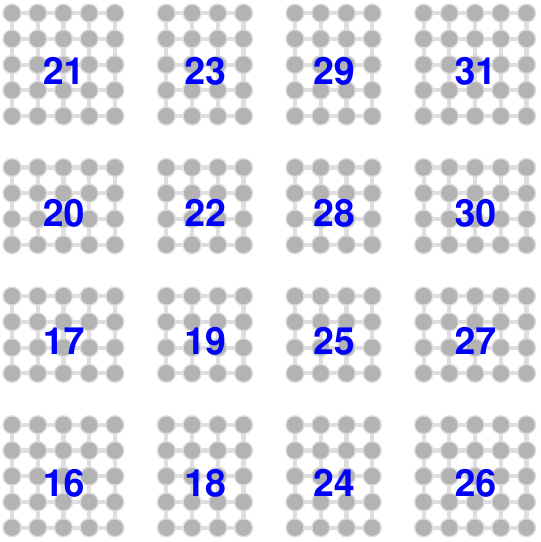} \\[4mm]
(a) & (b) & (c)
\end{tabular}
\end{center}
    \caption{Nested dissection ordering of a uniform grid on a square. In this example, the tree has five levels.}
    \label{fig:nd2}
\end{figure}

\begin{figure}[t!]
\centering

\setlength{\unitlength}{1mm}
\begin{picture}(120,73)
\put(00,37){\includegraphics[height=43mm]{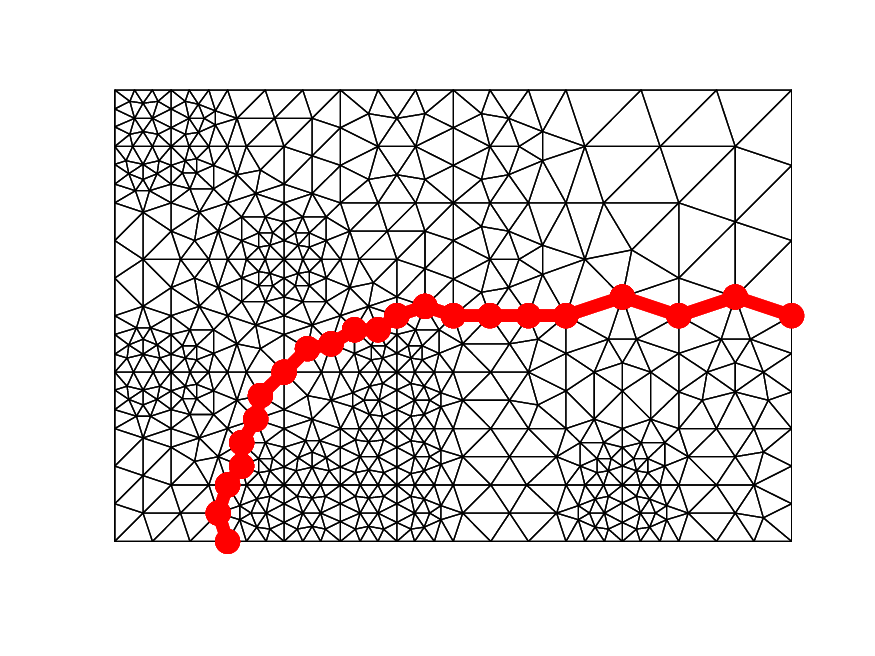}}
\put(60,37){\includegraphics[height=43mm]{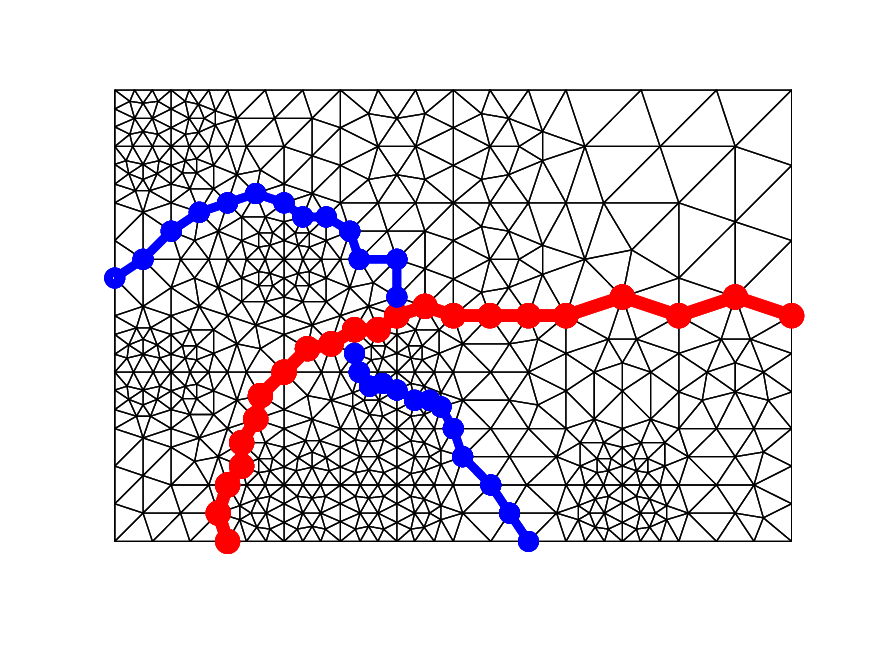}}
\put(00,-3){\includegraphics[height=43mm]{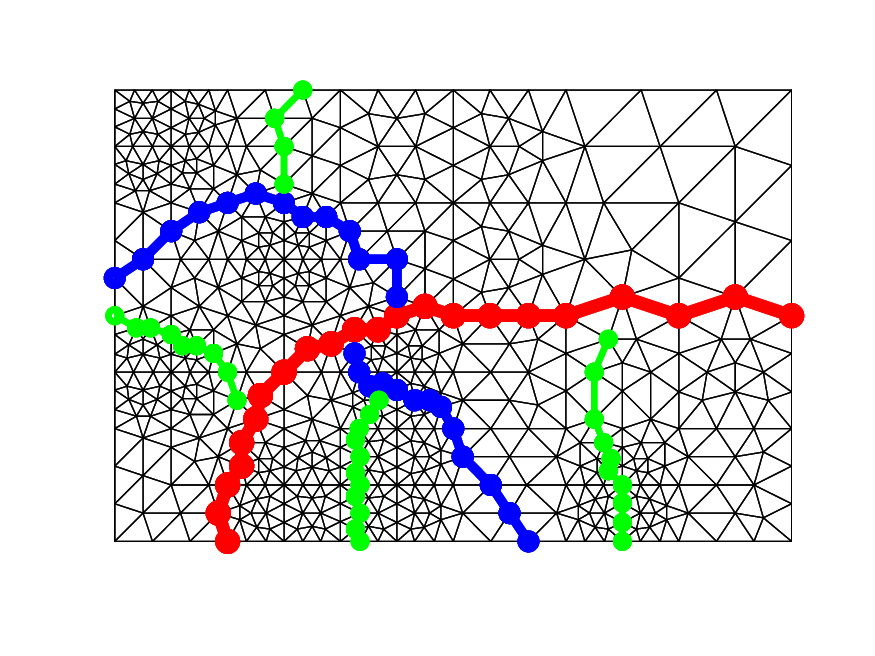}}
\put(60,-3){\includegraphics[height=43mm]{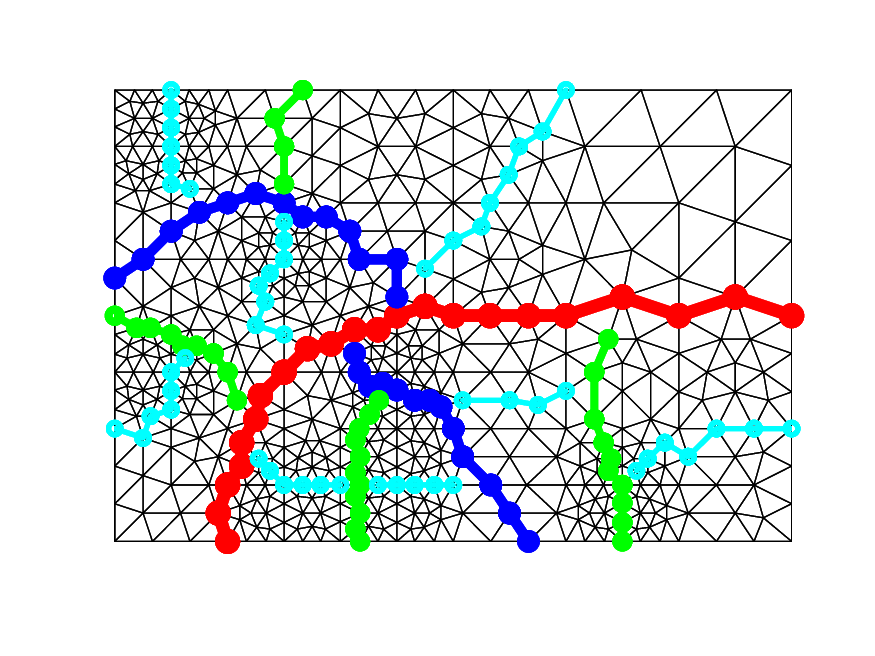}}
\put(10,40){\small(i) \textit{Top level separator}}
\put(70,40){\small(ii) \textit{Two levels of separators}}
\put(08,00){\small(iii) \textit{Three levels of separators}}
\put(70,00){\small(iv) \textit{Four levels of separators}}
\end{picture}

\caption{Nested dissection applied to an unstructured grid. Reproduced with permission from~\cite{2012_xia_robust_efficient_multifrontal}.}
\label{fig:nd3}
\end{figure}
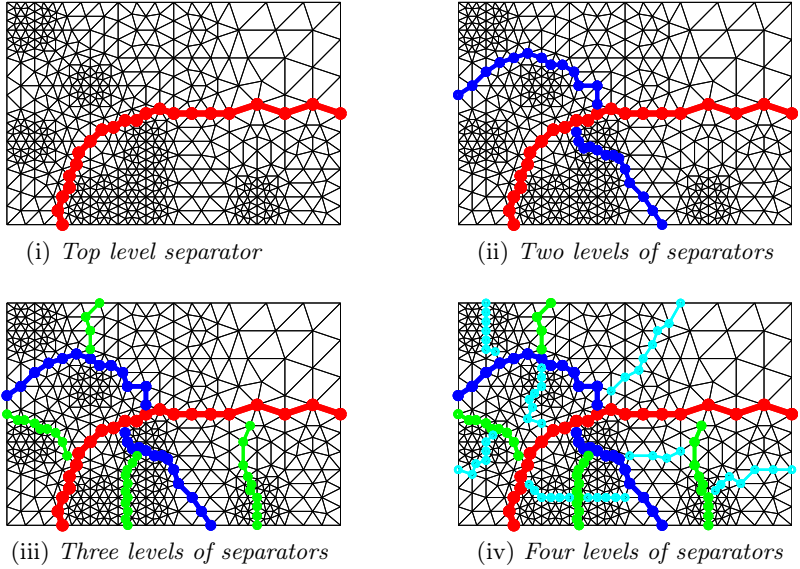

To illustrate the tree structure that results from nested dissection, Figure
\ref{fig:nd2} shows a slightly larger grid that is partitioned into 4 separate
levels. An example involving a less structured grid is shown in Figure
\ref{fig:nd3}.

\subsection{Asymptotic complexity of multifrontal methods}
\label{sec:sparsedirectcomplexity}

It is a relatively simple matter to work out the asymptotic complexity of the
scheme described. Let $L$ denote the number of levels in the tree, and let
$\ell \in [0,\,1,\,2,\,\dots,\,L]$ denote one specific level, with $\ell = 0$
corresponding to the root of the tree. For an $n\times n$ grid, there are
$m_{\ell} \sim 2^{-\ell/2}n$ points in any cell at level $\ell$, and there are
$2^{\ell}$ cells. The dominant cost in processing a cell is the dense
factorization of the $m_{\ell} \times m_{\ell}$ Schur complement. In
consequence, the total cost becomes
$$
T_{\rm build}
\sim \sum_{\ell=0}^{L}2^{\ell}m_{\ell}^{3}
\sim \sum_{\ell=0}^{L}2^{\ell}\bigl(n 2^{-\ell/2}\bigr)^{3}
\sim \sum_{\ell=0}^{L}2^{-\ell/2}\,n^{3}
\sim n^{3}
\sim N^{3/2},
$$
as $n \sim N^{1/2}$ in this case.

Let us next consider the analogous
problem in three dimensions with $N = n\times n \times n$ grid points.
In this case, $m_{\ell} \sim 2^{-2\ell/3}n^{2}$,
which means that the factorization has cost
$$
T_{\rm build}
\sim \sum_{\ell=0}^{L}2^{\ell}m_{\ell}^{3}
\sim \sum_{\ell=0}^{L}2^{\ell}\bigl(n^{2} 2^{-2\ell/2}\bigr)^{3}
\sim \sum_{\ell=0}^{L}2^{-\ell}\,n^{6}
\sim n^{6}
\sim N^{2},
$$
since $n \sim N^{1/3}$ in three dimensions.

To summarize, the costs of performing an LU factorization of the sparse
coefficient matrix resulting from a finite difference discretization of an
elliptic PDE on a box in two or three dimensions can -- by using a nested
dissection ordering -- be reduced to:

\vspace{1.5mm}

\begin{center}
\begin{tabular}{c|c|c|c|c}
& \textit{complexity of factorization} & \textit{complexity of solve} & \textit{storage} \\
\hline
2D & $\cO(N^{3/2})$ & $\cO(N\log N)$ & $\cO(N\log N)$ \\
3D & $\cO(N^{2})$ & $\cO(N^{4/3})$ & $\cO(N^{4/3})$.
\end{tabular}
\end{center}

\vspace{1.5mm}

\noindent These asymptotic costs are optimal in the sense that no other
ordering does better. While we considered only the special case of a finite
difference discretization on a uniform grid, the complexities cited are typical
for a range of different discretizations (finite elements, finite volumes,
etc.), and for most quasiuniform meshes.

Despite the unfavorable asymptotic complexities of sparse direct solvers, they
are popular among practitioners as they are very simple to use, and numerically
robust. For two dimensional problems, one can easily handle tens of millions of
degrees of freedom, and in three dimensions, about one million degrees of
freedom is generally not a problem. These numbers are large enough for many
practical applications.\footnote{The numbers quoted here refer to low order
discretizations. Most high order discretization methods lead to less sparse
coefficient matrices, which tends to greatly increase computational costs for
sparse direct solvers. Specialized methods designed to miminize fill-in do
better, however \cite{2012_spectralcomposite}.}

\begin{remark}[Pivoting and numerical stability]
  For LU factorization to be numerically stable, it is typically necessary to
  pivot. For dense matrices, it is well understood how to do this, and standard
  partial pivoting is in practice highly numerically stable, and very amenable
  to high performance implementations. For sparse matrices, the situation is
  more subtle. There is a difficult-to-navigate trade-off between choosing an
  ordering of the variables that minimizes fill-in, versus one that leads to
  minimal growth of the elements in the triangular factors. Additionally, high
  performance implementations often defer global updates, and work on only small
  parts of the matrix at one time. Some software packages pick an elimination
  order in advance, based only on an analysis of the connectivity graph, and
  will then deviate from the preset order only if unusually large growth factors
  are encountered.
\end{remark}

\subsection{Acceleration via rank structure}
\label{sec:fastLU}

We saw in Section \ref{sec:sparsedirectcomplexity} that the arithmetic cost of
sparse LU is dominated by the factorization of the dense Schur complements at
the top levels. Remarkably, it turns out that these dense blocks share many
characteristics of discretized integral operators. In particular, they are rank
structured and admit fast matrix arithmetic (see
Section~\ref{sec:dissip_frontal} for numerical examples).

To illustrate how the introduction of rank structured matrix algebra impacts the
computational costs of sparse direct solvers, let us revisit the problem of
computing the LU factorization of a finite difference matrix on a regular
$N = n\times n$ grid in 2D. Using the same notation as in Section
\ref{sec:sparsedirectcomplexity}, we again let $m_{\ell}$ denote the number of
points in a front at level $\ell$, so that $m_{\ell} \sim 2^{-\ell/2}n$.
Observe further that if we continue the dissection process until there are
$O(1)$ points in each leaf node, then the number of levels $L$ in the tree
satisfies $2^{L} \sim N$, so that $m_{\ell} \sim 2^{-(L+\ell)/2}\,N$. Suppose
now that our dense matrix operations have \textit{linear complexity}, so that
the cost of forming a Schur complement of size $m_{\ell} \times m_{\ell}$ and
then computing its LU factorization is $O(m_{\ell})$. Since there are $2^{\ell}$
fronts at level $\ell$, we then find that the cost to process level $\ell$ is
$$ C_{\ell} \sim 2^{\ell} \times 2^{-(L+\ell)/2}N = 2^{(L-\ell)/2}n. $$ Let us
compare the per level costs of the classical factorization methods to the one
accelerated by rank structured matrix algebra:

\vspace{1.5mm}

{\footnotesize
\begin{tabular}{l||c|c|c|c|c|c}
& $\ell=L$ & $\ell = L-1$ & $\ell = L-2$ & $\ell = L-3$ & $\cdots$ & $\ell = 0$\\
& (leaf level) &&&&& (root level) \\ \hline
Classical   & $N$ & $2^{ 1/2}N$ & $2^{ 2/2}N$ & $2^{ 3/2}N$ & $\cdots$ & $2^{ L/2}N \sim N^{3/2}$ \\
Accelerated & $N$ & $2^{-1/2}N$ & $2^{-2/2}N$ & $2^{-3/2}N$ & $\cdots$ & $2^{-L/2}N \sim N^{1/2}$
\end{tabular}
}

\vspace{1.5mm}

\noindent In other words, in the accelerated scheme, the costs
\textit{decrease} as we move towards the coarser levels, and the overall cost
is dominated by processing the leaves. In addition to obtaining overall
linear~$\cO(N)$ complexity, we have a secondary benefit in that the step that
is now the most expensive, is also the one that is the easiest to parallelize!
(For further details, including the analogous calculation for 3D problems, see
\cite{2013_xia_rankpatterns} and
\cite[Ch.~21]{2019_martinsson_fast_direct_solvers}.)

In practice, realizing linear complexity in factorizing the intermediate Schur
complements is challenging. The accelerated scheme involves complicated data
structures with nested hierarchies: an ``outer'' hierarchy provided
by the nested dissection ordering of the mesh, and then an ``inner''
hierarchy that organizes the points in a front into a tree that defines the rank
structured matrix algebra. Coding this efficiently turns out to be sufficiently
complex that many groups rely on simplified formats. For instance, Schur
complements are often formed as full dense matrices, and then compressed in
order to admit fast factorization, making the cost to process a
leaf~$\cO(m_{\ell}^{2})$ rather than~$\cO(m_{\ell})$. Such schemes still perform
well in practice, but do not attain true linear complexity.  On the other hand,
limiting the multi-level nature of such schemes in particular elongated or
\emph{slab} geometries can also result in extremely performan algorithms,
especially when coupled with GPU accelerated dense matrix
algebra~\cite{Yesypenko_2024}.

An important question for a practitioner is what ``format'' to use for the rank
structured matrix algebra. This depends on the underlying dimension and the
size of the problem. In 2D, simple data structures such as HODLR are almost
always sufficient. For large scale 3D problems, more sophisticated structures
are typically required. We will return to this question in
Section~\ref{sec:dissipation}, after we have introduced the notion of ``strong
admissibility'' in Section~\ref{sec:BIE_strong}. (Recall that for 2D problems,
the rank structured matrix ``lives'' on a front of dimension \textit{one},
which explains why simple formats with weak admissibility work well,
cf.~Section~\ref{sec:dimensionissue}. Analogously, in 3D, the fronts are for
practical purposes two dimensional.)

\subsection{Background and context}
Early work discussing how choosing a good mesh ordering greatly improves the
performance of sparse LU includes \cite{george_1973} and \cite{hoffman_1973},
which establish that $\cO(N^{3/2})$ for 2D and $\cO(N^2)$ for 3D are both
attainable and optimal. Subsequent work refined the technology
\cite{1989_directbook_duff,2006_davis_directsolverbook} and led to the
development of widely used software packages such as UMFPACK, SuiteSparse,
SuperLU, MUMPS, Pardiso, and Strumpack, to name just a few. The 2016 survey
\cite{2016_acta_sparse_direct_survey} provides details.

The realization that rank structured matrix algebra can be used to further
improve the asymptotic complexity dates back at least to
\cite{2007_leborne_HLU,2009_xia_superfast,2009_martinsson_FEM,2010_shiv_PDEranks}.
Subsequent years saw significant advances in the field \cite{2019_amestoy_flat,
  2019_amestoy_BLR_multifrontal, 2016_darve_block_low_rank_FEM,
  2017_pouransari_darve, 2017_pichon_sparse, 2011_xia_dehoop,
  2012_spectralcomposite}, including the development of black-box randomized
algorithms for constructing the rank-structured representations of Schur
complements~\cite{2011_martinsson_randomhudson, 2011_lin_lu_ying,
  2013_xia_randomized_FDS, 2016_ghysels_randomized_multifrontal,
  levitt2024linear}.  A related but slightly different approach was pursued in
\cite{2016_ying_HIFDE} where a method known as Hierarchical Interpolative
Factorization was presented. This technique attains linear or quasilinear
complexity without the double hierarchies, by exploiting rank structure to
further compress the frontal separators before moving up the hierarchical data
structure.


%% file: 07-BIE2D/BIE2D02.tex
In Section \ref{sec:onedimbvpIE}, we showed how to convert a one-dimensional
boundary value problem to an integral equation. This was useful because it gave
us a second-kind Fredholm equation that was inherently well-conditioned.  Now we
will show how to do the same thing for two-dimensional problems. For certain
types of equations (like Laplace or Helmholtz equations with constant
coefficients), we get an extra bonus: we only need to work on the boundary of
the region, not the entire area inside. This reduces the effective dimension of
the problem from two to one (or from three to two for three-dimensional
problems). The trade-off is that we end up with a dense matrix, but this matrix
is rank structured and amenable to FMM accelerated iterative solvers, or to
direct solvers such as those introduced in Section~\ref{sec:HODLR}.

\subsection{Formulation}
\label{sec:BIE_formulation}
To make matters concrete, let us focus on the basic interior Dirichlet problem
for the Laplace equation:
\begin{equation}
\label{eq:nicholas}
\begin{aligned}
-\Delta u(\pxx) =&\ 0,\qquad& \pxx \in \Omega,\\
u(\pxx) =&\ f(\pxx),\qquad&\pxx \in \Gamma,
\end{aligned}
\end{equation}
where $\Omega$ is a domain in two dimensions with boundary $\Gamma$. A thorough
introductory treatment of boundary integral equations can be found
in~\cite{colton2012, guenther1996partial}. The integral equation reformulation
we will describe is versatile and effective in situations involving complicated
geometries, but in order to avoid needless technical complications, we will assume
that~$\Omega$ is simply connected and that~$\Gamma$ is smooth.

There are several ways to reformulate~\eqref{eq:nicholas} as an integral
equation. The classical \emph{indirect} approach is to write~$u$ as a
\emph{double layer potential} with unknown density~$\sigma$:
\begin{equation}
  \label{eq:double}
  u(\pxx) =  \int_{\Gamma} d(\pxx,\pyy) \, \sigma(\pyy) \, ds(\pyy) 
    = \cD[\sigma] (\pxx),
\end{equation}
for~$\pxx \in \Omega$ and where~$d$ is the double layer potential
\begin{equation}
\label{eq:def_d}
  d(\pxx,\pyy) = \pvct{n}(\pyy)\cdot \nabla_y \phi(\pxx-\pyy),
\end{equation}
the function~$\phi$ is the free space fundamental solution associated
with the Laplace operator in two dimensions,
\begin{equation}
  \label{eq:g02d}
  \phi(\pxx-\pyy) = -\frac{1}{2\pi}\log|\pxx-\pyy|,
\end{equation}
and~$\pvct{n}$ is the outward unit normal vector to~$\Gamma$.
If we take the limit as~$\pxx \to \Gamma$, we arrive at the following second-kind
integral equation for the new unknown~$\sigma$:
\begin{equation}
  \label{eq:intdir}
  -\frac{1}{2} \sigma(\pxx) + 
\int_{\Gamma} d(\pxx,\pyy) \, \sigma(\pyy) \, ds(\pyy) = f(\pxx), \qquad \text{for } \pxx \in \Gamma.
\end{equation}
A full analysis of the integral equation (\ref{eq:intdir}) (e.g. uniqueness,
regularity, etc.) can be found in~\cite{colton2012,2013_colton_kress}, along
with various details and other boundary integral formulations for other related
boundary value problems (e.g. the Helmholtz equation, exterior problems, etc.).
The integral equation (\ref{eq:intdir}) will admit various hierarchical low-rank
factorizations after it is discretized, which we briefly discuss next.

\subsection{Discretization}

In order to numerically solve the BVP~(\ref{eq:nicholas}), one could either use
(\ref{eq:nicholas}) or an equivalent BIE, such as~(\ref{eq:intdir}), as a
starting point. The key advantages of the integral equation formulation are that
it involves a problem on a 1D domain rather than a 2D domain, and that it is
inherently much better conditioned (for instance, when~$\Gamma$ is
smooth,~(\ref{eq:intdir}) is a second-kind Fredholm equation
on~$L^{2}(\Gamma)$).  The key disadvantage of~(\ref{eq:intdir}) is that it
involves a global operator, and therefore leads to a linear system with a dense
coefficient matrix.

There are several ways to discretize the boundary integral 
equation~\eqref{eq:intdir}, cf.~\cite{brebbia2016boundary}, but we will use the simplest Nystr\"om 
scheme~\cite{atkinson1997} for
the following discussion.  For the specific problem introduced in Section \ref{sec:BIE_formulation}, it can be shown that the
double layer kernel~$d$ is not singular at~$\pxx=\pyy$, in fact for~$\pxx,\pyy \in \Gamma$,
\begin{equation}
  \lim_{\pxx \to \pyy} d(\pxx,\pyy) = \frac{1}{4\pi} \kappa(\pxx),
\end{equation}
where~$\kappa(\pxx)$ denotes the curvature of~$\Gamma$ at the point~$\pxx$, and that
there exists a well-known \emph{quadrature rule} for discretizing the continous
integral equation into a discrete one via a Nystr\"om
method~\cite{atkinson1997,colton2012}:
\begin{equation}
  -\frac{1}{2} \sigma(\pxx_{i}) + \sum_{j=1}^{N} w_{j} \, d(\pxx_{i}, \pxx_{j})
  \, \sigma(\pxx_{j}) = f(\pxx_{i}).
\end{equation}
If the curve~$\Gamma$ has length~$L$ and the points~$\pxx_{j}$ are chosen to be
equispaced in arclength along~$\Gamma$, then choosing~$w_{j} = L/N = h$
yields an efficient discretization of the equation~\cite{colton2012}.
The discretized finite dimensional linear
system can then be written as
\begin{equation}
  \label{eq:intdir2}
 \left(  -\frac{1}{2} \mtx{I} + h \mtx{A}  \right) \vct{s} = \vct{f},
\end{equation}
with
\begin{equation*}
  \vct{s}(i) \approx \sigma(\pxx_i), \qquad \vct{f}(i) = f(\pxx_i),
  \quad \text{and} \quad \mtx{A}(i,j) = d(\pxx_i,\pxx_j).
\end{equation*}
It remains to develop an efficient numerical method for solving the linear
system (\ref{eq:intdir2}).

\begin{remark}
The double layer potential (\ref{eq:def_d}) is an unusual kernel in that
it is smooth at the diagonal. Almost all other integral
equation formulations involve singular kernels that require 
specialized quadrature rules in order to obtain high accuracy.
\end{remark}

\subsection{Efficient compression and inversion}

For many of the same reasons that the linear system arising from the
discretization of the integral equation corresponding to the 1D boundary
value problem in Section~\ref{sec:onedimbvpIE} is data-sparse, the
matrix~$\mtx{A}$ in~\eqref{eq:intdir2} is also data sparse. In fact, using
analytical properties of the kernel~$d$ of the double layer potential, the
matrix~$\mtx{A}$ can be compressed in hierarchical low-rank format and applied
to vectors in~$\mathcal O(N)$ time, as well as inverted in the same asymptotic
complexity.

A more detailed description of this analytic compression is given in
Section~\ref{sec:compression}.

\begin{remark}[Linear algebraic vs.~analytic dimension reduction]
  The boundary integral equation reformulation described in this section and the
  nested dissection-based LU factorization described in Section~\ref{sec:FDSPDE}
  are similar in that for both of them, the dominant computational cost concerns
  the factorization of a dense matrix that is defined on a domain of lower
  dimensionality than the original BVP. One advantage to doing this using
  analytic tools (instead of agnostic numerical ones) is that superior numerical
  stability can be attained by doing the conversion exactly. However, the
  analysis-based approach is not always available, as only some differential
  equations can, in a practical way, be reformulated as integral equations.
\end{remark}

We now take an opportunity to discuss a few other common regimes for which we
can derive well-conditioned integral equations for PDEs.

\subsection{Other type of integral equations}

Integral equation formulations can be derived for a broad class of differential
equations, we describe some of the more popular ones below.  Each of the
following classes of problem is amenable to the construction of an associated
fast direct solver.

\subsubsection*{Exterior problems}

While the boundary value problem in~\eqref{eq:nicholas} concerns finding the
solution in a bounded region, often times in classical mathematical physics the
solution is required in the unbounded exterior of a domain~$\Omega$, and the
boundary data (i.e. the function~$f$ above) is obtained from an
\textit{incoming} potential or wave. In the case of electrostatics, this might
be the electrostatic potential due to a collection of charges at some distance
from~$\Omega$. In the acoustic (i.e. Helmholtz) or electromagnetic
(i.e. Maxwell) case, the data might be obtained from evaluating a plane wave on
the boundary~$\Gamma$. Exterior boundary value problems result in integral
equations which are often formally the adjoint of interior problems, and
analogous analysis holds in both situations. Second-kind integral equation
formulations exist for these problems with kernels that are merely linear
combinations of derivatives of the fundamental solution.

\subsubsection*{Inhomogenous problems}

Nearly identical integral equations can be derived for inhomogeneous boundary
value problems as well, for example the standard interior Poisson problem:
\begin{equation}
\label{eq:poisson}
\begin{aligned}
-\Delta u(\pxx) =&\ g(\pxx),\qquad&\pxx \in \Omega,\\
u(\pxx) =&\ f(\pxx),\qquad&\pxx \in \Gamma.
\end{aligned}
\end{equation}
In this case, a slightly different ansatz for the solution~$u$ must be made
which explicitly satisfies the inhomogeneity:
\begin{equation}
  \label{eq:vp}
  u(\pxx) = \cD[\sigma] (\pxx) + \cV[g] (\pxx)  
\qquad \text{for } \pxx \in \Omega,
\end{equation}
where~$\cD$ and~$\cV$ are the double layer and volume potentials previously defined.
Again, the volume potential term automatically
satisfies the inhomogeneity, and the resulting integral equation is the same as
for the interior Dirichlet problem except with a different right hand side:
\begin{equation}
  \label{eq:intpos}
  -\frac{1}{2} \sigma(\pxx) +  \cD[\sigma](\pxx) = 
      f(\pxx) - \cV[g](\pxx), \qquad \text{for } \pxx \in \Gamma.
\end{equation}
The integral operator is identical to that derived for the interior Dirichlet
problem. One additional numerical challenge that remains is to accurately
evaluate the volume potential $\mathcal{V}$, but this does not affect
the compression or inversion of the equation.

\subsubsection*{Variable coefficient problems}

For some variable coefficent problems it is possible to reformulate them into
what are referred to as \emph{volume} integral equations. In this case, there
is no dimension reduction of the problem (which is expected, due to the
variable coefficients), but the resulting integral equation is often very
well-conditioned and inherits many of the same properties of boundary integral
equations~\cite{2013_colton_kress, colton1998inverse,colton1988inverse}.

A simple example of such a variable coefficient PDE is the Helmholtz equation
with background wavenumber~$\kappa$ augmented with a compactly supported sound speed
perturbation~$q > -1$, and compactly supported right hand side~$f$:
\begin{equation*}
\left(  \Delta + \kappa^2(1+q) \right) u = f, \qquad \text{in } \bbR^2.
\end{equation*}
Writing the solution~$u$ using the Green's function~$\phi_\kappa$
for the Helmholtz equation,
\begin{equation*}
    u(\pxx) = \int_{\bbR^2} \phi_\kappa(\pxx-\pyy) \, \sigma(\pyy) \, da(\pyy) 
    = \cV_{\kappa}[\sigma](\pxx) 
\end{equation*}
yields what is known as the Lippmann-Schwinger equation:
\begin{equation*}
  \sigma(\pxx) + \kappa^2 q(\pxx) \cV_\kappa[\sigma](\pxx) = f(\pxx).
\end{equation*}
Note that outside of the support of~$q$ and~$f$, it is clear that the
solution~$\sigma$ is identically zero, and therefore the above volume integral
equation only needs to be discretized on a finite domain. The compressibility
of the integral operator $\mathcal{V}_{\kappa}$ requires a bit of extra care when compared to
ones merely defined on curves (as opposed to volume filling regions), but the
main considerations are the same~\cite{gillman2014direct,
  2012_spectralcomposite, ambikasaran2016fast,
  2014_corona_martinsson_BIE_plane, gopal2022accelerated, 2016_ying_HIF_IE}.


%% file: 08-nested/nested_new03.tex
In Section \ref{sec:BIE2D}, we described how to reformulate a boundary value
problem such as (\ref{eq:nicholas}) into a boundary integral equation, which is
amenable to ``almost linear'' complexity direct solvers using the HODLR matrix
techniques described in Section \ref{sec:HODLR}. In this section, we describe a
refinement on the HODLR format that will enable us to attain true $\cO(N)$ linear
complexity (as opposed to the~$\cO(N \log^2 N)$ complexity of HODLR
factorization) for a set of important problems in two dimensions. (In Section
\ref{sec:BIE_strong}, an additional refinement that is necessary to get linear
complexity in 3D will be described.)

The new concept that we introduce, ``nested basis matrices'', is designed to
reduce the overhead associated with storing and manipulating the basis matrices
associated with large blocks in the HODLR format. The idea is to use recursive
representations, so that the basis matrices on one level are expressed in terms
of the basis matrices on the next finer level, in a manner similar to how one
passes from one level to the next in multigrid, or when using wavelet based
methods. When this notion of nested bases is introduced to the HODLR format, we
obtain a more efficient format that is often referred to as
\textit{Hierarchically   Semi-Separable (HSS)} matrices in the literature, and
there is indeed a connection to the semi-separable matrices that we saw in
Section \ref{sec:onedimbvp}. However, there are distinct differences, and some
authors prefer the term \textit{Hierarchically Block Separable (HBS)}, which is
closely related to the notion of an ``$\mathcal{H}^{2}$-matrix with weak
admissibility''.

Our presentation starts with a single-level (non-hierarchical) method in
Section \ref{sec:singlelevel}, and then proceeds to a hierarchical one in
Section~\ref{sec:HBSinformal}. Next, Section~\ref{sec:HBSstability} discusses
numerical stability, and Section~\ref{sec:HBSbackground} provides connections to
the broader literature.

\subsection{Inversion of a block separable matrix} 
\label{sec:singlelevel}
Let us first consider an $N\times N$ matrix $\mtx{A}$ that is partitioned into $8\times 8$ blocks, each of size $n\times n$, so that $N = 8n$. We assume that the off-diagonal blocks all have ranks at most $k$, where $k < n$, and further assume that $\mtx{A}$ admits a factorization of the form
\begin{equation}
\label{eq:BS1}
\begin{array}{cccccccccc}
\mtx{A}         & = & \mtx{U}       & \tilde{\mtx{A}}       & \mtx{V}^{*}   & + & \mtx{D},\\
\includegraphics[scale=0.2]{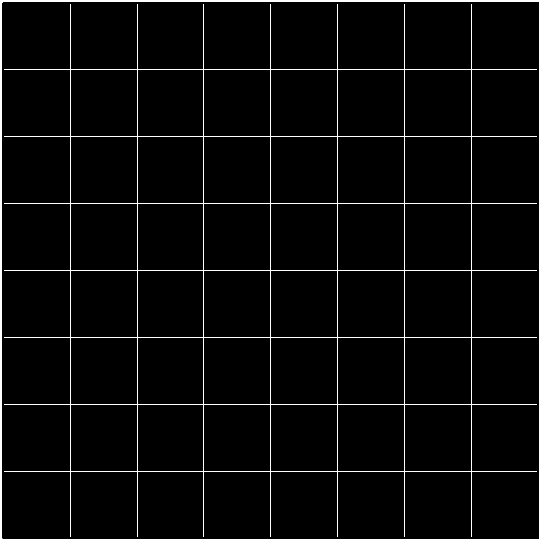} &&
\includegraphics[scale=0.2]{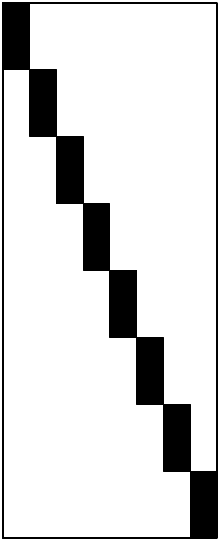} &
\raisebox{9mm}[0mm]{\includegraphics[scale=0.1]{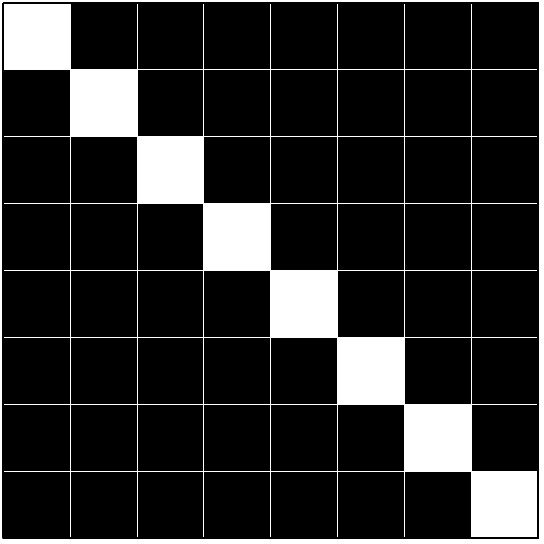}} &
\raisebox{9mm}[0mm]{\includegraphics[scale=0.2]{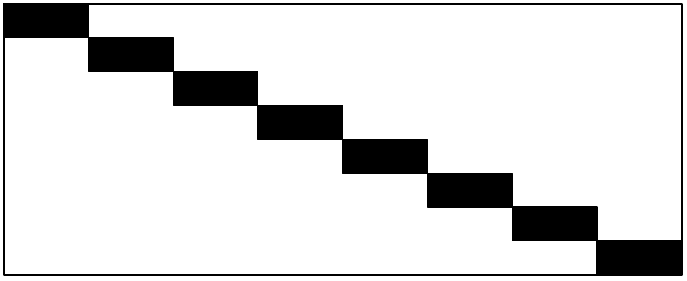}} &&
\includegraphics[scale=0.2]{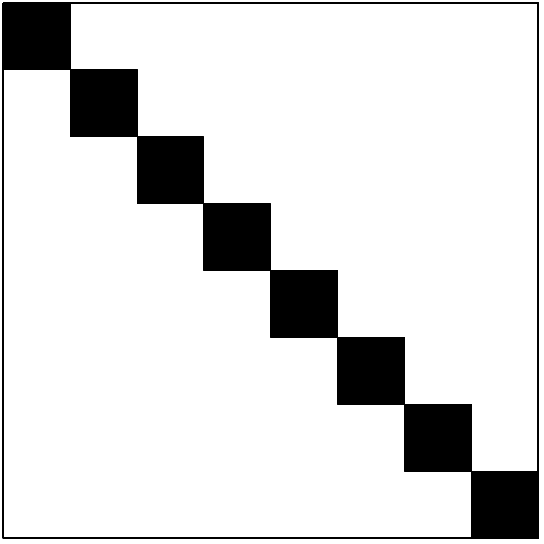}\\
8n \times 8n && 8n \times 8k & 8k \times 8k & 8k \times 8n && 8n \times 8n 
\end{array}
\end{equation}
To be precise, if the partitioning is defined via some index vectors~$I_{\tau}$
such that
$$
1:N = I_{8} \cup I_{9} \cup I_{10} \cup \cdots \cup I_{15},
$$
then our central assumption is that each off-diagonal
block~$\mtx{A}(I_{\alpha},I_{\beta})$ admits a factorization
\begin{equation}
\label{eq:defBS}
\begin{array}{cccccccc}
\mtx{A}(I_{\alpha},I_{\beta}) &=& \mtx{U}_{\alpha} & \tilde{\mtx{A}}_{\alpha,\beta} & \mtx{V}_{\beta}^{*},&\qquad\mbox{when}\ \alpha \neq \beta.\\
n \times n && n \times k & k\times k & k \times n
\end{array}
\end{equation}
The reason we start indexing at ``8'' is that we will soon organize the index
vectors into a hierarchical tree, as in Section~\ref{sec:trees}.  Observe in
particular that we use the same basis matrix $\mtx{U}_{\alpha}$ for the column
space of every off-diagonal block in the ``$\alpha$-row'', and likewise the same
basis matrix $\mtx{V}_{\beta}$ for the row space of every block in the
``$\beta$-column''. The block diagonal matrices $\mtx{U}$ and $\mtx{V}$ of
course hold these basis matrices, while $\mtx{D}$ holds the diagonal blocks. We
say that a matrix that obeys condition (\ref{eq:defBS}) is \textit{block
  separable}.

When a block separable matrix is invertible, its inverse is also block separable.
To be precise, $\mtx{A}^{-1}$ admits a factorization of the form
\begin{equation}
\label{eq:BS2}
\begin{array}{cccccccccc}
\mtx{A}^{-1} & = & \mtx{E} & \bigl(\tilde{\mtx{A}} + \hat{\mtx{D}}\bigr)^{-1} & \mtx{F}^{*}   & + & \mtx{G}.  \\
\includegraphics[scale=0.2]{08-nested/figs/fig_wood_A8.pdf} &&
\includegraphics[scale=0.2]{08-nested/figs/fig_wood_U8.pdf} &
\raisebox{9mm}[0mm]{\includegraphics[scale=0.1]{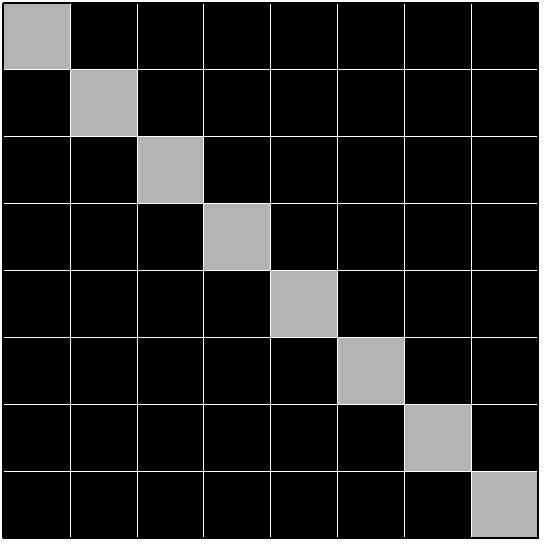}} &
\raisebox{9mm}[0mm]{\includegraphics[scale=0.2]{08-nested/figs/fig_wood_V8.pdf}} &&
\includegraphics[scale=0.2]{08-nested/figs/fig_wood_D8.pdf}\\
8n \times 8n && 8n \times 8k & 8k \times 8k & 8k \times 8n && 8n \times 8n 
\end{array}
\end{equation}
To make this claim precise, we observe that the decomposition of $\mtx{A}$ takes the form
$$
\mtx{A} = 
\underbrace{\mtx{U}\tilde{\mtx{A}}\mtx{V}^{*}}_{\textit{``Low rank''}}
+
\underbrace{\mtx{D}}_{\textit{``Easy to invert''}}.
$$
This brings the Woodbury formula for matrix inversion to mind, and indeed, there is a variation of this formula that provides exactly what we need:

\begin{lemma}[Variation of Woodbury]
\label{lemma:woodbury}
Suppose that $\mtx{A}$ is an invertible $N\times N$ matrix, that $K$ is a positive
integer such that $K < N$, and that $\mtx{A}$ admits the factorization
\begin{equation}
\label{eq:woodbury_assump}
\begin{array}{cccccccc}
\mtx{A} &=& \mtx{U} & \tilde{\mtx{A}} & \mtx{V}^{*} &+& \mtx{D}.\\
N \times N && N \times K & K \times K & K \times N && N \times N
\end{array}
\end{equation}
Then
\begin{equation}
\label{eq:woodbury}
\begin{array}{cccccccc}
\mtx{A}^{-1} &=& \mtx{E} & \bigl(\tilde{\mtx{A}} + \hat{\mtx{D}}\bigr)^{-1} & \mtx{F}^{*} &+& \mtx{G},\\
N \times N && N \times K & K \times K & K \times N && N \times N
\end{array}
\end{equation}
where
\begin{align}
\label{eq:def_muhD}
\hat{\mtx{D}} =&\ \bigl(\mtx{V}^{*}\mtx{D}^{-1}\mtx{U}\bigr)^{-1},\\
\label{eq:def_muE}
\mtx{E}       =&\ \mtx{D}^{-1}\mtx{U}\hat{\mtx{D}},\\
\label{eq:def_muF}
\mtx{F}       =&\ (\hat{\mtx{D}}\mtx{V}^{*}\mtx{D}^{-1})^{*},\\
\label{eq:def_muG}
\mtx{G}       =&\ \mtx{D}^{-1} - \mtx{D}^{-1}\mtx{U}\hat{\mtx{D}}\mtx{V}^{*}\mtx{D}^{-1},
\end{align}
provided all inverses that appear in the formulas exist.
Moreover, $\mbox{rank}(\mtx{G}) = N-K$.
\end{lemma}

Lemma \ref{lemma:woodbury} appears as Lemma 13.1 in \cite{2019_martinsson_fast_direct_solvers}, where a proof is provided.

Observe that the formulas (\ref{eq:def_muhD}) -- (\ref{eq:def_muG}) all involve block diagonal matrices, which means that they can be evaluated efficiently. The upshot is that we have reduced the task of inverting the $8n \times 8n$ matrix $\mtx{A}$ 
in (\ref{eq:BS1}) to, in essence, the task of inverting a matrix of size $8k \times 8k$.

Observe further that Lemma \ref{lemma:woodbury} does not guarantee the numerical
\textit{stability} of the inversion formula, which is a question we will return
to in Section \ref{sec:HBSstability}.


\subsection{A multilevel format}
\label{sec:HBSinformal}

The inversion formula described in Section \ref{sec:singlelevel} is a powerful
tool in its own right \cite{2005_martinsson_skel}, but to unleash its full
potential, we need to apply it recursively. To illustrate, let us continue to
discuss the~$8\times 8$ block matrix~$\mtx{A}$ that we introduced in
Section~\ref{sec:singlelevel}, but further add the assumption that~$\mtx{A}$ is
also a HODLR matrix with respect to the tree structure shown in
Figure~\ref{fig:tree}. To reveal the implications of the HODLR assumption, let
us relabel $\tilde{\mtx{A}} \rightarrow \mtx{A}^{(2)}$, so that (\ref{eq:BS1})
becomes
\begin{equation}
\label{eq:HBS1}
\begin{array}{cccccccccc}
\mtx{A} & = & \mtx{U}^{(3)} & \mtx{A}^{(2)} & (\mtx{V}^{(3)})^{*} & + & \mtx{D}^{(3)}.\\
\includegraphics[scale=0.2]{08-nested/figs/fig_wood_A8.pdf} &&
\includegraphics[scale=0.2]{08-nested/figs/fig_wood_U8.pdf} &
\raisebox{9mm}[0mm]{\includegraphics[scale=0.1]{08-nested/figs/fig_wood_AT8.pdf}} &
\raisebox{9mm}[0mm]{\includegraphics[scale=0.2]{08-nested/figs/fig_wood_V8.pdf}} &&
\includegraphics[scale=0.2]{08-nested/figs/fig_wood_D8.pdf}
\end{array}
\end{equation}
The fact that $\mtx{A}$ is a HODLR matrix now implies that there are additional rank-deficiencies in the off-diagonal blocks of $\mtx{A}^{(2)}$. To be precise, $\mtx{A}^{(2)}$ admits a factorization of the form
\begin{equation}
\label{eq:HBS2}
\begin{array}{cccccccccc}
\mtx{A}^{(2)} & = & \mtx{U}^{(2)} & \mtx{A}^{(1)} & (\mtx{V}^{(2)})^{*} & + & \mtx{D}^{(2)},\\
\includegraphics[scale=0.1]{08-nested/figs/fig_wood_AT8.pdf} &&
\includegraphics[scale=0.1]{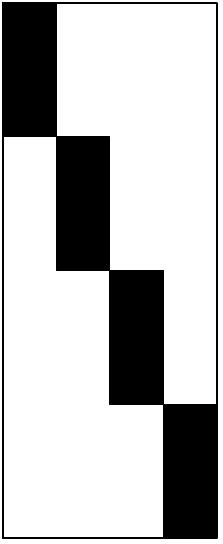} &
\raisebox{4.5mm}[0mm]{\includegraphics[scale=0.05]{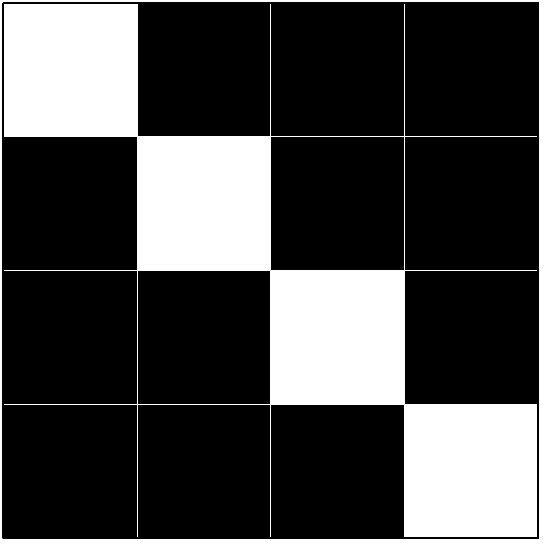}} &
\raisebox{4.5mm}[0mm]{\includegraphics[scale=0.1]{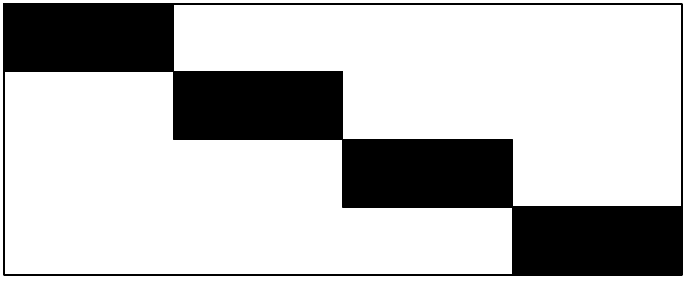}} &&
\includegraphics[scale=0.1]{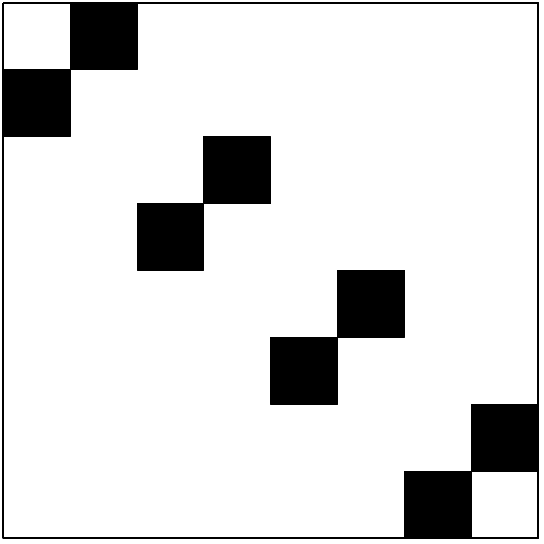}
\end{array}
\end{equation}
and similarly
\begin{equation}
\label{eq:HBS3}
\begin{array}{cccccccccc}
\mtx{A}^{(1)} & = & \mtx{U}^{(1)} & \mtx{A}^{(0)} & (\mtx{V}^{(1)})^{*} & + & \mtx{D}^{(1)}.\\
\includegraphics[scale=0.05]{08-nested/figs/fig_wood_AT4.pdf} &&
\includegraphics[scale=0.05]{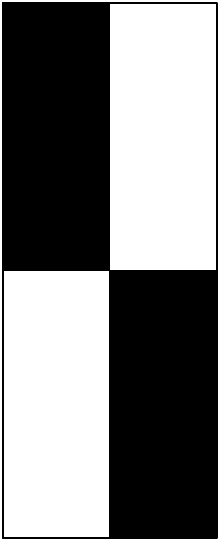} &
\raisebox{2.25mm}[0mm]{\includegraphics[scale=0.025]{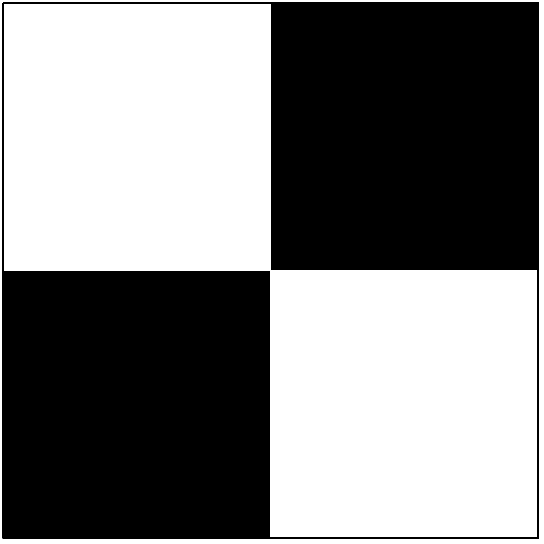}} &
\raisebox{2.25mm}[0mm]{\includegraphics[scale=0.05]{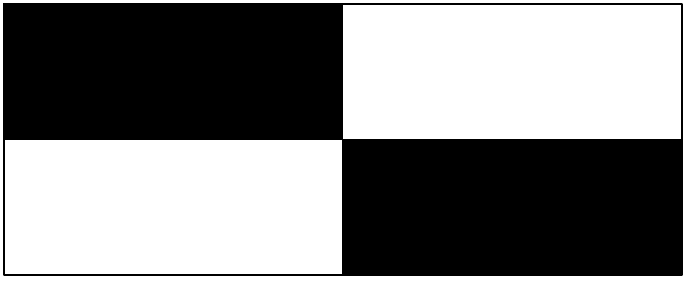}} &&
\includegraphics[scale=0.05]{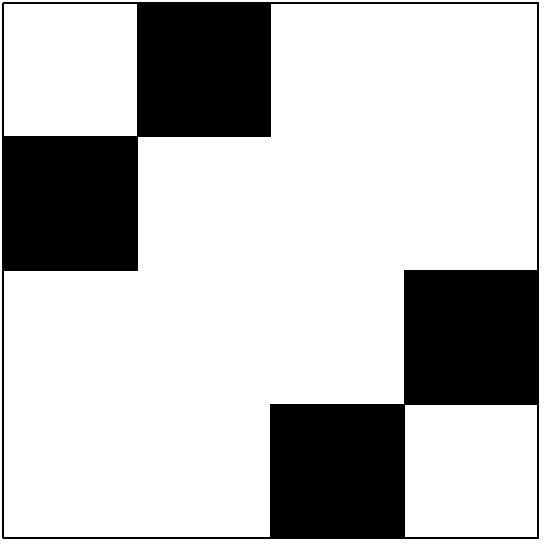}
\end{array}
\end{equation}
Combining (\ref{eq:HBS1}), (\ref{eq:HBS2}), and (\ref{eq:HBS3}) leads to a telescoping factorization of the form
\begin{equation}
\label{eq:HBS4}
\mA = \mU^{(3)}\bigl(\mU^{(2)}\bigl(\mU^{(1)}\,\mA^{(0)}\,(\mV^{(1)})^{*} + \mD^{(1)}\bigr)
(\mV^{(2)})^{*} + \mD^{(2)}\bigr)(\mV^{(3)})^{*} + \mD^{(3)},
\end{equation}
where the sparsity patterns of the factors are as follows:
\begin{center}
\raisebox{00.00mm}[20mm]{\includegraphics[scale=0.200]{08-nested/figs/fig_wood_U8.pdf}} \hspace{2mm} 
\raisebox{09.00mm}[20mm]{\includegraphics[scale=0.100]{08-nested/figs/fig_wood_U4.pdf}} \hspace{2mm} 
\raisebox{13.50mm}[20mm]{\includegraphics[scale=0.050]{08-nested/figs/fig_wood_U2.pdf}} \hspace{2mm} 
\raisebox{15.75mm}[20mm]{\includegraphics[scale=0.025]{08-nested/figs/fig_wood_AT2.pdf}}\hspace{2mm} 
\raisebox{15.75mm}[20mm]{\includegraphics[scale=0.050]{08-nested/figs/fig_wood_V2.pdf}} \hspace{2mm} 
\raisebox{13.50mm}[20mm]{\includegraphics[scale=0.050]{08-nested/figs/fig_wood_D2.pdf}} \hspace{2mm} 
\raisebox{13.50mm}[20mm]{\includegraphics[scale=0.100]{08-nested/figs/fig_wood_V4.pdf}} \hspace{2mm} 
\raisebox{09.00mm}[20mm]{\includegraphics[scale=0.100]{08-nested/figs/fig_wood_D4.pdf}} \hspace{2mm} 
\raisebox{09.00mm}[20mm]{\includegraphics[scale=0.200]{08-nested/figs/fig_wood_V8.pdf}} \hspace{2mm} 
\raisebox{00.00mm}[20mm]{\includegraphics[scale=0.200]{08-nested/figs/fig_wood_D8.pdf}} \hspace{2mm} 
\end{center}
We refer to a matrix admitting a representation of the form (\ref{eq:HBS4}) as a \textit{hierarchically block separable (HBS)} matrix.

The key point of this section is that HBS matrices are remarkably simple to invert. Before discussing inversion, though, let us pause to consider some implications of the telescoping factorization (\ref{eq:HBS4}). The first is that the blocks that need to be stored remain of the same size at every level.
Since the number of blocks to be stored gets cut in half at each level, we find
that the amount of memory required also gets halved each time we move to a
coarser level.  In consequence, the total storage cost can be expressed as a
geometric sum dominated by the leaves, and we find that the overall storage cost
is $\cO(kN)$, regardless of the number of levels in the tree.  This is in
contrast the $\cO(k N \log N)$ cost of storing a HODLR matrix, which was
incurred since the basis matrices get larger as we move towards the root of the
tree.  By an analogous argument, the cost of applying an HBS matrix to a vector
is also $\cO(kN)$.

Turning next to our main topic of inversion, recall from the previous section
that the bottleneck in evaluating the formula (\ref{eq:woodbury}) is the
inversion of the dense matrix $\tilde{\mtx{A}} + \hat{\mtx{D}}$ (since the
factors $\mtx{E}$, $\mtx{F}$, and $\mtx{G}$ all involve block diagonal
matrices). Relabeling the terms $\tilde{\mtx{A}} \rightarrow \mtx{A}^{(2)}$ and
$\hat{\mtx{D}} \rightarrow \hat{\mtx{D}}^{(2)}$, and exploiting formula
(\ref{eq:HBS2}), we find that this matrix takes the form
\begin{equation}
\label{eq:HBS5}
\begin{array}{cccccccccc}
\mtx{A}^{(2)} + \hat{\mtx{D}}^{(2)} & = & \mtx{U}^{(2)} & \tilde{\mtx{A}}^{(1)} & (\mtx{V}^{(2)})^{*}  & + & \mtx{D}^{(2)} + \hat{\mtx{D}}^{(2)},\\
\includegraphics[scale=0.1]{08-nested/figs/fig_wood_ATD8.pdf} &&
\includegraphics[scale=0.1]{08-nested/figs/fig_wood_U4.pdf} &
\raisebox{4.5mm}[0mm]{\includegraphics[scale=0.05]{08-nested/figs/fig_wood_AT4.pdf}} &
\raisebox{4.5mm}[0mm]{\includegraphics[scale=0.1]{08-nested/figs/fig_wood_V4.pdf}} &&
\includegraphics[scale=0.1]{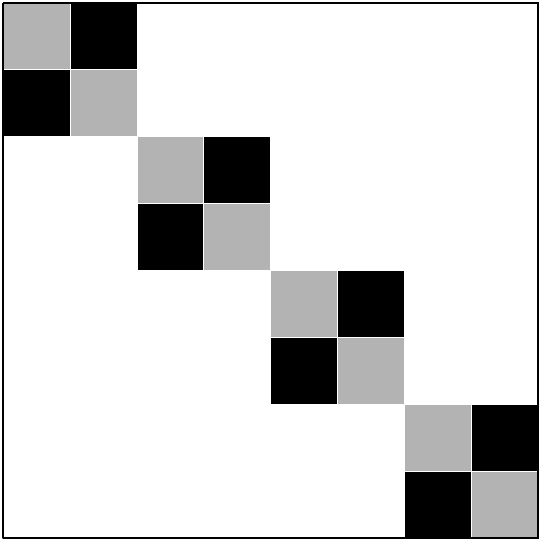}\\
8k \times 8k && 8k \times 4k & 4k \times 4k & 4k \times 8k && 8k \times 8k 
\end{array}
\end{equation}
which brings us to a point where we can simply apply Lemma \ref{lemma:woodbury} again!

Putting everything together, we end up with a remarkably simple algorithm for inverting an HBS matrix shown in Figure \ref{fig:HSSinversion}. Once the factors $\mtx{E}_{\tau},\,\mtx{F}_{\tau},\,\mtx{G}_{\tau}$ have all been constructed for every node $\tau$ in the tree, the inverse $\mtx{A}^{-1}$ can be applied to a vector via the algorithm shown in Figure \ref{fig:HSSinvapply}.
\begin{figure}
\begin{center}
\fbox{
\begin{minipage}{.9\textwidth}
\begin{tabbing}
\hspace{5mm} \= \hspace{5mm} \= \hspace{5mm} \= \kill
\textbf{loop} over all levels, finer to coarser, $\ell = L,\,L-1,\,\dots,1$\\
\> \textbf{loop} over all boxes $\tau$ on level $\ell$,\\
\> \> \textbf{if} $\tau$ is a leaf node\\
\> \> \> $\mtD_{\tau} = \mD_{\tau}$.\\
\> \> \textbf{else}\\
\> \> \> Let $\alpha$ and $\beta$ denote the children of $\tau$.\\
\> \> \> $\mtD_{\tau} = \mtwo{\mhD_{\alpha}}{\mtA_{\alpha,\beta}}{\mtA_{\beta,\alpha}}{\mhD_{\beta}}$.\\
\> \> \textbf{end if}\\
\> \> $\mhD_{\tau} = \bigl(\mV_{\tau}^{*}\,\mtD_{\tau}^{-1}\,\mU_{\tau}\bigr)^{-1}$.\\
\> \> $\mE_{\tau} = \mtD_{\tau}^{-1}\,\mU_{\tau}\,\mhD_{\tau}$.\\
\> \> $\mF_{\tau}^{*} = \mhD_{\tau}\,\mV_{\tau}^{*}\,\mtD_{\tau}^{-1}$.\\
\> \> $\mG_{\tau} = \mhD_{\tau} - \mtD_{\tau}^{-1}\,\mU_{\tau}\,\mhD_{\tau}\,\mV_{\tau}^{*}\,\mtD_{\tau}^{-1}$.\\
\> \textbf{end loop}\\
\textbf{end loop}\\
$\mG_{1} = \mtwo{\mhD_{2}}{\mtA_{2,3}}{\mtA_{3,2}}{\mhD_{3}}^{-1}$.
\end{tabbing}
\end{minipage}}
\end{center}
\caption{\label{fig:HSSinversion}Inversion of an HBS matrix $\mtx{A}$. The algorithm takes as its input all the blocks that together define the HBS representation (to be precise, the basis matrices $\mtx{U}_{\tau}$ and $\mtx{V}_{\tau}$ for every node $\tau$, the diagonal blocks $\mtx{D}_{\tau}$ for every leaf node $\tau$, and the ``sibling interaction matrices'' $\tilde{\mtx{A}}_{\alpha,\beta}$ and $\tilde{\mtx{A}}_{\beta,\alpha}$ for every sibling pair $\{\alpha,\beta\}$). The output is a collection of small matrices $\{\mtx{E}_{\tau},\mtx{F}_{\tau},\mtx{G}_{\tau}\}$ that together define the inverse, which can then be applied to vectors via the method in Figure \ref{fig:HSSinvapply}.}
\end{figure}

\begin{figure}
\begin{center}
\fbox{
\begin{minipage}{.9\textwidth}
\begin{tabbing}
\hspace{5mm} \= \hspace{5mm} \= \hspace{5mm} \= \kill
\textbf{loop} over all leaf boxes $\tau$\\
\> $\hat{\vct{u}}_{\tau} = \mF_{\tau}^{*}\,\vct{u}(I_{\tau})$.\\
\textbf{end loop}\\[1.25mm]
\textbf{loop} over all levels, finer to coarser, $\ell = L,\,L-1,\,\dots,1$\\
\> \textbf{loop} over all parent boxes $\tau$ on level $\ell$,\\
\> \> Let $\alpha$ and $\beta$ denote the children of $\tau$.\\
\> \> $\hat{\vct{u}}_{\tau} = \mF_{\tau}^{*}\,\vtwo{\hat{\vct{u}}_{\alpha}}{\hat{\vct{u}}_{\beta}}$.\\
\> \textbf{end loop}\\
\textbf{end loop}\\[1.25mm]
$\vtwo{\hat{\vct{q}}_{2}}{\hat{\vct{q}}_{3}} = \mhG_{1}\,\vtwo{\hat{\vct{u}}_{2}}{\hat{\vct{u}}_{3}}$.\\[1.25mm]
\textbf{loop} over all levels, coarser to finer, $\ell = 1,\,2,\,\dots,\,L-1$\\
\> \textbf{loop} over all parent boxes $\tau$ on level $\ell$\\
\> \> Let $\alpha$ and $\beta$ denote the children of $\tau$.\\
\> \> $\vtwo{\hat{\vct{q}}_{\alpha}}{\hat{\vct{q}}_{\beta}} =
       \mE_{\tau}\,\hat{\vct{u}}_{\tau} +
       \mG_{\tau}\,\vtwo{\hat{\vct{u}}_{\alpha}}{\hat{\vct{u}}_{\beta}}$.\\
\> \textbf{end loop}\\
\textbf{end loop}\\[1.25mm]
\textbf{loop} over all leaf boxes $\tau$\\
\> $\vct{q}(I_{\tau}) = \mE_{\tau}\,\hat{\vct{q}}_{\tau} + \mG_{\tau}\,\vct{u}(I_{\tau})$.\\
\textbf{end loop}
\end{tabbing}
\end{minipage}}
\end{center}
\caption{\label{fig:HSSinvapply} Application of the inverse of an HBS matrix. Given a vector $\vct{u}$, compute $\vct{q} = \mA^{-1}\,\vct{u}$ using the compact representation of $\mA^{-1}$ resulting from the algorithm in Figure \ref{fig:HSSinversion}.}
\end{figure}

\begin{remark}[Recursive skeletonization]
A particularly efficient subcategory of HBS matrices is obtained when the \textit{interpolatory decomposition} is used to factorize all the off-diagonal blocks. In this case, the matrix $\tilde{\mtx{A}}$ in (\ref{eq:BS1}) becomes a submatrix of $\mtx{A}$ itself, which means that these blocks never need to be computed explicitly. Additional savings are obtained in that the basis matrices $\mtx{U}_{\tau}$ and $\mtx{V}_{\tau}$ each contain $k\times k$ identity matrices \cite{2005_martinsson_skel,2012_greengard_ho_recursive_skeletonization,2012_martinsson_FDS_survey}. The resulting compression and inversion schemes are often referred to as \textit{recursive skeletonization}.
\end{remark}


\subsection{Numerical stability}
\label{sec:HBSstability}

The technique for inverting an HBS matrix that was presented in Sections
\ref{sec:singlelevel} and \ref{sec:HBSinformal} is simple to describe and to
implement, but it is not optimal either from the point of view of arithmetic
work or of numerical stability~\cite{amsel2025}.  Most importantly, it turns out
to be possible to avoid many of the matrix inversions in the intermediate steps
\cite[Ch.~18]{2019_martinsson_fast_direct_solvers}, which improves both speed
and accuracy.  Alternatively, one can form the rank structured equivalent of a
``ULV'' factorization where a matrix is brought to triangular form through a
sequence of unitary
transformations~\cite{2006_chandrasekaran_HSS_ULV}.  This
method is also of true linear complexity, is ideal from the point of view of
numerical stability, and is (like the other methods described) amenable to
highly optimized implementations \cite{ma2024inherently}.  Other
reformulations that often improve conditioning in practice are described in
\cite{gopal2022accelerated, 2012_greengard_ho_recursive_skeletonization}.

In terms of rigorous theory, the stability of the various factorization 
schemes is poorly understood. Some results exist for block low-rank
formats~\cite{higham2022solving}, but less is known about the hierarchical
setting~\cite{2016_xia_stability}.  There is in general no guarantee that the
intermediate matrices that need to be inverted are well-conditioned, even when
$\mtx{A}$ itself is well conditioned. One exception is the case where $\mtx{A}$
is symmetric and positive definite, in which case stability can often be proven, cf.~\cite[Sec.~3]{2012_martinsson_FDS_survey} and \cite{2016_xia_stability}.

\subsection{Background}
\label{sec:HBSbackground}

The ideas underpinning the direct solvers described in this section trace back
at least to
\cite{starr_rokhlin,1991_greengard_rokhlin_twopoint,1996_mich_elongated}, with
the specific version described here first published in
\cite{2005_martinsson_skel} for the single level scheme, and
\cite{2005_martinsson_fastdirect} for the multi level scheme.  These works
approached the question specifically from the point of view of fast solvers for
integral equations.  Later refinements of these solvers include
\cite{2012_martinsson_FDS_survey,2012_greengard_ho_recursive_skeletonization,
  2013_martinsson_smooth_BIE,2014_corona_martinsson_BIE_plane,2009_martinsson_ACTA}.
Similar ideas have been known in the linear algebra community as well, with
influential early works including \cite{2010_gu_xia_HSS,2005_gu_HSS}.  A
discussion of the relationship between HODLR and HSS can be found in
\cite{2016_kressner_rankstructured_review}.

The HSS/HBS structure can also be viewed as a simplified version of the Fast Multipole Method where interactions between \textit{all} boxes are compressed. In this view, the translation table between the two formalisms becomes:

\vspace{1mm}

\begin{center}
\begin{tabular}{l|c|l}
& \textit{HBS object} & \textit{FMM equivalent} \\ \hline
For a parent node $\tau$: & $\mtx{V}_{\tau}$ & Mpole-to-mpole \textit{and} sources-to-mpole.\\
For a parent node $\tau$: & $\mtx{U}_{\tau}$ & Local-to-local \textit{and} local-to-potentials.\\
For a sibling pair $\{\alpha,\beta\}$: & $\tilde{\mtx{A}}_{\alpha,\beta}$ & Mpole-to-local.\\
For a leaf node $\tau$: & $\mtx{D}_{\tau}$ & Sources-to-potential.
\end{tabular}
\end{center}

\vspace{1mm}

Finally, the HBS/HSS formalism fits naturally into the taxonomy of hierarchical matrices, where they would be classified as ``$\mathcal{H}^{2}$-matrices with weak admissibility''. (We note that in the early literature on $\mathcal{H}$-matrices, the term ``weak admissibility'' had a slightly different meaning. These papers would refer to what we call ``weak admissibility'' as ``off-diagonal'' admissibility.) We defer a more detailed discussion of hierarchical matrices until Section \ref{sec:strong_discussion}.



%% file: 09-strong/strong_new02.tex
In this section, we introduce the notion of \textit{strong admissibility}, which
is the final piece of machinery that is required to obtain linear complexity
algorithms for general problems in three dimensions.  The idea of strong
admissibility is most easily discussed in the context of compressing a dense
matrix arising from the discretization of an integral equation, such as those
previously discussed in Section~\ref{sec:BIE2D}.  The motivation here is
that if the formats we used in Sections~\ref{sec:HODLR} and~\ref{sec:BIE_nested}
(which are based on \textit{weak admissibility}) were to be applied to a problem
involving a general distribution of points in three dimensions, then the
numerical ranks of the off-diagonal blocks would grow large. 
To fight back against this effect, we will treat more
blocks in the matrix as unstructured (dense). The price to pay is that we end
up with additional data structures, more dense blocks to keep track of, and
substantially more complicated algorithms.

\subsection{Strong versus weak admissibility}

\begin{figure}[t]
    \centering
    \includegraphics[width=115mm]{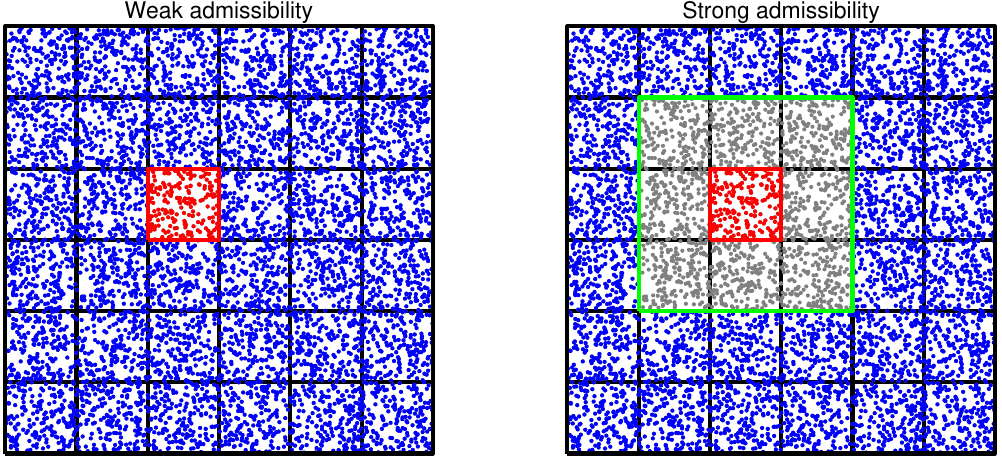}\\
    (a) \hspace{60mm} (b)
    \caption{A square domain $\Omega$ holding a set of
    point sources is divided into $6\times 6$ subdomains. One of them,
    $\Omega_{\tau}$, is highlighted in red. (a) In weak admissibility,
    the interactions between the points inside~$\Omega_{\tau}$ (red) and \textit{all}
    external points (blue) are compressed. (b) In strong admissibility,
    only the interactions with the points that are \textit{well-separated},
    cf.~Remark \ref{remark:wellsep}, are compressed. The well-separated points
    (blue) are collected in the index vector $I_{\tau}^{\rm far}$.}
    \label{fig:strongweak1}
\end{figure}

The difference between weak and strong admissibility comes down to which
subblocks of the global matrix are represented through low rank approximations.
To illustrate, let us consider a square domain $\Omega$ that has been
partitioned into $6\times 6$ subblocks, as shown in
Figure~\ref{fig:strongweak1}. We can imagine that the interaction between each
pair of points -- i.e. each matrix entry -- is given by a Green's function for
an elliptic PDE.

In weak admissibility, we compress the interactions between all pairs of
non-identical boxes, as shown in Figure~\ref{fig:strongweak1}(a). In
strong admissibility, we compress only those interactions corresponding to boxes
that do not touch, as shown in Figure~\ref{fig:strongweak1}(b).

\begin{figure}[t]
    \centering
    \includegraphics[width=115mm]{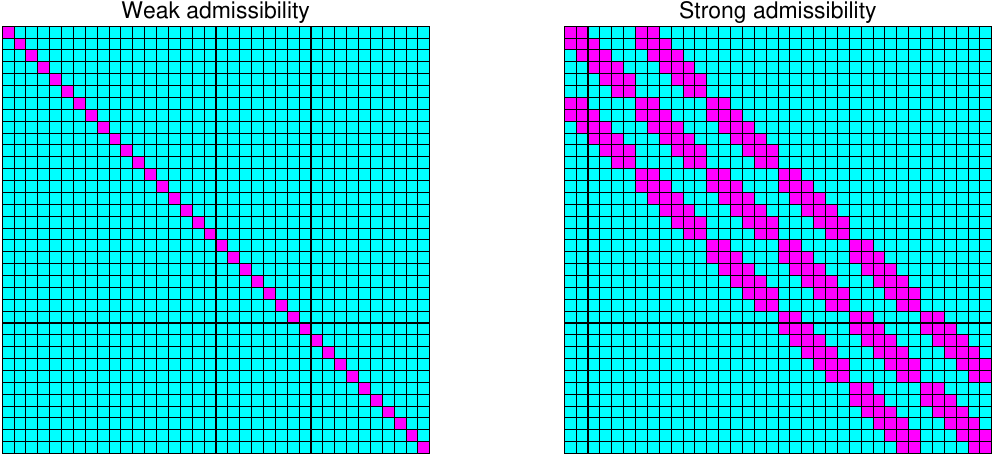}

    (a) \hspace{60mm} (b) \caption{The matrix associated with the domain shown
    in Figure \ref{fig:strongweak1} is split into $36\times 36$ blocks. The cyan
    blocks are represented as low rank. Observe how strong admissibility leads
    to far more dense blocks.} \label{fig:strongweak2} \end{figure}

Let us next consider the matrix representing interactions between the $36$
boxes, as shown in Figure \ref{fig:strongweak2}. Weak admissibility is to the
left, and strong to the right, with all low-rank blocks in cyan. We see that
weak admissibility leads to far more blocks being represented in a low rank
format, and a simple tessellation pattern. In contrast, when compression is
based on a strong admissibility criterion, up to nine blocks per row must be
represented as dense.

\begin{figure}[t]
    \centering
    \includegraphics[width=115mm]{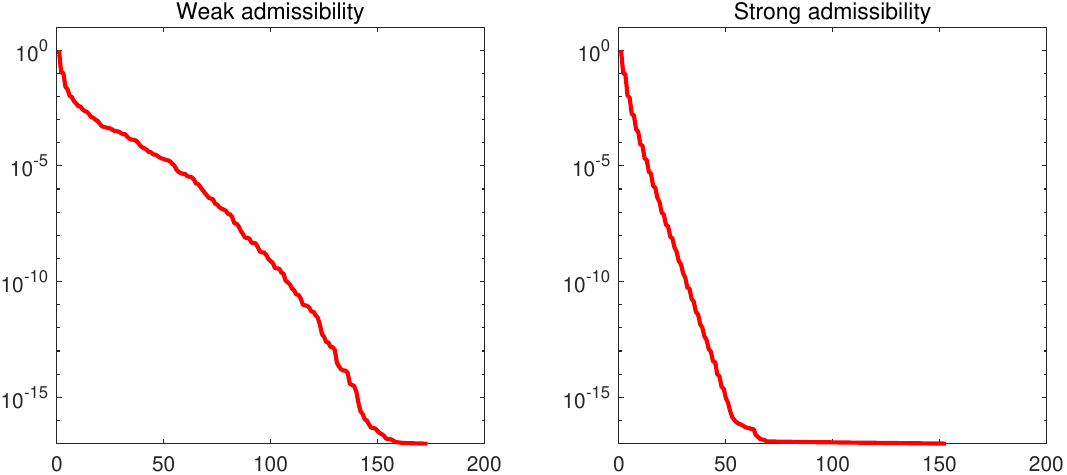}

    (a) \hspace{60mm} (b)
    \caption{The singular values of the matrix that needs to be compressed to
    compute the basis vectors associated with the box $\tau$ marked in
    Figure \ref{fig:strongweak1} (when a ``uniform basis'' is used).
    (a) The normalized singular values of the matrix
    $\mtx{A}(I_{\tau},I_{\tau}^{\rm c})$.
    (b) The normalized singular values of the matrix
    $\mtx{A}(I_{\tau},I_{\tau}^{\rm far})$.
    Observe how strong admissibility leads to lower ranks.}
    \label{fig:strongweak3}
\end{figure}

The advantage of strong admissibility is that the ranks of interaction are far
smaller, as illustrated in Figure \ref{fig:strongweak3}.  To be precise, the
figure on the left shows the singular values of the matrix
$\mtx{A}(I_{\tau},I_{\tau}^{\rm c})$, where $I_{\tau}$ is the index vector
identifying the block $\Omega_{\tau}$ in Figure \ref{fig:strongweak1}, and where
$I_{\tau}^{\rm c} = (1:N) \backslash I_{\tau}$ is the \textit{complement
  vector}.  Indeed,~$\mtx{A}(I_{\tau},I_{\tau}^{\rm c})$ is the matrix that we
need to numerically compress in order to build the matrix $\mtx{U}_{\tau}$ in
the factorization (\ref{eq:defBS}) in Section~\ref{sec:singlelevel}.  In
contrast, the figure on the right shows the singular values of the matrix
$\mtx{A}(I_{\tau},I_{\tau}^{\rm far})$, where $I_{\tau}^{\rm far}$ identifies
all points that are \textit{well separated} from the box $\Omega_{\tau}$,
cf.~Remark~\ref{remark:wellsep}.  The
matrix~$\mtx{A}(I_{\tau},I_{\tau}^{\rm far})$ corresponds to the interaction
between the red box and the blue boxes in the figure on the right in
Figure~\ref{fig:strongweak1}.  We see that strong admissibility leads to
significantly faster decay of the singular values, and consequently to much
lower ranks for any given precision~$\varepsilon$. Furthermore, when strong
admissibility is used, there is a universal bound on the numerical rank that is
\textit{independent} of the number of points in any of the boxes.  In contrast,
when weak admissibility is used, the numerical rank slowly grows as more and
more points are added due to the increased number of interacting points which
are actually nearby one another, see Remark \ref{remark:rankgrowth} for details.

\begin{remark} \label{remark:wellsep} In the FMM, the technical
condition that is used to demarcate which blocks should be treated as low rank
is based on the notion of points and boxes being \textit{well separated}. If $\tau$
is a box with side length~$a$ centered at a point~$\pcc_{\tau}$, then
we say that any other box (or point) is well-separated from~$\tau$ if it does
not intersect the box of sidelength $3a$ centered at $\pcc_{\tau}$. This
larger box is shown in green in Figure \ref{fig:strongweak1}.
\end{remark}

\begin{figure}[t] \centering
  \includegraphics[width=.75\textwidth]{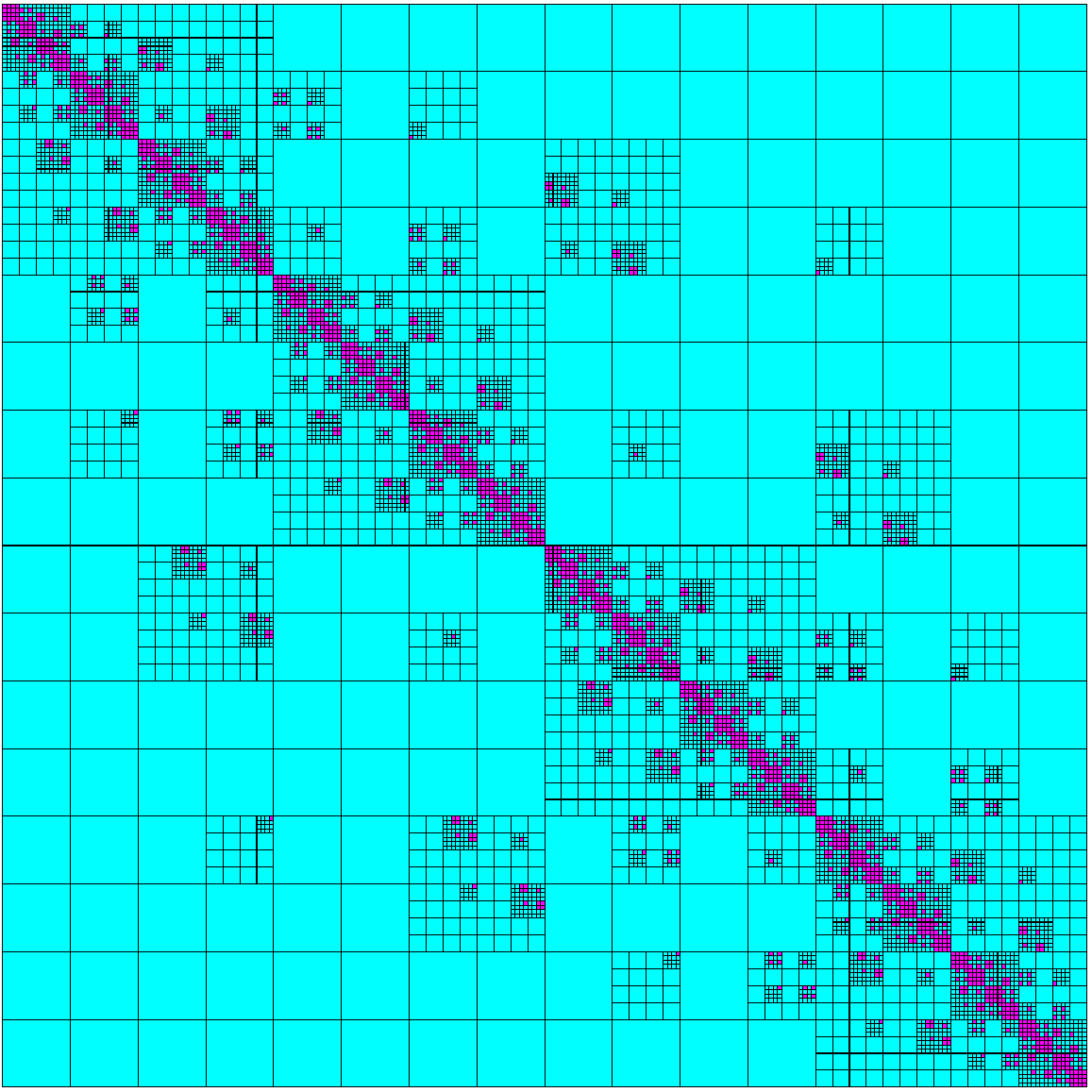}
  \caption{Tessellation pattern of a rank structured matrix when strong
  admissibility is used. The underlying geometry is a square that is
  partitioned into $16 \times 16$ leaf boxes. Cyan blocks are low rank, and
  magenta block are represented as dense.} \label{fig:fmm_tessellation}
\end{figure}

So far, our discussion has focused on the \textit{single level} or \textit{flat}
environment, where the matrix is tessellated into more or less equi-sized
submatrices. The concept directly carries over to hierarchical formats, with
only blocks that are well-separated being compressed at each level.
Figure~\ref{fig:fmm_tessellation} illustrates a typical tessellation pattern for
a simple geometry in two dimensions when strong admissibility is used.

To summarize:
\begin{itemize}
\item \textit{Weak admissibility:} All off-diagonal blocks are compressed. This leads to simple tessellation patterns and more blocks represented in low rank format, but higher numerical ranks.
\item \textit{Strong admissibility:} Only blocks that represent boxes that are well separated are compressed. This leads to lower ranks, but far more dense blocks, more complex tessellation patterns, and more complex factorization algorithms.
\end{itemize}

\begin{remark}[Rank growth in weak admissibility]
\label{remark:rankgrowth}
It is relatively well understood how rapidly the numerical ranks grow if weak
admissibility is used in higher dimensions.  For instance, if a set of
collocation points is distributed relatively uniformly in a two dimensional
domain (as in Figure \ref{fig:strongweak1}, or if a uniform grid is used), then
a box that holds $m$ points typically has numerical interaction rank that scales
as $\mathcal{O}(m^{1/2})$, resulting in overall $\mathcal{O}(N^{3/2})$
complexity for standard direct solvers, see~\cite{gopal2022accelerated} for one
example.  Similar scalings are also observed for boundary integral equations on
regular surfaces in three dimensions.  The worst case for weak scaling would be
a volume integral equation in 3D, where the interaction rank would scale as
$\mathcal{O}(m^{2/3})$, resulting in a $\mathcal{O}(N^2)$ overall
complexity.
\end{remark}

\subsection{The matrix vector multiplication}

The strong admissibility condition has been used in fast summation schemes for a
long time, forming an integral part of the Fast Multipole Method (FMM)
of Greengard and Rokhlin~\cite{rokhlin1987,rokhlin1985} of the 1980s, as well as
in the contemporaneous Barnes-Hut and tree-code
algoritms~\cite{1986_barnes_hut}. The data structures required are relatively
straight-forward in the case of uniform trees. For non-uniform point
distributions that rely on adaptive trees, some additional machinery is required
\cite{1988_CarrierGreengardRokhlin}, but no major complications arise. The
classical Barnes-Hut algorithm does not fully implement the notion of nested
bases (cf.~Section \ref{sec:BIE_nested}) and attains a computational complexity
of~$O(kN\log N)$ for the matrix vector multiplication. In contrast, the FMM does
rely on nested bases and incurs a computational complexity~$O(kN)$ with rigorous
accuracy guarantees.

The extension of techniques such as the FMM and Barnes-Hut to a more general
linear algebraic framework was initiated by Hackbusch through the $\mathcal{H}$-
and $\mathcal{H}^{2}$-matrix formalisms \cite{hackbusch,2002_hackbusch_H2}. The
distinction between the two formats concerns whether nested basis matrices are
used, see Section~\ref{sec:BIE_nested}.  The~$\mathcal{H}$-matrix format is
similar to Barnes-Hut in that nested bases are not used, and the associated
algorithms achieve an~$O(kN\log(N))$ complexity. The~$\mathcal{H}^{2}$-matrix
formalism generalizes the FMM and attains~$O(kN)$ complexity for the
matrix-vector multiplication through the use of both upward and downward passes
through the hierarchy, as well as the use of nested bases.

\subsection{Matrix factorization and inversion} 

Rank structured matrices associated with weak admissibility conditions often
admit exact inversion formulas that are well suited to high performance
implementations; we saw this illustrated for the HODLR format in
Section~\ref{sec:HODLRinv} and for HBS matrices in Section~\ref{sec:HBSinformal}.
In contrast, existing algorithms for inverting and factorizing matrices
associated with strong admissibility typically involve more complex recursive
algorithms where the ranks of the off-diagonal blocks successively grow as the
algorithm proceeds \cite{2010_borm_book,2016_kressner_rankstructured_review}
(analogous to the recursive method for inverting a HODLR matrix in
Section \ref{sec:HODLRwoodbury}). However, it so happens that most of the
matrices that arise in applications are highly amenable to rank structured
representations, so in practice, recompression can often successfully combat the
rank growth.

Several papers that describe non-recursive algorithms for factorization and
inversion have been published. 
In~\cite{2016_ying_HIF_IE}, an algorithm was proposed that used two data
structures, one a geometrically shifted version of the other. Matrices are
compressed with respect to one data structure, and then re-compressed with
respect to the other. This effectively shifts discretization points from the
near field to the far field, achieving additional compression.
In~\cite{2014_corona_martinsson_BIE_plane}, a scheme was proposed for inverting
volume integral equations in two dimensions by re-compressing certain matrices
in a hierarchical format, not low-rank format -- this results in somewhat of a
``double hierarchy'' data structure, but yields a linearly scaling algorithm for
constructing the inverse.

Perhaps the most promising method for attaining linear complexity inversion
while avoiding either recursions or double hierarchies is the \textit{strong
  skeletonization} technique~\cite{2017_ho_ying_strong_RS}, which is related to
the~\textit{Inverse FMM}~\cite{2014_darve_IFMM, coulier2017inverse}.
Both of these methods compute triangular factorizations of a given matrix
through a sequence of block Gaussian eliminations. They work within the data
structures provided by the Fast Multipole Method, albeit with extended
interaction lists \cite{sushnikova2023fmm}. These methods have high storage
requirements, but appear to be easier to implement than many competing
techniques and well-suited for parallelization. However, they do still rely on
recompression of modified blocks that arise during the computation and do not
work on ``general''~$\mathcal{H}^{2}$-matrices -- they are geared toward
problems arising from elliptic PDEs. Recently, \cite{yesypenko2023randomized}
proposed a technique based on randomized compression that reduces both
bookkeeping and storage requirements in the method of
\cite{2017_ho_ying_strong_RS}.

Usually, since the compressibility properties of an integral operator are based
on geometric considerations, spatial partitioning is best done using quadtrees
in two dimensions and octrees in three dimensions. These trees can be
constructed in both uniform and adaptive fashions, depending on the particular
distribution of discretization points. Algorithms based on both weakly
admissible compression schemes and strongly admissible compression schemes are
compatible with these data structures, but of course the implementations and
implicit constants in the algorithms are different in each case.

\subsection{Relationships between different rank structured formats}
\label{sec:strong_discussion}

Having described the concepts of nested basis functions in
Section~\ref{sec:BIE_nested} and strong vs.~weak admissibility in
Section~\ref{sec:BIE_strong}, we can now organize the different rank structured
formats that we have described into a table:

\vspace{2mm}

\begin{center}
\begin{tabular}{p{39mm}|p{38mm}|p{40mm}}
& \textit{Weak admissibility.} & \textit{Strong admissibility.}\\
& Simpler data structures.     & Lower ranks. \\
& Good for 1D and 2D.          & Longer interaction lists. \\ \hline
\textit{General basis matrices:}
& HODLR
& $\mathcal{H}$ matrices. \\
Lower ranks. &&  Barnes-Hut. \\ \hline
\textit{Nested basis matrices:} & HBS / HSS. & $\mathcal{H}^{2}$ matrices. \\
Enable $O(N)$ complexity. && FMM and its descendants.
\end{tabular}
\end{center}

\vspace{2mm}

\noindent
Technically speaking, \textit{all} formats in the above table
are~$\mathcal{H}$-matrices, but when the term is used it is generally implied
that a strong admissibility condition is enforced.

The full hierarchical formats are sometimes unnecessarily complex. In practice,
flat tessellations, as illustrated in Figures~\ref{fig:strongweak1}(a)
and~\ref{fig:strongweak2}(b), are often preferred, using either strong or weak
admissibility. These formats are sometimes referred to as \textit{block low
  rank}~\cite{2017_amestoy_BLR_sparse_direct, 2019_amestoy_BLR_multifrontal},
\textit{tile low rank}~\cite{2024_keyes_tile_low_rank}, \textit{block
  separable}~\cite{2012_martinsson_FDS_survey}, or
\textit{skeletonization}~\cite{2005_martinsson_skel}.

Finally, let us note that the concept of strong admissibility is 
central to the Fast Multipole Method
\cite{rokhlin1985,rokhlin1987,1988_CarrierGreengardRokhlin,1988_greengard_dissertation},
which forms the intellectual foundation of Fast Direct Solvers. It was
also an integral part of the Barnes-Hut method~\cite{1986_barnes_hut}.
Additional discussions of strong versus weak admissibility can be found in~\cite{ashcraft2021,amsel2025}.


%% file: operator_algebra3.tex
The development of fast direct solvers has unlocked the ability to explicitly
compute functions of operators, and to directly solve coupled linear systems
involving different operators or different discretizations. The ability to
perform operator algebra also opens up a whole host of opportunities in terms of
developing analysis-based preconditioning techniques for a wide variety of
problems.  To illustrate, consider the large family of problems where the
solution involves applying a function $f = f(\cA)$ of some elliptic
operator~$\cA$.  Examples include, among many others,
\begin{itemize}
\item evolution operators for parabolic PDEs where $f(\cA) = \exp(-t\cA)$ and
  hyperbolic PDEs where $f(\cA) = \exp(it\sqrt{\cA})$,
\item Cayley transforms with $f(\cA) = (I + i\cA)^{-1}(I - i\cA)$, and
\item the solution of fractional PDEs where $f(\cA) = \cA^{-\alpha}$
  for $\alpha \in (0,1)$).
\end{itemize}

Now, since we have established techniques for inverting elliptic operators, we
can also \textit{explicitly} form an approximation to~$f(\cA)$ through a
rational approximation.  In the simplest setting, suppose that~$\cA$ is a
diagonalizable operator and that we can build a rational approximation to~$f$:
\begin{equation}
  \label{eq:rationalapprox}
  f(t) \approx \sum_{j=1}^{J} \frac{a_{j}}{1 - b_{j}t},\qquad t \in \Lambda,
\end{equation}
where $\Lambda$ holds the ``relevant'' part of the spectrum of $\cA$ (what this
means specifically is highly problem dependent).  Then, we
can explicitly form an approximation
\[
  f(\cA) \approx \sum_{j=1}^{J} a_{j}\,\left( \cI - b_{j}\cA \right)^{-1}.
\]
See~\cite{2008_higham_matrix_functions} for a discussion on the
convergence of matrix functions and related topics, and~\cite{cortinovis2022dnc,casulli2024symmetric}
for a discussion of some algorithms for computing these
matrix functions.  Next, let~$\mtx{A}$ denote the matrix obtained from a discretization
of the operator~$\cA$.  When~$\mtx{A}$ is a rank structured matrix, fast direct
solver techniques often allow us to efficiently compute each
term~$\bigl(\mtx{I} - b_{j}\mtx{A}\bigr)^{-1}$~\cite{2014_haut_hyperbolic}.  In
practice, both the derivation of rational approximations and an analysis of the
errors is done through the study of Cauchy integrals
\cite[Ch.~14]{hackbusch2015hierarchical}.  

\begin{remark}
If an approximation such as (\ref{eq:rationalapprox}) involves a large
number of terms, it could become expensive to construct all 
the inverses $\{\bigl(\mtx{I} - b_{j}\mtx{A}\bigr)^{-1}\}_{j=1}^{J}$,
and infeasible to store them all given memory constraints. 
However, it is often possible to exploit the fact that the matrices that
are inverted are identical except for a diagonal additive term to efficiently
compute a whole family of inverses jointly \cite{2023_martinsson_gopal_broadband}.
\end{remark}

\begin{figure}[t!]
\begin{center}
\textit{The geometry:}
\hspace{32mm}
\textit{The linear system:}\\
\includegraphics[height=30mm]{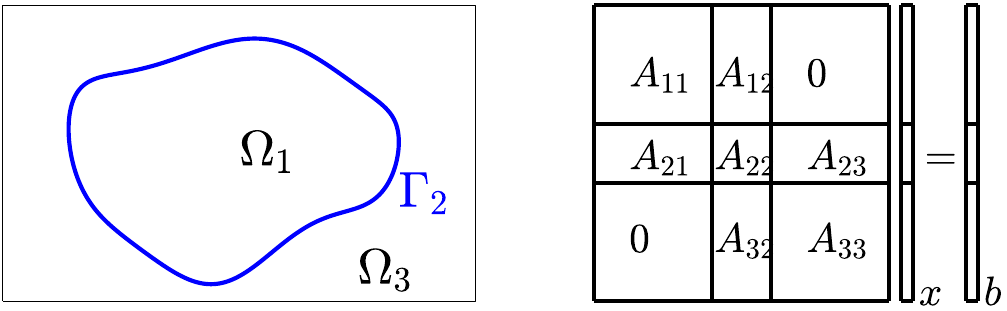}
\end{center}
\caption{A coupled problem, as described in Section \ref{sec:operator_algebra}.}
\label{fig:twodomain}
\end{figure}

To illustrate the power of being able to perform operator algebra directly, let
us next consider a toy problem involving two adjacent and interacting domains,
e.g., a wave propagation problem involving both a fluid and a solid
medium. Suppose the two domains~$\Omega_{1}$ and~$\Omega_{3}$ are connected
through a boundary $\Gamma_{2}$, as shown in Figure~\ref{fig:twodomain}. Upon
discretization, we would typically obtain a linear system of the form
\begin{equation}
\label{eq:coupledsystem}
\left[\begin{array}{rrr}
\mtx{A}_{11} & \mtx{A}_{12} & \mtx{0} \\
\mtx{A}_{21} & \mtx{A}_{22} & \mtx{A}_{23} \\
\mtx{0}      & \mtx{A}_{32} & \mtx{A}_{33}
\end{array}\right]
\left[\begin{array}{r}
\mtx{x}_{1} \\ \mtx{x}_{2} \\ \mtx{x}_{3}
\end{array}\right]
=
\left[\begin{array}{r}
\mtx{b}_{1} \\ \mtx{b}_{2} \\ \mtx{b}_{3}
\end{array}\right],
\end{equation}
where $\mtx{A}_{11}$ and $\mtx{A}_{33}$ are discretizations of the equilibrium
operators for the subdomains~$\Omega_{1}$ and~$\Omega_{3}$.  We seek to solve
the full system that incorporates the interface conditions along $\Gamma_{2}$.
Coupled systems such as (\ref{eq:coupledsystem}) often pose significant
difficulties to iterative methods, even in situations where good individual
preconditioners for~$\mtx{A}_{11}$ and~$\mtx{A}_{33}$ are available.  With the
ability to execute operator algebra directly, we can side-step these problems by
explicitly forming a new linear system that is defined on the joint interface
alone
\begin{equation}
\label{eq:interfaceeq}
\bigl(\mtx{A}_{22} - \mtx{A}_{21}\mtx{A}_{11}^{-1}\mtx{A}_{12}
                   - \mtx{A}_{23}\mtx{A}_{33}^{-1}\mtx{A}_{32}\bigr)\vct{x}_{2}
=
\vct{b}_{2} - \mtx{A}_{21}\mtx{A}_{11}^{-1}\vct{b}_{1} - \mtx{A}_{23}\mtx{A}_{33}^{-1}\vct{b}_{3},
\end{equation}
and then directly factorize the coefficient matrix. 
Similar ideas can be applied to multiphysics problems 
in a single domain~$\Omega$, as well as to problems 
which are time dependent and require implicit
time integration methods for stability reasons.

Lastly, we note that in practice, the methods discussed in this survey are often
executed at lower precisions, with the resulting approximate inverses deployed
as preconditioners alongside efficient matrix-vector multiplication
algorithms. When solving a system only once, these hybrid approaches can
outperform standalone direct solvers, since the computational cost of FDS
increases sharply with the requested accuracy. Preconditioning remains valuable
even in the ``multiple right-hand sides'' environment, as memory limitations may
prohibit storing a high-precision matrix inverse.

Techniques for using approximate LU decompositions of $\mathcal{H}$-matrices as
preconditioners are explored in
\cite{2007_leborne_HLU,2006_leborne_H_matrix_precond,2008_grasedyck_HLU}, see
also \cite[Ch.~13]{2016_hackbusch_iterative_book}.  Analogous techniques in the
context of HSS-accelerated multifrontal methods are described in
\cite{2009_xia_superfast,2012_xia_robust_efficient_multifrontal}, with a related
approach in \cite{2017_pouransari_darve}.  Rather than using a low computational
tolerance, one can halt a high- or medium-precision factorization process before
all nodes in a hierarchical tree have been processed (in a manner somewhat
analogous to incomplete LU for sparse direct solvers).  This is particularly
effective when used in combination with weak admissibility, as the rank growth
that eventually dooms weak admissibility for large scale problems in 3D is
limited to the nodes in the tree that are closer to the root
\cite{2013_martinsson_smooth_BIE,gopal2022accelerated}.  For high frequency
problems, rank growth presents a challenge in the general case. However, for
geometries which can be tesselated into thin slabs, techniques based on
``sweeping'' have proven successful in many
environments~\cite{engquist2011sweeping,2020_demanet_L_sweeps,2019_gander_iterative_survey}.
In the context of integral equations, ill-conditioning can sometimes be
introduced by local refinement of corners and edges; this issue can be resolved
be exploiting numerical rank-deficiencies in combination with local direct
inversion
\cite{2013_martinsson_smooth_BIE,2013_martinsson_corner_short,2008_helsing_corner_BIE,2013_helsing_RCIP_tutorial}
(this problem can also be avoided by using specialized quadratures such
as~\cite{2012_bremer_nystrom, 2012_bremer_direct_BIE_corners, 2019hoskins,
  serkh2016solution}).

%% file: 10-dissipation/dissipation_new03.tex
A recurring theme of this survey has been the claim that many dense matrices
representing the solution operators for elliptic PDEs are rank structured. The
reason that off-diagonal blocks of these matrices have low numerical rank is
that even though they represent global operators in which every point is
connected to every other point, the amount of \textit{information} that is
communicated from one subdomain to another is very limited.
In this section, we will attempt to quantify this statement by illustrating what
the typical numerical ranks actually are.
Section~\ref{sec:dissip_fast_summation} discusses a particularly simple
situation where we can prove the boundedness of the ranks directly from known
analytic expansions of the fundamental solution of the elliptic PDE. Then, we
will demonstrate through empirical examples that the same type of information
loss is ubiquitous across a diverse set of computational environments.

\subsection{The potential evaluation problem}
\label{sec:dissip_fast_summation}

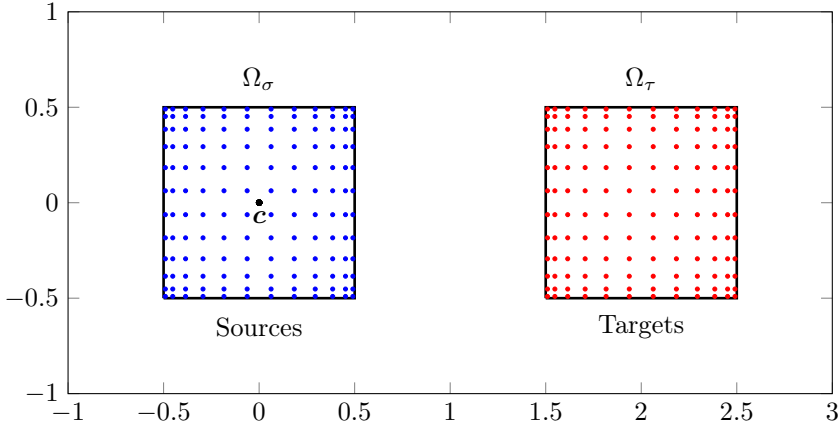
\begin{figure}[!t]
  \centering
  \begin{tikzpicture}
    \begin{axis}[x=1in, y=1in, xmin=-1, xmax=3, ymin=-1, ymax=1 ]
      \addplot[color=black, mark=none, line width=1pt] 
      coordinates {
        (-.5,-.5) (.5,-.5) (.5,.5) (-.5,.5) (-.5,-.5)};
        \addplot[color=black, mark=none, line width=1pt] 
        coordinates {
          (1.5,-.5) (2.5,-.5) (2.5,.5) (1.5,.5) (1.5,-.5)};
      \addplot[only marks, mark size=0.75pt, color=blue]
      table {code/data/sources.dat};
      \addplot[only marks, mark size=0.75pt, color=red]
      table {code/data/targets.dat};
      \node[] at (axis cs: 0,-.65) {Sources};
      \node[] at (axis cs: 0,.65) {\(\Omega_\sigma\)};
      \node[] at (axis cs: 2,-.65) {Targets};
      \node[] at (axis cs: 2,.65) {\(\Omega_\tau\)};
      \addplot[only marks, mark size=1pt, color=black] (0,0);
      \node[] at (axis cs: 0,-.075) {\(\pcc\)};
    \end{axis}
  \end{tikzpicture}
\caption{Geometry of the potential evaluation problem of
  Sections~\ref{sec:intro_diss} and~\ref{sec:dissip_fast_summation}. A continuum
  electric charge density~$q$ in the source box~$\Omega_\sigma$ sampled on a
  tensor-product Gauss-Legendre grid (blue) induces an electric potential~$u$ in
  the target box~$\Omega_{\tau}$, also sampled on an analogous grid (red),
  cf.~(\ref{eq:poteval2}). The point $\pcc$ is the center of the source box. The
  kernel~$\phi$ can also be discretized using the same tensor-product grids,
  from which the singular values of~$\cV$ can be estimated.}
  \label{fig:sources_targets}
\end{figure}

Let us first return to the simple potential evaluation problem that we
introduced in Section \ref{sec:intro_diss}, except in two dimensions.
Specifically, let us consider the integral operator $\cV$ that maps a given
electric charge density~$q$ in a source box $\Omega_\sigma$ to a potential $u$
in a target box~$\Omega_{\tau}$, as illustrated in
Figure~\ref{fig:sources_targets}. With~$\phi$ denoting the fundamental solution
of the Laplace operator in two dimensions, as before in
equation~\eqref{eq:g02d},~$\cV$ takes the form
\begin{equation}
  \label{eq:poteval2}
  u(\pxx) = \cV[q](\pxx) = \int_{\Omega_\sigma}\phi(\pxx - \pyy)
\, q(\pyy) \, d\pyy, \qquad \pxx \in \Omega_{\tau}.
\end{equation}
The key property of $\cV$ that we seek to highlight is that it can, to
arbitrarily high accuracy, be approximated by an operator of finite (low) rank.
To make this precise, let us consider $\cV$ as an operator from
$L^{2}(\Omega_\sigma)$ to $L^{2}(\Omega_{\tau})$. To approximate the singular
values of this operator, we can write down a discretization rule that
approximates the~\mbox{$L^2$-action} of the operator as
an~$\ell^2$-embedding~\cite{2012_bremer_nystrom}. If~$\Omega_\sigma$
and~$\Omega_\tau$ are well-separated, then the kernel~$\phi$ is smooth, and high-order interpolation and quadrature rules can readily be used. If we
discretize the above potential evaluation using weights~$w'_j$ and nodes~$\pxx'_j$,
then
\begin{equation*}
u(\pxx) \approx \sum_j w'_j \, \phi(\pxx-\pxx'_j) \, q_j,
\end{equation*}
where~$q_j = q(\pxx'_j)$. Selecting a set of nodes~$\pxx_i$ and weights~$w_i$
for the target domain, setting~$u_i = u(\pxx_i)$, and scaling the above equation
accordingly, we see that
\begin{equation}
\label{eq:sqrtsc}
\sqrt{w_i} \, u_i \approx \sum_j \sqrt{w_i w'_j} \, \phi(x_i-x'_j) \,  \sqrt{w'_j} q_j.
\end{equation}
The scaling in (\ref{eq:sqrtsc}) ensures that the~$\ell^2$ norm of the vector with
entries~$\sqrt{w'_j} q_j$ converges to the~$L^2$ norm of the continuous
function~$q$ over the domain~$\Omega_\sigma$, i.e.
\begin{equation}
  \left \Vert \sqrt{w'_j} q_j \right \Vert^2_{\ell^2} = \sum_j w'_j \, |q_j|^2 \approx 
  \int_{\Omega_\sigma} |q(\pyy)|^2 \, d\pyy = \Vert q \Vert^2_{L^{2}(\Omega)}.
\end{equation}
Furthermore, it can be shown that 
the spectrum of the matrix~$\mtx{V}$ with entries
\begin{equation}
  \mtx{V}(i,j) =
  \sqrt{w_i w'_j} \, \phi(\pxx_i-\pxx'_j)
\end{equation} converges to the spectrum of the continuum
operator~$\cV$ in (\ref{eq:poteval2}).

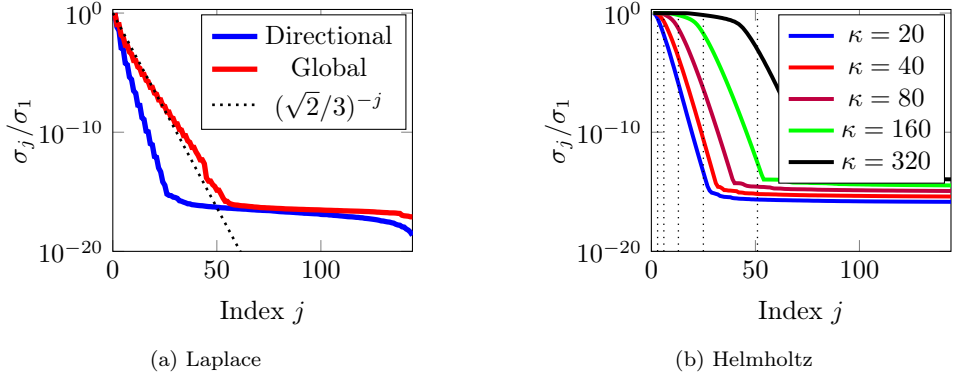
\begin{figure}[t!]
  \centering
  \begin{subfigure}[b]{.45\linewidth}
    \centering
    \begin{tikzpicture}
      \begin{semilogyaxis}[ylabel=\(\sigma_j/\sigma_1 \),
        xlabel=Index \(j\),
        xmin=0, xmax=144, ymin=1.0e-20, ymax=2, width=.95\linewidth]
      \addplot[mark=none, color=blue, line width=2pt] 
        file {code/data/laplace_pot_eval.dat};
        \addlegendentry{Directional}
      \addplot[mark=none, color=red, line width=2pt] 
        file {code/data/laplace_pot_eval_global.dat};
        \addlegendentry{Global}
        \addplot[color=black,dotted, line width=1pt,domain=0:144] 
        plot (\x,0.4714^\x);
        \addlegendentry{$(\sqrt{2}/3)^{-j}$}
      \end{semilogyaxis}
    \end{tikzpicture}
      \caption{Laplace}
      \label{fig:potspec0}
  \end{subfigure}
  \hfill
  \begin{subfigure}[b]{.45\linewidth}
    \centering
      \begin{tikzpicture}
      \begin{semilogyaxis}[ylabel=\(\sigma_j/\sigma_1 \),
        xlabel=Index \(j\),
        xmin=0, xmax=144, ymin=1.0e-20, ymax=2, width=.95\linewidth]
        \addplot[mark=none, color=blue, line width=1.5pt] 
        file {code/data/helmholtz20_pot_eval.dat};
        \addlegendentry{\(\kappa = 20\)}
        \addplot[mark=none, color=red, line width=1.5pt] 
        file {code/data/helmholtz40_pot_eval.dat};
        \addlegendentry{\(\kappa = 40\)}
        \addplot[mark=none, color=purple, line width=1.5pt] 
        file {code/data/helmholtz80_pot_eval.dat};
        \addlegendentry{\(\kappa = 80\)}
        \addplot[mark=none, color=green, line width=1.5pt] 
        file {code/data/helmholtz160_pot_eval.dat};
        \addlegendentry{\(\kappa = 160\)}
        \addplot[mark=none, color=black, line width=1.5pt] 
        file {code/data/helmholtz320_pot_eval.dat};
        \addlegendentry{\(\kappa = 320\)}
        \addplot[color=black, dotted, line width=0.5pt] 
          coordinates {(51,1.0e-20) (51,1.0)};
        \addplot[color=black, dotted, line width=0.5pt] 
          coordinates {(25,1.0e-20) (25,1.0)};
        \addplot[color=black, dotted, line width=0.5pt] 
          coordinates {(13,1.0e-20) (13,1.0)};
        \addplot[color=black, dotted, line width=0.5pt] 
          coordinates {(6,1.0e-20) (6,1.0)};
        \addplot[color=black, dotted, line width=0.5pt] 
          coordinates {(3,1.0e-20) (3,1.0)};
      \end{semilogyaxis}
    \end{tikzpicture}
      \caption{Helmholtz}
      \label{fig:potspeck}
  \end{subfigure}
    \caption{The decay of singular values of the potential evaluation maps
    described in Section~\ref{sec:dissip_fast_summation}. On the left, shown is
    the spectrum for the operators~$\cV$ in~\eqref{eq:poteval2}, and on
    the right is shown the spectrum for various values of~$\kappa$ for the
    operator~$\cV_k$ in~\eqref{eq:poteval3}. Vertical dotted lines are
    placed at predictions for the number of large singular values given by the  
    ``two points per wavelength'' heuristic, as shown in equation~\eqref{eq:nyquist}.}
    \label{fig:potspec}
\end{figure}

In order to estimate the spectrum of the operator~$\cV$, we merely build the
matrix~$\mtx{V}$ and compute its singular value decomposition.
Setting~$\Omega_\sigma$ to be the unit box centered at the origin
and~$\Omega_\tau$ to be the unit box centered at the point~$(2,0)$, see
Figure~\ref{fig:sources_targets}, and choosing a $k \times k$ tensor-product
Gauss-Legendre discretization and quadrature for each domain, we can estimate
the spectrum of~$\cV$. With~$k = 12$ the operator is discretized with 144 points
in each of the source and target domain, and Figure~\ref{fig:potspec0}
``Directional'' shows a
plot of the estimated scaled~spectrum of the operator~$\cV$, i.e. $\sigma_j/\sigma_1$,
where~$\sigma_j$ is the $j$th singular value of the matrix~$\mtx{V}$. Indeed, we see
that to precision~$\varepsilon = 10^{-10}$, its rank is only 17!


The behavior of $\cV$ is mathematically well-understood. A simple proof of the
fact that the singular values of $\cV$ decay exponentially fast follows from a
\textit{multipole expansion} of $\phi$, which in the present context takes the
form
\begin{equation}
\label{eq:mpole_laplace}
\log|\pxx - \pyy| =  \sum_{j=1}^p b_j(\pyy) \, c_j(\pxx) + e_p(\pxx,\pyy).
\end{equation}
In (\ref{eq:mpole_laplace}),  $b_{j}$ and $c_{j}$ represent two sets of
functions\footnote{\label{foot:mpole}The exact form of the functions $b_{j}$ and
$c_{j}$ is not important for our discussion here, but let us nonetheless provide
them for completeness. Letting $\pcc$ denote the center of the source box
$\Omega_{\sigma}$, and introducing polar coordinates
$$
\pxx - \pcc = r(\cos\theta,\,\sin\theta), \qquad\mbox{and}
\qquad \pyy - \pcc =
r'(\cos\theta',\,\sin\theta'),
$$
these functions take the form
$b_{1}(\pxx) = -\log(r)$, $c_{1}(\pyy) = 1$, and then for $j=1,2,3,\dots$
$$
\begin{array}{lrllr}
b_{2j}(\pxx) =& (1/j)\,r^{-j}\cos(j\theta),
&\qquad&
c_{2j}(\pyy) = (r')^{j}\cos(j\theta'),\\
b_{2j+1}(\pxx) =& (1/j)\,r^{-j}\sin(j\theta),
&\qquad&
c_{2j+1}(\pyy) = (r')^{j}\sin(j\theta').
\end{array}
$$
See \cite[Sec.~6.2]{2019_martinsson_fast_direct_solvers} for a detailed
derivation. } that span the range and the co-range of the operator~$\cV$
in~\eqref{eq:poteval2}. Furthermore, one can prove that the remainder
term~$e_{p}(\pxx,\pyy)$ decays exponentially fast whenever the box containing all
source points~$\pyy$ and the box containing all target points~$\pxx$ are
\emph{well-separated} (cf.~\cite{rokhlin1987} and Remark \ref{remark:wellsep}). 
By choosing the number of terms in the sum (\ref{eq:mpole_laplace}) to satisfy~$p \sim \log_2 \epsilon$, we can meet any requested tolerance $\epsilon$. 
For instance, if both boxes are squares with sidelength~$a$,
and the distance between the boxes is also $a$, then
\begin{equation}
\label{eq:mpole_bound}
\sup \{|e_{p}(\pxx,\pyy)|\,\colon\,\pxx \in \Omega_{\tau},\ \pyy \in \Omega_\sigma\}
\leq C\left(\frac{\sqrt{2}}{3}\right)^{p}
\end{equation}
This bound is drawn with a dashed line in
Figure~\ref{fig:potspec0}.  Observe that the bound is far from sharp. This bound
represents the error in the multipole expansion with respect to targets located
in \emph{all} directions, not just a single direction as depicted in
Figure~\ref{fig:sources_targets}.  Adding targets in additional directions would
increase the numerical rank slightly, but not much -- in fact the~$\epsilon$
rank grows to 33.  Figure~\ref{fig:potspec0} ``Global'' shows a plot of the
estimated scaled~spectrum of the operator~$\cV$ when the source box is
surrounded by well-separated target boxes, similar to the configuration in
Figure~\ref{fig:strongweak1}b.
Here, the agreement with the rate~$(\sqrt{2}/3)^{-j}$  is extremely close.

Let us next consider what happens if we replace the fundamental solution $\phi$
associated with the Laplace equation by the fundamental solution $\phi_{\kappa}$
associated with the Helmholtz equation in two dimensions,
\begin{equation}
\label{eq:def_phikappa}
\phi_{\kappa}(\pxx-\pyy) = \frac{i}{4}H_{0}^{(1)}(\kappa|\pxx-\pyy|),
\end{equation}
where $H_{0}^{(1)}$ is the zeroth order Hankel function of the first kind. In
other words, we consider the operator
\begin{equation}
\label{eq:poteval3}
\cV_{\kappa}[\sigma](\pxx) = \int_{\Omega_\sigma}
\phi_{\kappa}(\pxx-\pyy) \, q(\pyy) \, d\pyy, \qquad \pxx \in \Omega_{\tau},
\end{equation}
again as an operator from $L^{2}(\Omega_\sigma)$ to $L^{2}(\Omega_{\tau})$.
Figure~\ref{fig:potspeck} shows the singular spectrum of $\cV_{\kappa}$
for~$\kappa = $ 20, 40, 80, 160, 320. In this case, to a relative accuracy
of~$10^{-10}$, the ranks are 19, 24, 31, 45, and 70, respectively. Due to the
oscillatory nature of the kernel, these ranks were obtained using a
finer~$64 \times 64$ Gauss-Legendre discretization. We again find that the
singular values decay exponentially, but with the crucial
distinction that the exponential decay starts \emph{only after} a relatively
flat plateau where they remain of roughly the same size. This behavior is well
understood~\cite{osipov2013prolate, 2006_rokhlin_wideband,
  crutchfield2006remarks, michielssen1996multilevel}, and we know that the
number of constant singular values is approximately bounded by
\begin{equation}
\label{eq:nyquist}
R_{\kappa} \approx \frac{\kappa D}{2\pi},
\end{equation}
where $D$ is the diameter of the box~\cite{bucci1987spatial,bucci1989degrees}.
The estimate~\eqref{eq:nyquist} is one manifestation of the ubiquitous principle
that at least ``two points per wavelength'' is required to fully resolve an
oscillatory field. In Figure~\ref{fig:potspeck}, we indicate the predicted
location of the ``knee''~\eqref{eq:nyquist} by vertical dotted lines; the match is
remarkable. These results all follow from an analysis of the known properties of
the fundamental solution~$\phi_{\kappa}$. In particular, one can derive a
separation of variables that is analogous of the multipole expansion described
in Footnote~\ref{foot:mpole} for the Laplace equation. Assuming
that~$|\pxx| > |\pxx'|$, one such expansion takes the form of the famous Graf addition
formula~\cite{olver2010nist}:
\begin{equation}
  H_0^{(1)}(\kappa|\pxx-\pxx'|) = \sum_{\ell=-\infty}^\infty J_\ell(\kappa r')
  \, e^{-i\ell \theta'} \, H^{(1)}_\ell(\kappa r) \, e^{i\ell \theta},
\end{equation}
where the polar coordinates of~$\pxx$ and~$\pxx'$ are given as~$(r,\theta)$
and~$(r',\theta')$, respectively, and~$J_n$ and~$H^{(1)}_n$ are Bessel functions of
the first and third kind, respectively. It is also possible to tighten the
bound by replacing the ``diameter'' with the minimum viewing angle of the source
and target boxes, respectively~\cite{2018_engquist_zhao_helmholtz}.

\begin{remark}
  Many of the salient features regarding the dissipation of information with
regard to integral equations follows directly from the previous discussion on
potential evaluations -- most of the time, the kernel of an integral equation
is the Green's function, or various derivatives of the Green's function, for
the associated PDE. Upon discretization of the integral equation, submatrices
corresponding to well-separated interactions take the form of~$N$-body
operators, i.e. merely non-singular potential evaluation problems. The more
accurate the discretization is, the more accurately the spectrum of these
operators mimics the true underlying spectrum of the continous operator, but
nonetheless, the ranks often decay exponentially fast for low-order
discretizations and to below machine precision in many cases.
\end{remark}

\subsection{Schur complement matrices in sparse LU factorization}
\label{sec:dissip_frontal}

Let us next consider large sparse matrices that arise upon discretization of
BVPs using, for example, finite element or finite difference methods. Given such a
matrix $\mtx{A}$, we seek an LU factorization $\mtx{A}(I,J) = \mtx{LU}$.  We
recall from Section~\ref{sec:FDSPDE} that the core technical obstacle is that in
the process of the factorization, sparsity is gradually lost as larger and
larger dense matrices arise. The loss of sparsity can be
slowed down by picking good orderings $I$ and $J$ of the rows and columns, but
remains a serious hurdle even with optimal re-orderings.  In this section, we seek
to substantiate our claim that these matrices have exploitable rank structure,
and specifically, that they have very similar rank structure to the discretized 
integral operators in Section \ref{sec:dissip_fast_summation}.

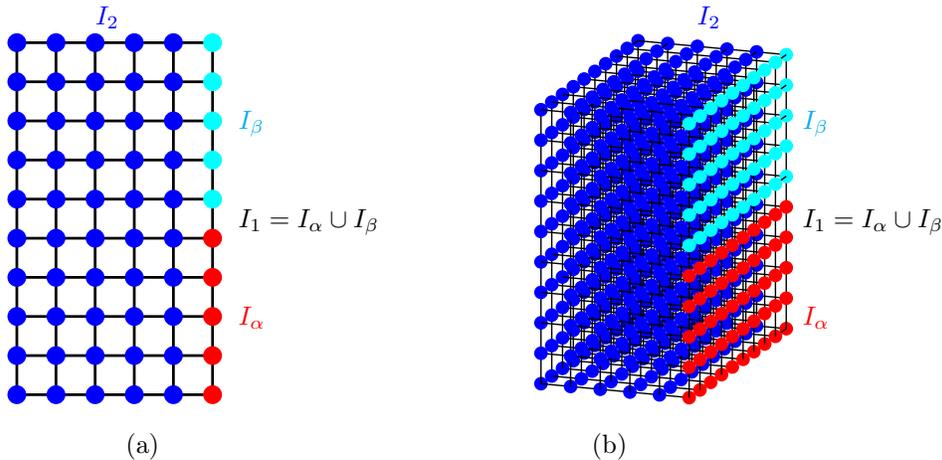
\begin{figure}
\centering
\setlength{\unitlength}{1mm}
\begin{picture}(120,60)
\put(00,06){\includegraphics[height=50mm]{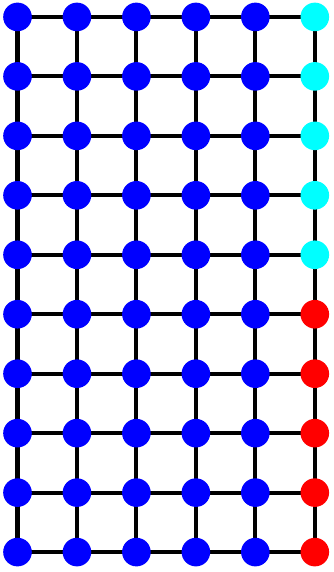}}
\put(12,57){\color{blue}$I_{2}$}
\put(31,17){\color{red}$I_{\alpha}$}
\put(31,30){$I_{1} = I_{\alpha} \cup I_{\beta}$}
\put(31,43){\color{cyan}$I_{\beta}$}
\put(70,06){\includegraphics[height=50mm]{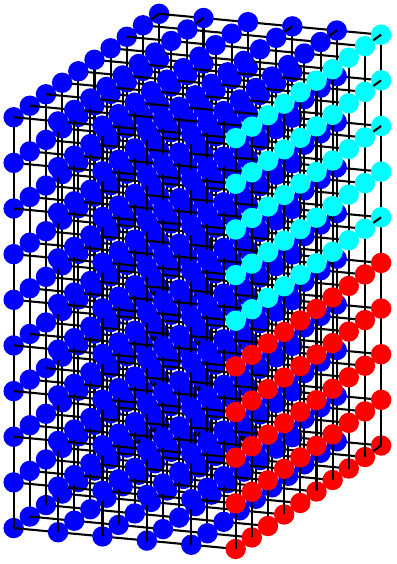}}
\put(92,57){\color{blue}$I_{2}$}
\put(106,17){\color{red}$I_{\alpha}$}
\put(106,30){$I_{1} = I_{\alpha} \cup I_{\beta}$}
\put(106,43){\color{cyan}$I_{\beta}$}
\put(16,00){(a)}
\put(78,00){(b)}
\end{picture}

\caption{The geometry discussed in Section \ref{sec:dissip_frontal},
shown for $n=10$.
(a) A grid in two dimension. The edge points in $I_{1}$ are partitioned into two
sets $I_{\alpha}$ (red) and $I_{\beta}$ (cyan).
(b) The analogous partitioning of the points in three dimensions.}
\label{fig:Salphabeta}
\end{figure}

To keep things simple, we return to the model problem introduced in
Section~\ref{sec:NDLU} where we used the standard five-point difference stencil
to discretize the 2D boundary value problem (\ref{eq:buck1}) to obtain a linear
system $\mtx{A}\vct{u}=\vct{b}$.  Recall that we partitioned the domain into
three sets $I_{1}$, $I_{2}$, and $I_{3}$ that represent the middle, the left, and the right nodes in the mesh, see~Figure~\ref{fig:nd1}a. Subsequently, we saw
that if we use a divide-and-conquer strategy for computing an LU factorization
of~$\mtx{A}$, the dominant cost was the factorization of the Schur complement
$ \mtx{S} = \mtx{A}_{11} - \mtx{A}_{12}\mtx{A}_{22}^{-1}\mtx{A}_{21} -
\mtx{A}_{13}\mtx{A}_{33}^{-1}\mtx{A}_{31}$ (which is a dense matrix). We asserted 
in Section~\ref{sec:fastLU} that this dense matrix is rank structured, and
furthermore, that this fact can be exploited to obtain an LU factorization in
only~$\cO(N)$ operations.

\begin{figure}[t]
    \centering
    \textit{\small{Laplace (2D)}}
    \hspace{25mm}
    \textit{\small{Helmholtz (2D)}}
    \hspace{25mm}
    \textit{\small{Laplace (3D)}}
    
    \includegraphics[width=\textwidth]{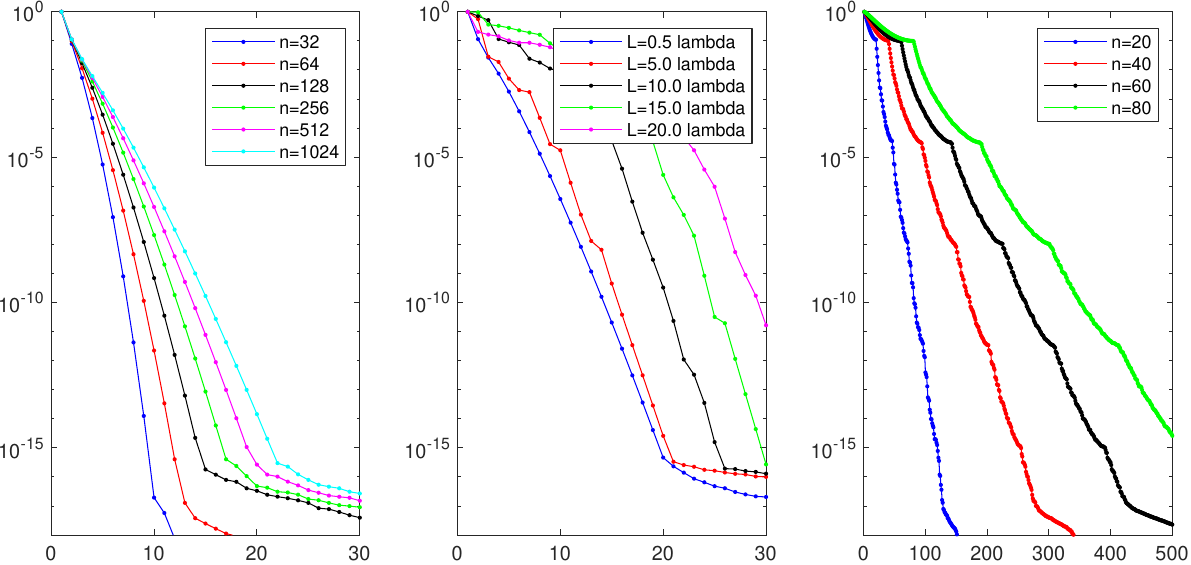}

    (a)\hspace{40mm}(b)\hspace{40mm}(c)
    
    \caption{The rate of decay in the singular values of an off-diagonal block
    in a Schur complement arising in the LU factorization of a finite difference
    matrix, as described in Section~\ref{sec:dissip_frontal}. The graphs show the
    relative singular values~$\sigma_{j}/\sigma_1$ of the blocks~$\mtx{S}_{\alpha,\beta}$
    defined by~\eqref{eq:Salphabeta}.}
    \label{fig:schur_ranks}
\end{figure}

In this section, we will look more closely at the singular values of the
off-diagonal blocks of the matrix $\mtx{S}$ to get a sense of what the numerical
ranks are in practice. To be precise, let us split the index vector $I_{1}$ into
two parts, $I_{1} = I_{\alpha} \cup I_{\beta}$, as shown in
Figure~\ref{fig:Salphabeta}(a), and then define
\begin{equation}
\label{eq:Salphabeta}
\mtx{S}_{\alpha,\beta} = 
\mtx{A}(I_{\alpha},I_{2})\mtx{A}_{22}^{-1}\mtx{A}(I_{2},I_{\beta}).
\end{equation}
The matrix~$\mtx{S}_{\alpha,\beta}$ is the largest off-diagonal
block in a HODLR representation of $\mtx{A}_{12}\mtx{A}_{22}^{-1}\mtx{A}_{21}$.
Figure~\ref{fig:schur_ranks} shows the singular values of $\mtx{S}_{\alpha,\beta}$
for three different test cases: 
\renewcommand{\labelenumi}{(\alph{enumi})}
\begin{enumerate}
    \item Laplace's equation $-\Delta u = f$ on a uniform 2D mesh of size
    $n\times n$, shown for increasing values of $n$.
    \item Helmholtz equation $-\Delta u - \kappa^{2} u = f$ on a uniform 2D mesh of size
    $512\times 512$. The lines correspond to different value of $\kappa$, 
    chosen so that the number of oscillations across the $L\times L$ square 
    are $0.5,\,5,\,10,\,15,\,20$, respectively.
    \item Laplace's equation $-\Delta u = f$ on a uniform 3D mesh of size
    $n\times n \times n$, shown for increasing values of $n$.
\end{enumerate}

Figure~\ref{fig:schur_ranks}(a) shows that for a 2D Laplace problem, rapid
exponential decay occurs for all problem sizes.  The numerical ranks increase as
the mesh size is refined, but only very slowly.  Figure~\ref{fig:schur_ranks}(b)
shows that in the case of problems with highly oscillatory solutions, almost no
decay is observed until the waves have been resolved to about two points per
wavelength. After that, the same rapid decay that was observed in the Laplace
case kicks in.  Lastly, Figure~\ref{fig:schur_ranks}(c) shows that for a 3D Laplace
problem, we still see exponential decay, but at a slower rate resulting
in significantly higher numerical ranks for any given precision. (This can be
remedied by switching from HODLR to a format that relies on ``strong
admissibility'', see Section~\ref{sec:BIE_strong}.)

It is of particular interest that in all cases the decay continues through
fifteen orders of magnitude all the way down to machine precision. This
indicates that the property we see is inherent to the matrix~$\mtx{A}$ since the
precision in discretizing the PDE is a few digits (at best) for a basic finite
difference discretization. The underlying \emph{continuous} Schur complement
operator shares this behavior as well, but we need not discretize it to
high-order in order to realize the behavior numerically. We will also remark
that the overall rank patterns are remarkably similar to what we observed for
potential evaluation problems in boundary integral equations in
Section~\ref{sec:dissip_fast_summation}.

\subsection{Domain decomposition via Poincar\'e-Steklov operators}
\label{sec:dissip_FEMBEM}

Boundary-to-boundary operators such as the Dirichlet-to-Neumann operator are
powerful tools for ``gluing together'' different subdomains in a computational
simulation.  The domains may involve regions with different physical
constitutive equations (e.g.~fluid-solid coupled problems), or, they may involve
the same equations but discretized using different numerical methods (e.g.,
coupling local direct solvers with an iteration for the global system).  The
Poincar\'e-Steklov operators provide a convenient way to represent continuity
conditions at the interfaces that ensure uniqueness of the solution.

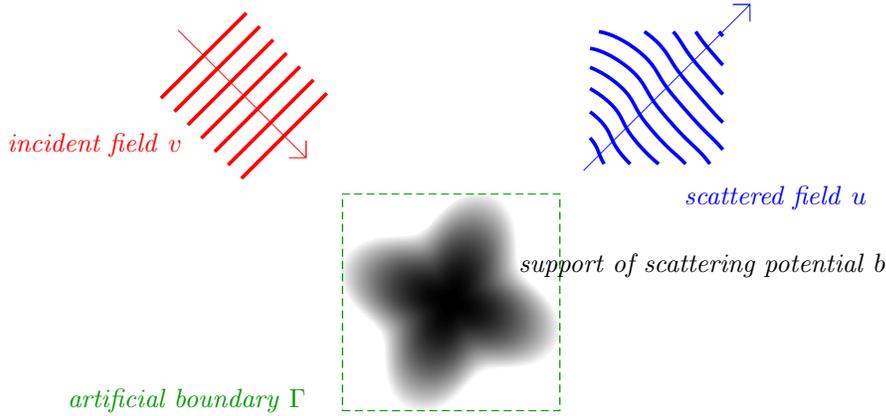
\begin{figure}
\centering
\setlength{\unitlength}{1mm}
\begin{picture}(116,56)
\put(20,00){\includegraphics[height=55mm]{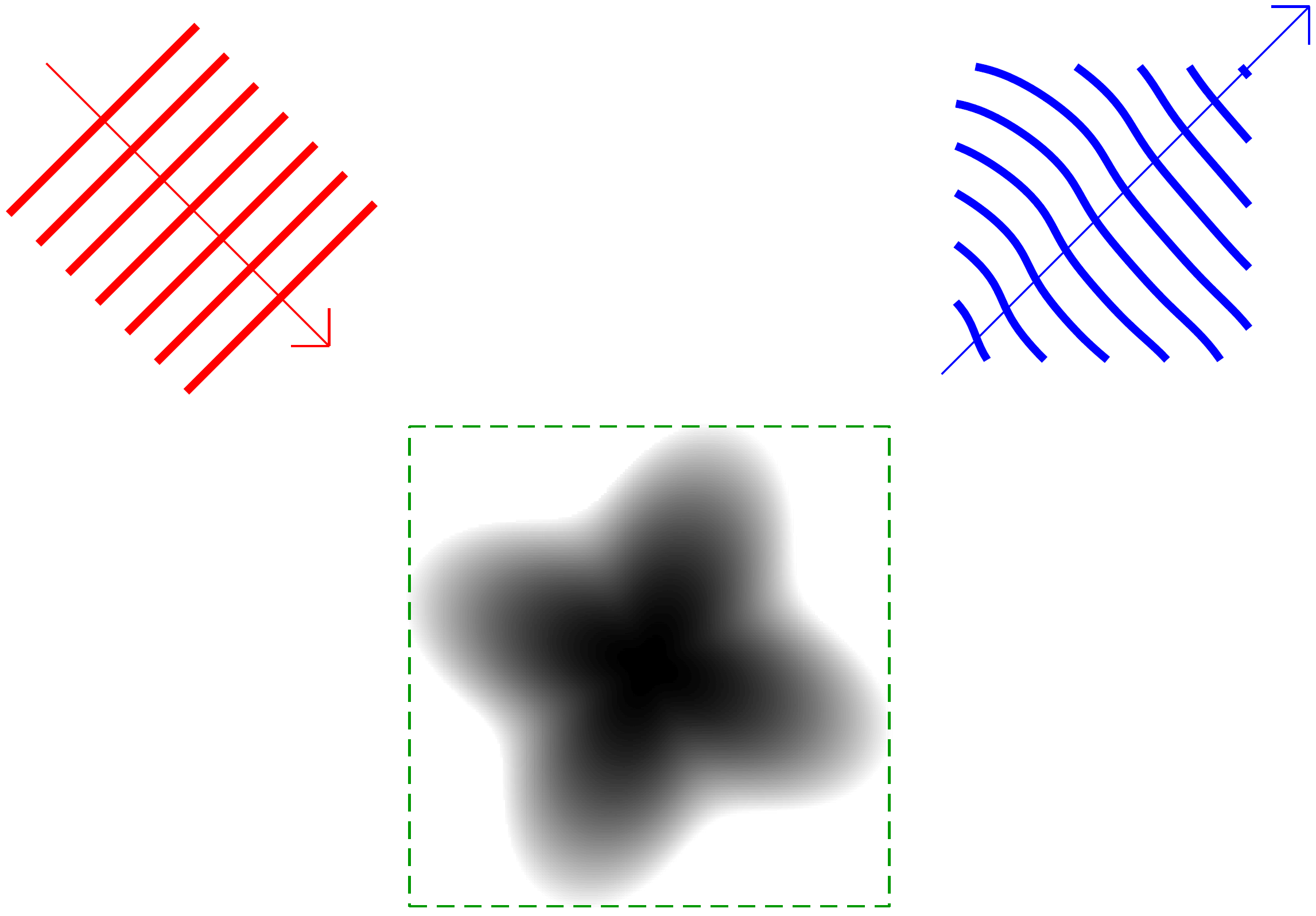}}
\put(00,35){\color{red}\textit{incident field} $v$}
\put(90,28){\color{blue}\textit{scattered field} $u$}
\put(68,19){\textit{support of scattering potential} $b$}
\put(08,01){\color{MyDarkGreen}\textit{artificial boundary} $\Gamma$}
\end{picture}
\caption{The scattering problem described in Section \ref{sec:dissip_FEMBEM}.
  Waves propagate in a medium with constant wave number $\kappa$ everywhere
  except in the area shaded gray, where the wave number is
  $\kappa\,\sqrt{1-b(\pxx)}$ for a given smooth compactly supported scattering
  potential $b = b(\pxx)$.  An incident field $v$ hits the scattering potential
  and induces the scattered field $u$.  The dashed green line marks an artificial
  boundary $\Gamma$ which encloses the support of the scattering potential.}
\label{fig:ItI_problem}
\end{figure}

In our final illustration of dissipation of information, we investigate the rank
property of a boundary-to-boundary operator that is used to glue together two
subdomains that have been discretized using different techniques.  Specifically,
we investigate an acoustic wave-propagation problem where an incoming field hits
a region of variable wave-speed. Our task is to compute the scattered field,
cf.~Figure \ref{fig:ItI_problem}.  One approach is to impose an artificial box
$\Omega$ with boundary $\Gamma$ on the domain, in such a way that it encloses
the region with variable wave-speed. Then use a BIE formulation to model the
constant coefficient problem exterior to $\Gamma$, and a direct discretization
of the PDE that models the variable coefficient problem inside $\Gamma$. To glue the
two domains together, we enforce continuity of both potentials and their normal
derivatives using the impedance-to-impedance (ItI) map
\cite{1994_monk_ItI, 2013_martinsson_ItI}.  (The ItI operator maps one
impedance boundary condition $u + iu_{n}$ on $\Gamma$ to its conjugate
$u - iu_{n}$.) 

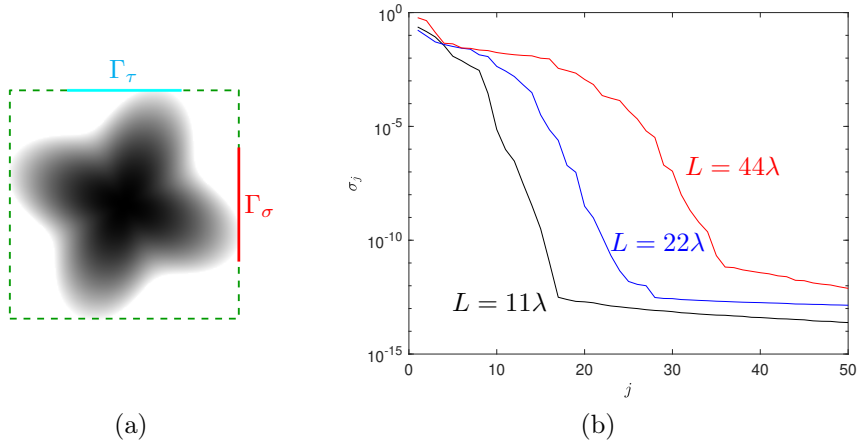
\begin{figure}
\begin{center}
\begin{tabular}{ccc}
\setlength{\unitlength}{1mm}
\begin{picture}(35,55)
\put(00,10){\includegraphics[width=32mm]{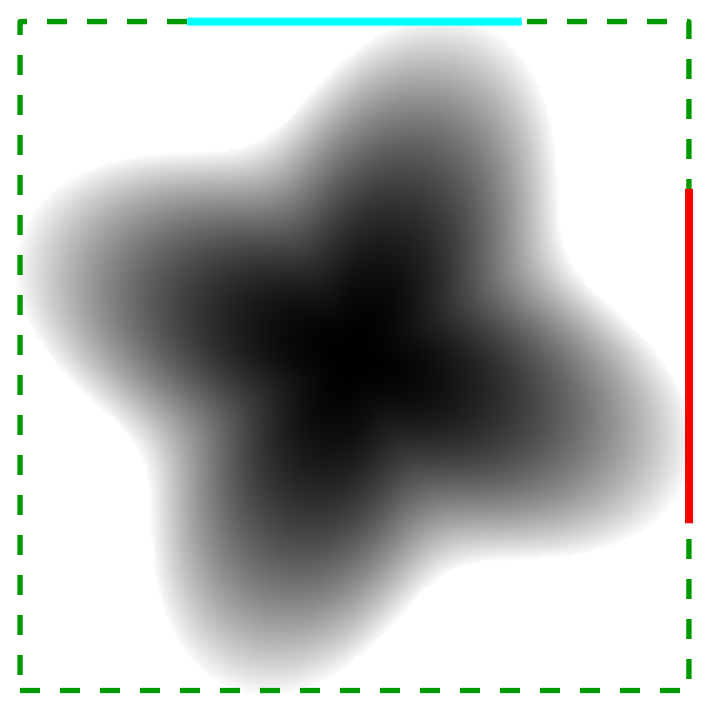}}
\put(14,43){\color{cyan}$\Gamma_{\tau}$}
\put(32,25){\color{red} $\Gamma_{\sigma}$}
\end{picture}
&\mbox{}\hspace{1mm}\mbox{}&
\setlength{\unitlength}{1mm}
\begin{picture}(68,55)
\put(00,00){\includegraphics[height=53mm]{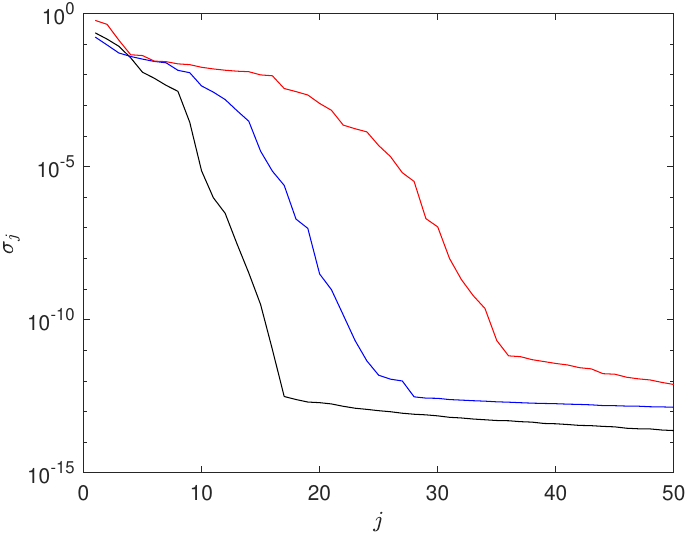}}
\put(14,12){\color{black}$L=11\lambda$}
\put(35,20){\color{blue}$L=22\lambda$}
\put(45,30){\color{red}$L=44\lambda$}
\end{picture}
\\
(a) && (b)
\end{tabular}
\end{center}
\caption{Numerical ranks of the ItI operator described in Section \ref{sec:dissip_FEMBEM}. (a) The two sub-regions $\Gamma_{\sigma}$ and $\Gamma_{\tau}$ of the problem described in Figure \ref{fig:ItI_problem}. (b)
The singular values $\sigma_{j}$ of the matrix $\mtx{A}(I_{\tau},I_{\sigma})$, where $\mtx{A}$  is an approximation to the ItI operator for the interior Helmholtz problem $-\Delta u - (2\pi/\lambda)^{2}(1-b)u=0$, where $b$ is the scattering potential shown in (a) (see Section \ref{sec:dissip_FEMBEM}). For a box size of $L\times L$, we picked three different wave-lengths $\lambda$, chosen so that $L = 11\lambda,\,22\lambda,44\lambda$, respectively.}
\label{fig:ItI_spectrum}
\end{figure}

The dense matrix $\mtx{A}$ that we consider in this case is an approximation to
the ItI operator that is computed through a high order discretization of the
interior BVP. Specifically, we use an $18^{\rm th}$-order accurate multidomain
spectral collocation method~\cite{2013_martinsson_ItI}. To study the rank
behavior of $\mtx{A}$, we isolate two subdomains~$\Gamma_{\sigma}$
and~$\Gamma_{\tau}$ within~$\Gamma$, and plot the singular values of the
corresponding restriction of $\mtx{A}$ in Figure \ref{fig:ItI_spectrum}. We
again observe basically the same rank patterns as we saw in
Section~\ref{sec:dissip_frontal}. (Similar rank behaviors would appear in other
applications when using similar spectral element methods, such as for surface
PDEs in~\cite{fortunato2024high} .)

\subsection{Summary}
In this section we have illustrated, through a diverse set of examples, that a
variety of dense matrices arising in the context of fast solvers for elliptic
PDEs all share the property of being \emph{rank-structured}. Specifically, the
singular values of off-diagonal blocks of the matrix often decay exponentially
fast.

The rank-structured property is persistent across different discretization
methods, whether through finite element or finite difference discretizations of
the original PDE, or through collocation or boundary element methods for the
associated integral equations. Remarkably, the fast rank decay often continues 
all the way down to machine precision, even in cases where the discretization 
of the continuum problem is accurate only to a couple of digits.


%% file: 11-compression/compression02.tex
At this point, we have described how dense but rank-structured matrices arise in
a number of different applications, and how this structure can be exploited to
accelerate many standard linear algebraic tasks. Next, we will describe a
number of techniques for obtaining the data sparse representation in
the first place. In other words, given a rank-structured matrix $\mtx{A}$ of
size $N \times N$, how does one find the collection of small dense matrices and low-rank factors that define it? One sometimes has direct access to the
entries of~$\mtx{A}$ (e.g.~in the integral equation setting), but often the access is implicit in that we can only observe the action of $\mtx{A}$ on vectors. Each of these situations requires a different approach to matrix compression.

Throughout this section, we assume that the tree that defines the off-diagonal
blocks is known.  Furthermore, we let $k$ denote the rank of the off-diagonal
blocks. In somewhat of an abuse of notation, we will use the terms ``rank'' and
``$\epsilon$-rank'' interchangeably.
What follows is a brief
description of some of the techniques that can be used to obtain low-rank
representations of submatrices of these hierarchically compressible operators.

\subsection{Brute force compression}
\label{sec:comp_bruteforce}

The most straightforward method is to directly form the entire matrix $\mtx{A}$
and store it as an array of $N\times N$ numbers. Then, one can loop over all
compressible blocks and directly compress them using standard factorization
methods such as forming a full SVD and truncating it, executing column pivoted
QR with a stopping criterion, or using a randomized method
\cite{2007_martinsson_PNAS,2011_martinsson_randomsurvey,2020_martinsson_acta,buttari2025truncated}.
Such an approach is attractive in that close-to-optimal compression rates are
easy to obtain, but have the obvious drawback that the overall cost incurred is
no less than~$\cO(N^2)$, and is typically $\cO(kN^2)$, where~$k$ is the
approximate rank of the compressible blocks.

A standard use case for brute force compression is when rank-structured matrices
are used to accelerate sparse direct solvers. In this context, the Schur
complements that arise are often formed via complex operations that involve a
mix of rank-structured, dense, and sparse blocks. While it might in theory be
possible to form the Schur complement by working strictly with data sparse
representations, the practical complexity is such that practitioners often fall
back on the much simpler brute force compression approach. (Sometimes with the
compression step accelerated using randomized
methods~\cite{2023_malik_li_HSS_randomized}.)

Additional compression can sometimes be achieved by using
mixed-precision formats or after-the-fact by mere
compression of floating point numbers (which often contain unneeded
information with regard to the low-rank accuracy)~\cite{kriemann2025hierarchical,carson2025mixed,higham2022mixed}.

\subsection{Interpolative black box algorithms}
\label{sec:comp_ibb}

In the case where one has access to the matrix entries in the form of a smooth
kernel~$\Phi$, i.e.~$\mtx{A}(i,j) = \Phi(\pxx_i,\pxx_j)$, such as in the
integral equation setting, a direct \emph{continuous} interpolative
approximation can be made using polynomials, as described
in~\cite{fong2009black}:
\begin{equation}
\label{eq:polyc}
\Phi(\pxx,\pyy) = \sum_{i=1}^k \sum_{j=1}^k 
\Phi(\pxx_i,\pyy_j) \, p_i(\pxx) \, q_j(\pyy) 
=   \sum_{i=1}^k p_i(\pxx) \sum_{j=1}^k 
\Phi(\pxx_i,\pyy_j) \,  q_j(\pyy),
\end{equation}
where~$p_i$ and~$q_j$ are the interpolating functions, usually chosen to be some
combination of Chebyshev polynomials in order to facilitate a subsequent fast
algorithm for computing the sum (\ref{eq:polyc}). Explicit formulas can be obtained
for~$p_i$ and~$q_j$, much in the same way that explicit forms exist for Lagrange
or barycentric interpolation, and for non-oscillatory kernels~$\Phi$ the
resulting approximation can be performed to high order. Note that the above form
of the interpolant shows that the kernel~$\Phi$ has been factored as the product
of two rank-$k$ operators. The rank $k$ can be set \emph{a priori}, or
determined adaptively, on-the-fly, in order to obtain an~$\epsilon$-accurate
approximation. This adaptive approach involves increasing the number of
interpolation nodes and stopping when sufficient agreement has been achieved in
evaluating~$\Phi$ off-grid.

The approximation method for $\Phi$ that we described was used in~\cite{fong2009black} as a ``black box fast multipole method'' for computing fast matrix-vector multiplications, but it also applies in the context of fast direct solvers for general kernels. 
While this technique is very general, it often sub-optimal in that it leads to excessively large ranks for any given approximation tolerance~$\epsilon$;
additional numerical re-compression can be used to obtain tighter factorizations.

\subsection{The proxy surface technique}
\label{sec:comp_proxy}

\begin{figure}[t]
  \centering
  \begin{subfigure}[t]{.55\linewidth}
    \centering \includegraphics[width=.95\linewidth]{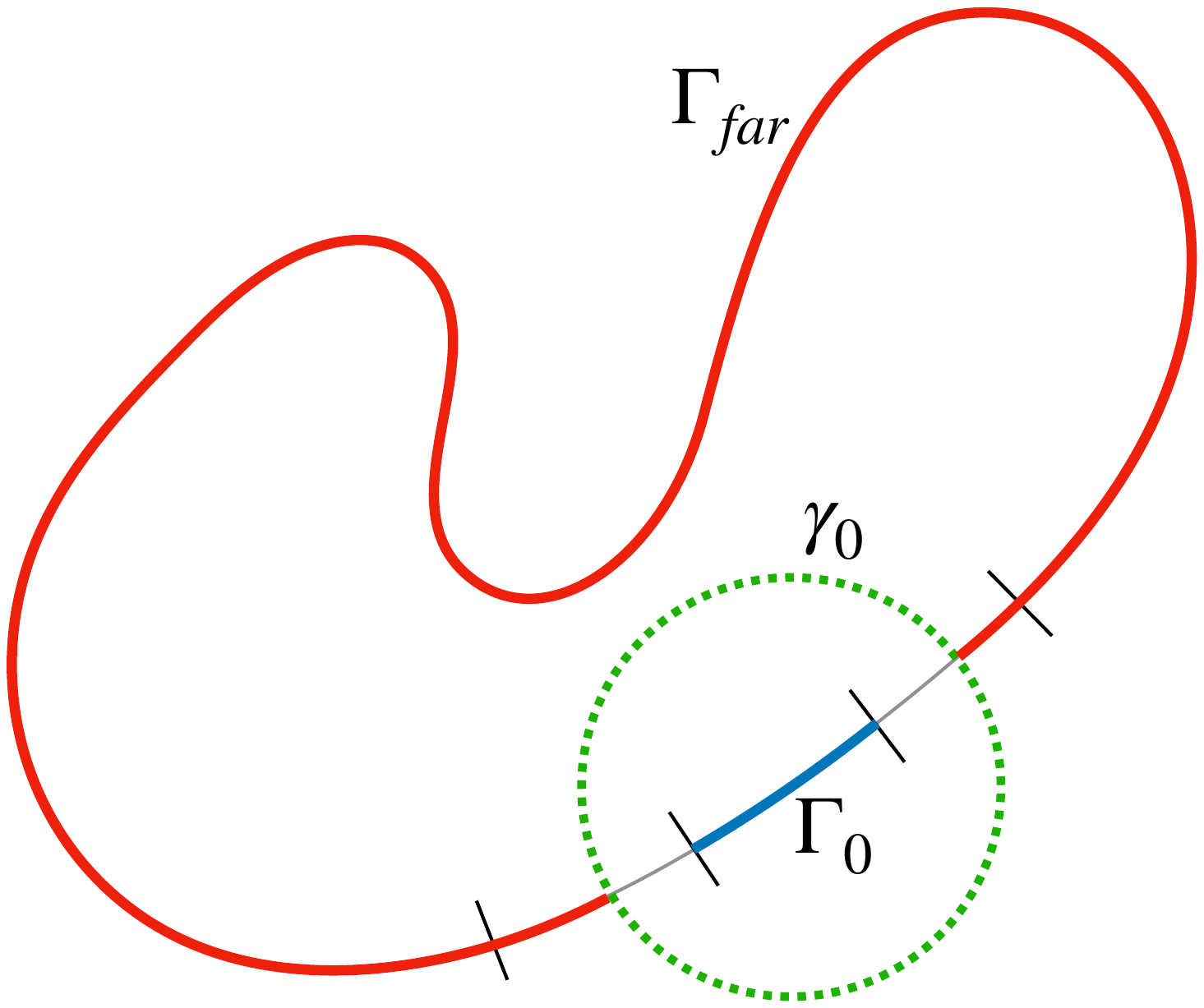}
    \caption{A discretization of a curve in 2D, along with source and target
      panels and an enclosing proxy surface.}
    \label{fig:proxy_curve}
  \end{subfigure}
  \quad
  \begin{subfigure}[t]{.4\linewidth}
    \includegraphics[width=.95\linewidth]{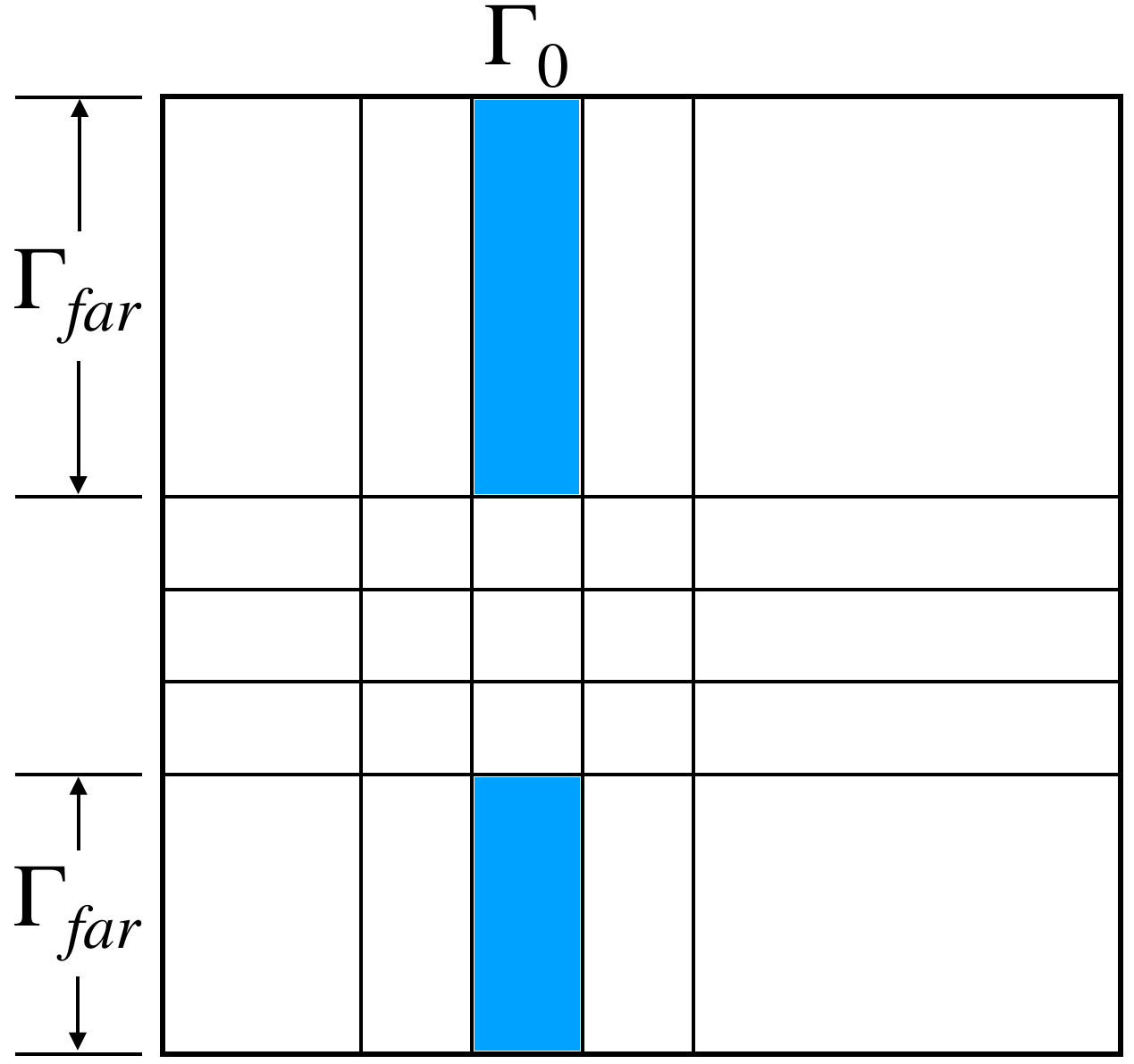}
    \caption{The corresponding matrix blocks to be compressed.}
    \label{fig:proxy_matrix}
  \end{subfigure}
  \caption{Setup of a proxy compression.}
  \label{fig:proxy}
\end{figure}

In the case where the matrix at hand comes from the discretization of an
integral equation, compressed forms of off-diagonal submatrices can be obtained
using what are known as \emph{proxy surfaces} and some simple theory supported
by Green's reproducing identity. Consider the integral equation
in~\eqref{eq:intdir}, previously discussed in Section~\ref{sec:BIE2D}, for the
interior Dirichlet problem for the Laplace equation in two dimensions:
\begin{equation}
  -\frac{1}{2} \sigma(\pxx) + \int_\Gamma d(\pxx,\pyy) \, \sigma(\pyy) \,
  ds(\pyy)
  = f(\pxx), \qquad \text{for } \pxx \in \Gamma,
\end{equation}
where as before the kernel~$d$ is known as the double layer kernel
\begin{equation}
d(\pxx,\pyy) = \pnn(\pyy) \cdot \nabla_y \phi(\pxx-\pyy)
\end{equation}
and the potential~$u$ has been represented as a double
layer,~${u = \mathcal D[\sigma]}$.  After (a Nystr\"om) discretization and
imposing a hierarchical structure on the matrix, our task is to find a
low-rank representation of blocks of the matrix containing entries of the form
\begin{equation}
  \mtx{A}(i,j) = w_j \, d(\pxx_i,\pxx_j),
\end{equation}
where~$w_j$ is a quadrature weight and where the
\emph{targets}~$\pxx_i$ are clearly separated from the
\emph{sources}~$\pxx_j$. Considering as an example the geometry shown in
Figure~\ref{fig:proxy}, let us seek to construct a low-rank
approximation of the submatrix encoding the interaction of
sources~$\pxx_j$ located on panel~$\Gamma_0$ with all of the targets located on
panel~$\Gamma_{\text{far}}$.  Denoting the indices of the points on~$\Gamma_0$
as~$J$ and those on~$\Gamma_{\text{far}}$ as~$I$, we thus seek a low-rank
representation of~$\mtx{A}(I,J)$, which is highlighted in blue in
Figure~\ref{fig:proxy_matrix}.

The main fact that allows for an efficient compression scheme is the existence
of Green's Third Identity~\cite{colton2012}. For~$\pyy \in \Gamma_0$, the
function~$d(\cdot,\pyy)$ is harmonic outside any region enclosing~$\Gamma_0$,
namely the \emph{proxy surface}~$\gamma_0$ depicted in
Figure~\ref{fig:proxy_curve}. This means that
for~$\pxx \in \Gamma_{\text{far}}$,~$d(\cdot,\pyy)$ can be written in terms of
its boundary values on~$\gamma_0$:
\begin{equation}
\label{eq:GreenThird}
  \begin{aligned}
    d(\pxx,\pyy) &= \int_{\gamma_0} \left( \frac{\partial }{\partial n_w} 
                   \phi(\pxx,\pww) \, d(\pww,\pyy) - \phi(\pxx,\pww) \frac{\partial }{\partial n_w}
                   d(\pww,\pyy)
                   \right)
                   ds(\pww) \\
                 &= \int_{\gamma_0} \left( d(\pxx,\pww) \, d(\pww,\pyy) - s(\pxx,\pww)  \, d'(\pww,\pyy)
                   \right)
                   ds(\pww),
  \end{aligned}
\end{equation}
where~$d'(\pww,\pyy)$ denotes the directional derivative of~$d(\cdot,\pyy)$ at
the point~$\pww \in \gamma_0$ in the direction of~$\pnn(\pww)$, the outward unit
normal to~$\gamma_0$ at~$\pww$. Formula (\ref{eq:GreenThird}) effectively decouples the
sources on~$\Gamma_0$ from the targets on~$\Gamma_{\text{far}}$, which
for~$\pxx \in \Gamma_{\text{far}}$, results in a factorization of the potential:
\begin{equation}
\label{eq:mathfact}
  \begin{aligned}
  u(\pxx) &= \int_{\Gamma_0} d(\pxx,\pyy) \, \sigma(\pyy) \, ds(\pyy) \\
  &= \int_{\Gamma_0} \int_{\gamma_0} 
  \left( d(\pxx,\pww) \, d(\pww,\pyy) - s(\pxx,\pww)  \, d'(\pww,\pyy) \right)
  ds(\pww) \, \sigma(\pyy) \, ds(\pyy) \\
  &= \int_{\gamma_0}     
  \left( d(\pxx,\pww)  - s(\pxx,\pww)  \frac{\partial }{\partial n}  \right) \int_{\Gamma_0}  d(\pww,\pyy) \, 
   \sigma(\pyy) \, ds(\pyy) \, ds(\pww) \\
   &= \left( \mathcal G_0 \circ \mathcal D \right) [\sigma](\pxx).
  \end{aligned}
\end{equation}
After a discretization of the integrals in (\ref{eq:mathfact}), we have that
\begin{equation}
\vct{u}(I) = \mtx{G}_0(I,P) \, \mtx{A}(P,J),
\end{equation}
where~$P$ denotes indices of \emph{proxy points} along the curve~$\gamma_0$
and~$\mtx{G}_0$ denotes a discretization of Green's Theorem mapping sources
on~$\gamma_0$ to potentials on~$\Gamma_{\text{far}}$. 

The take-home message of this calculation is that the original
submatrix~$\mtx{A}(I,J)$ can be analytically factored into two pieces which
couple only through a linear operator defined on a fictitious proxy surface.
The number of discretization points on~$\gamma_0$ is completely independent of
the number of discretization points on the original
geometry~$\Gamma_{\text{far}}$. The matrix~$\mtx{G}_0$ may still be large, i.e.
$\cO(N \times n_P)$ with~$n_P \sim \cO(1)$, but the matrix~$\mtx{A}(P,J)$ is
small, i.e. of size~$n_P \times n_J$.  A low-rank factorization of this latter
matrix, matrix~$\mtx{A}(P,J)$, can be obtained using dense linear algebra in a
complexity that has nothing to do with~$N$, the discretization size of the
overall problem. An analogous procedure and factorization could be performed
for the dual compression, i.e. sources located on~$\Gamma_{\text{far}}$ and
targets located on~$\Gamma_0$.

In practice, the full Green's Identity representation can often be replaced by
merely a ``single layer representation'' without significant loss of accuracy in
the compression. Inside of a hierarchical algorithm, the previous proxy
procedure can be repeated by ``merging'' compressed forms of
neighboring~$\mtx{A}(P,J)$ matrices and surrounding the associated pieces of
geometry with a larger proxy surface.  See~\cite{sushnikova2023fmm,
  2017_ho_ying_strong_RS, 2012_greengard_ho_recursive_skeletonization,
  2009_martinsson_ACTA, 2005_martinsson_fastdirect, xing2020interpolative} for
different implementations of such ideas, and how to efficiently integrate proxy
compression into hierarchical compression schemes. Green's identity methods can
even be used to simultaneously compress families of parameter dependent
operators (such as a Helmholtz integral operator for a broad interval of wave
numbers), often at surprisingly small cost in terms of increased ranks
\cite{2023_martinsson_gopal_broadband}.

\subsection{Adaptive cross approximation}
\label{sec:comp_aca}

A very common technique to obtain low-rank factorizations in the hierarchical
matrix literature is known as adaptive cross approximation
(ACA)~\cite{goreinov1997theory,tyrtyshnikov1996mosaic,bebendorf2011adaptive,borm2005hybrid,
  bebendorf2003adaptive}. The idea behind both cross approximation and
interpolative
decompositions~\cite{martinsson2007interpolation,2005_martinsson_skel} is very
similar: find a selection of rows and/or columns that form bases for the row
and/or column space, and then use those rows/columns to express the remaining
ones. 
Theses methods generally differ significantly in the
algorithms used to obtain the factorization, and in the overall structure of the
factorization.

Recently there has been work toward addressing both the algorithmic and accuracy
properties of ACA, both in the context of compressing smooth (non-oscillatory)
kernels~\cite{borm2005hybrid,borm2016approximation, tamayo2011multilevel} as
well as its use in the oscillatory case (although the approach is slightly
different)~\cite{borm2017approximation}. As in the interpolative decomposition,
the main challenge is to find a selection of pivot elements which then determine
the selection of rows and columns. Often times pivot elements are chosen
heuristically based on their magnitude relative to other elements in a
particular row or column. The pivot with the globally largest magnitude can
never be determined while maintaining efficient asymptotic speed since each
element of the matrix must be accessed, incurring a cost of~$\cO(MN)$ when
compressing an~$M \times N$ matrix. Depending on the particular algorithmic
choice of ACA, the overall compression cost of an~$N \times N$ rank-structured
matrix may vary from~$\cO(N)$ to~$\cO(N \log^p N)$ to~$\cO(N^2)$.
In~\cite{liu2020parallel}, an algorithm which improves the robustness of ACA
using blocking and hierarchical partitioning is presented in the context of
boundary integral equations and kernel matrices.

\subsection{Randomized black box algorithms}
\label{sec:comp_randblackbox}

In some situations, direct access to the matrix elements is either not available
or computationally expensive. If, however, we have the ability to apply
the given matrix~$\mtx{A}$ to vectors, there exist randomized methods that can
reconstruct the full rank-structured representation (all dense blocks and all
low rank factors) by inspecting sets of input-output pairs. In a typical
environment, a method of this nature constructs two ``test matrices''
$\mtx{\Omega},\mtx{\Psi} \in \mathbb{R}^{N\times \ell}$ and then collects two
``sample matrices'' $\mtx{Y},\,\mtx{Z} \in \mathbb{R}^{N\times \ell}$ through
the formulas
\begin{equation}
\mtx{Y} = \mtx{A}\mtx{\Omega},
\qquad\mbox{and}\qquad
\mtx{Z} = \mtx{A}^{*}\mtx{\Psi}.
\end{equation}
The method then reconstructs all the information required to define $\mtx{A}$
from the information provided in the set
$\{\mtx{Y},\,\mtx{\Omega},\,\mtx{Z},\,\mtx{\Psi}\}$.

A particularly natural application of a randomized black box algorithm is the
compression of a \textit{product} of structured matrices. Directly multiplying
two such matrices together (in their hierarchically compressed formats) is
difficult to execute efficiently in practice, but applying such a product to a
vector is easy. A related application is the compression of Schur complements
that arise in sparse direct solvers.  As mentioned in Section
\ref{sec:comp_bruteforce}, the formulas that define such Schur complements
typically involve a mix of sparse and rank-structured
matrices~\cite{2013_xia_randomized_FDS, ghysels2017robust}.  In the context of
compressing {integral operators} or \textit{products of integral operators}, we
often have access to legacy methods (such as the Fast Multipole
Method~\cite{rokhlin1987} or FFT-based methods) for applying the operator to a
matrix.  In such cases, a randomized black box compression method can be used to
obtain a representation of the operator in a format (such as HBS, cf.~Section
\ref{sec:BIE_nested}) that easily lends itself to direct inversion of the
matrix.

The most efficient randomized black box algorithms pertain to the HBS format we
described in Section \ref{sec:BIE_nested}. (This is perhaps not surprising,
given that the HBS format is also exceptionally efficient for inversion and
for matrix-vector multiplication.) An early
method~\cite{2011_martinsson_randomhudson} demonstrated the ability to attain
high speed and excellent numerical stability, but suffered from not being fully
black box, as this method required access to individual entries of the matrix,
in addition to matrix-vector products. A fully black box version is described in
\cite{levitt2024linear,2022_levitt_dissertation}; another
black-box adaptive scheme amenable to acceleration on GPUs is presented
in~\cite{boukaram2025adaptive}. Thorough and illuminating error bounds
on the appoximations made through these algorithms is provided
in~\cite{amsel2025}; in particular, in most cases, the approximations made are
within logarithmic factors of the best possible fixed-rank HSS or HBS
approximations, as measured in Frobenius norms.

For formats based on strong admissibility, the problem of extracting a data
sparse representation of a matrix from its action on vectors alone was in
principle solved by the ``peeling algorithm'' of \cite{2011_lin_lu_ying}. The
method has been successively improved over the years
\cite{2016_martinsson_hudson2,2022_martinsson_graph_coloring}, but its practical
computational cost remains relatively high. A faster approach based on
``tagging'' of off-diagonal blocks was presented in~\cite{levitt2024linear}, but
this faster version currently applies only to ``flat'' tessellations of the
matrix.  A recent ``data-driven'' alternative based on various sub-sampling
methods is described in~\cite{2024_chow_xi_data_driven_h_matrix}.

An interesting twist on randomized compression in the context of strong
admissibility is the recent technique~\cite{yesypenko2023randomized} that
simultaneously compresses and factorizes a matrix into the \textit{strong
  recursive skeletonization} representation of
\cite{2017_ho_ying_strong_RS}. Some recent work has been done on proving
worst-case approximation estimates to these randomized methods for hierarchical
matrices, an analysis of the HODLR format is presented
in~\cite{chen2025near, 2024_halikias_blackbox_recovery}.  (In practice, for
problems coming from PDEs or integral equations, the approximation accuracy is
many orders of magnitude better than the predicted worst-case scenario.)
Techniques that exploit that the matrices under consideration approximate
continuum operators were presented in
\cite{2023_townsend_boulle_learning_green,2022_townsend_learning_green}.
Integrating knowledge of the kernel matrix into the accuracy estimates is an
area rich for future work.

The reader interested in learning more will find an introduction to randomized
compression of rank-structured matrices in~\cite[Sec.~20]{2020_martinsson_acta}
(with an extended version available in~\cite{martinsson2020randomized}).


%% file: 12-statistics/statistics01.tex
So far in what we have discussed, the ability to numerically compress various
matrices in a hierarchical fashion has been a consequence of where the matrix
originated: as a discretization of a global operator associated with an 
elliptic PDE. In many cases it is easy to
prove that such matrices should, in fact, be compressible based on PDE or
integral equation theory. However, it turns out that many matrices appearing in
data science, machine learning, and statistics have similar compressible
structures despite their entries \emph{not} originating from elliptic
PDEs~\cite{chen2023linear,chen2017hierarchically}. The main category of these
matrices are known as \emph{kernel matrices} with matrix entries given
by~$k(\pxx_{i},\pxx_{j})$, for some kernel function~$k$. Depending on the
application, there may be restrictions put on the kernel function (e.g.
positive definiteness in the case of Gaussian processes) or specific functions
of the kernel matrix that need to be evaluated (e.g. factorizations,
determinants, etc.). In what follows and in order to be concise and concrete,
we give a brief overview of where these kernel matrices appear in the context
of statistical applications and the required operations, as well as some
methods for compressing non-PDE kernels. Additional applications of kernel
matrices exist in machine learning techniques, but a full discussion of all of
these topics is beyond the scope of this section.

Consider a sequence of normal random variables~$\cY_{1}, \ldots, \cY_{N}$ with
zero mean and covariance matrix with entries given
by~$\mtx{C}(i,j) = k(\pxx_{i},\pxx_{j})$.  Here, again to be concrete, we will assume
that the random variable~$\cY_{i}$ corresponds to some probabilistic process
occurring at some point~$\pxx_{i} \in \mathbb R^{d}$, for example the temperature
of air at some point in space~$\pxx_{i}$. Taken together, these random variables
form a multivariate normal distribution whose probability density function is
given by:
\begin{equation}
  p(\vct{y}) = \frac{1}{ (2\pi)^{N/2} \sqrt{  \det \mtx{C}} }
  e^{ - \vct{y}^T \mtx{C}^{{-1}} \vct{y} }.
\end{equation}
Merely evaluating the probability density function above requires the inversion
of the matrix~$\mtx{C}$ and the evaluation of its determinant. This probability
function also appears when performing maximum ($\log$-) likelihood optimizations
and in various regression problems~\cite{williams2006gaussian}.
Furthermore, denoting by~$\cZ_1, \ldots, \cZ_N$ a collection of independent unit
normal random variables, then
 one can generate samples from the distribution~$p$ by noting that
\begin{equation}
  \begin{pmatrix}
    \cY_{1} \\
    \vdots \\
    \cY_{N}
  \end{pmatrix}   \sim
  \sqrt{\mtx{C}}
  \begin{pmatrix}
    \cZ_{1} \\
    \vdots \\
    \cZ_{N}
  \end{pmatrix},
\end{equation}
where~$\sqrt{\mtx{C}}$ denotes the square-root of the matrix~$\mtx{C}$
and~$\sim$ denotes equivalence in distribution. The matrix square-root can be
obtained from an SVD or Cholesky factorization of~$\mtx{C}$.

In order for~$p$ to be a valid probability density function, the
matrix~$\mtx{C}$ must be positive definite.  In the case
where~$k(\pxx,\pyy) = k(\pxx-\pyy)$ (corresponding to an underlying stationary
random process), this constraint can be guaranteed by choosing~$k$ to be a
function whose Fourier transform is positive, a result known as Bochner's
Theorem~\cite{gikhman2004theory}.  Commonly used positive definite kernel functions
(sometimes also known as covariance functions) for stationary random
fields~\cite{williams2006gaussian} include
\begin{equation}
  k(r) = \frac{1}{\left( r^{2}/2\alpha \right)^{\alpha}} , \qquad
  k(r) = e^{-r^{2}}, \qquad
  k(r) = \frac{2^{1-\nu}}{\Gamma{(\nu)}}
  \left( \sqrt{2\nu} r  \right)^{\nu}
    K_{\nu} \left( \sqrt{2\nu} r \right),
\end{equation}
where~$r = |\pxx-\pyy|$ for $\pxx,\pyy \in \bbR^d$, and where it is assumed
that~$\alpha,\nu>0$ and~$K_{\nu}$ is the modified Bessel function of the second
kind~\cite{olver2010nist}. The last family of kernels above on the right are
known as Mat\'ern kernels, and are widely used in spatial statistics.

\subsection{Kernel matrix compression}

Many positive kernels used in computational statistics are not formally Green's
functions for elliptic (or non-elliptic) PDEs, but many admit the same types of
hierarchical matrix
factorizations~\cite{yu2017geometry,rebrova2018study,geoga2020scalable,2020_kriemann_keyes_hlibcov}. The
challenging part is mainly in the computation of the individual low-rank
factors, as there does not exist an underlying Green's theory as in the case for
elliptic PDEs. Some commonly used methods include, among others,
interpolation-based compression~\cite{fong2009black}, proxy point compression
for non-PDE kernels~\cite{xing2020interpolative,minden2017fast}, analytic
separation of variables, and adaptive cross
approximation~\cite{ambikasaran2015fast}.

A commonly used procedure for constructing compressed representations of
these kernel matrices is known as \emph{Vecchia approximation}~\cite{katzfuss2021} and relies on choosing a hierarchy of
points which can then be used to approximately factorize the probability
density~$p$ in terms of a sequence of conditional
probabilities. Algorithmically, Vecchia approximation closely resembles
algorithms based on recursive skeletonization or hierarchical matrix algebra,
but without the numerical approximation guarantees (since there is no underlying
Green's theory).

Lastly, since many covariance functions above are not formally PDE kernels,
there is no reason to believe that they would give rise to matrices for which
nested bases provably exist for the hierarchical partitioning. For this reason,
methods based purely on linear algebra are often preferred in the construction
of such factorizations. For example, using a \emph{modified proxy trick}
in~\cite{minden2017fast} and recursive
skeletonization~\cite{xing2020interpolative}, the authors are able to
obtain~$\cO(N)$ compression and factorization time for some kernel matrices used
in Gaussian process regression.

\subsection{Factorization and inversion}

After computing all the off-diagonal low-rank factors in the matrix~$\mtx{C}$
for whichever hierarchical format has been chosen (e.g. HODLR, HBS, etc.), it
is necessary to compute a factorization of the matrix which lends itself to the
tasks at hand: inversion, square-root factorization, and determinant
calculation. By and large, it is useful to use factorizations which
express~$\mtx{C}$ as a product of matrices instead of a nested sequence of
sums, as we saw, e.g., in Section \ref{sec:HODLRinv}. The reason for this is that
product factorizations naturally lend themselves to inversion and determinant
calculation, and when done in a particular fashion, allow for the symmetric
factorization or square-root factorization of the original
matrix~\cite{chen2014computing,ambikasaran2014fast}.

For example, imposing a HODLR structure on~$\mtx{C}$ and subsequently performing
a multiplicative factorization~\cite{ambikasaran2015fast,
  2013_darve_FDS, minden2017fast}
expresses~$\mtx{C}$ as a product of matrices, each of the
form~$\mtx{I} + \mtx{U}\mtx{V}^{T}$, where~$\mtx{U}$ and~$\mtx{V}$ are of
low-rank. The determinant of such a factorization can be computed easily using
Sylvester's Determinant Theorem~\cite{akritas1996various} on each
factor. Sylvester's Determinant Theorem states that
\begin{equation}
  \det (\mtx{I} + \mtx{U}\mtx{V}^{T} ) = \det (\mtx{I}+ \mtx{V}^T\mtx{U} ).
\end{equation}
If the original matrix~$\mtx{I} + \mtx{U}\mtx{V}^{T}$ was of size~$N \times N$,
and the low-rank factors are of rank~$k$, then the cost of~\emph{directly}
computing the determinant is reduced from~$\cO(N^3)$ to~$\cO(N\,k^2)$, where
hopefully~$k \ll N$.

\subsection{Future directions}

Kernel matrices are appearing more and more due to the adoption of machine
learning and statistical methods in the sciences, and there remains a lot of
work to do in order to efficiently manipulate them. On the one hand, many of
these kernels are smooth at the origin -- or at least have a few continuous
derivatives -- which allows for very efficient compression when the data
points~$\pxx_i$ lie in relatively low dimensions. However, often times,
particularly in applications appearing in neural network models, the ambient
dimension is quite high and compression becomes much more challenging if not
impossible \cite{yu2017geometry}.

Furthermore, the kernel generating these matrices can depend on several
parameters, known as \emph{hyperparameters}. Frequently an important step in the
overarching statistical calculation is to optimize these parameters so as to,
for example, perform a maximum likelihood calculation. In this case, what is
needed is the gradient of the kernel matrix~$\mtx{C}$ with respect to these
parameters. Rapidly computing an update to the matrix~$\mtx{C}$ using analytic
gradient information of the kernel and an existing hierarchical factorization is
ongoing research in many research groups.


%% file: 13-taxonomy/taxonomy03.tex
At this point, this survey has introduced a long list of different rank
structured formats, and the ones included are only a small fraction of what has
been proposed in the literature. This raises the question: Which technique or
approach should one choose for a given application? For problems in one and two
dimensions, it is typically not hard to find methods that are simple to
implement, fast in practice, and of linear or close-to-linear asymptotic
complexity. However, methods that attain linear complexity for problems in three
dimension tend to be much more complex. In practice, it may be worth
sacrificing an optimal~$\cO(N)$ algorithm in favor of one that scales
as~$\cO (N^{3/2})$, for example, due to the maximum~$N$ under consideration and
the constant implicit in the~$\cO (\cdot)$ notation. In what follows we try to
provide a rough guide for selecting methods for some canonical problems, keeping
in mind the desired goal of some standard ``users''. We also try to provide
pointers to some existing commonly-used software for each use case.

\subsection{PDE discretizations}

PDE discretizations can include finite differences, finite elements, or spectral
element methods. In all cases, the resulting matrix is sparse (element-wise or
block-wise sparse, depending on the particular type of discretization).  The
complexity of classical direct solvers that construct a triangular factorization
of the coefficient matrix depends on the ordering of the nodes in the mesh, as
discussed in Section~\ref{sec:FDSPDE}.  In most cases, a ``nested dissection''
ordering is optimal, and leads to complexities of~$\cO(N^{3/2})$
and~$\cO(N^{2})$ in two and three dimensions, respectively.  We recall that by
exploiting rank structure in the dense Schur complements that arise, linear
complexity can be attained in both two and three dimensions.

Due to the inherent complexity of the double hierarchical data structures that
are required, it would rarely make sense for a user to invest the work required
to code a fast direct solver from scratch in this environment.  Instead, we
would suggest investigating existing software packages such as
STRUMPACK~\cite{10046092} or MUMPS~\cite{2019_amestoy_BLR_multifrontal}, linked
to below.  Related methods based on spectral element discretizations, and known
as ``hierarchical Poincar\'e-Steklov'' methods, have been shown to scale quite
efficiently when coupled with GPU acceleration in two and three
dimensions~\cite{melia2025hardware}. These methods have also been shown to be
compatible with unstructured meshes~\cite{fortunato2021ultraspherical}.

\subsection{Boundary integral equations}

Due to the required numerical machinery for discretizing boundary integral equations
as compared to finite differences, for example, the available software is
somewhat more limited, and currently under development by several groups.

\paragraph{Curves}

Since weakly-admissible compression and inversion schemes lead to linearly
scaling algorithms for boundary integral equations on curves in two dimensions, 
there is little reason to deploy algorithms based on strong admissibility 
conditions. It may, however, be advantageous to use a scheme which takes 
advantage of nested bases in order to reduce the memory footprint of the 
compressed matrix. (The challenges in the ``high frequency'' regime of
course remain, cf.~Remark \ref{remark:highfreq}; however, simple methods
remain effective up to domains that are hundreds of wavelengths across.)

\paragraph{Surfaces}
For medium-scale problems, direct solvers for boundary integral equations
defined along surfaces based on weak admissibility provide an attractive and
straightforward option since they are based on simple and exact inversion
formulas. However, such methods scale as~$\cO(N^{3/2})$ at best for problems in
three dimensions.

For large-scale problems, methods that are based on strong admissibility data
structures become necessary, see Section~\ref{sec:BIE_strong}. Techniques based
on $\mathcal{H}$- and $\mathcal{H}^{2}$-matrices are commonly used, and several
software packages are available (see Section~\ref{sec:libraries}). Due to the
need for nested recursions and re-compression, techniques of this nature become
quite expensive if high accuracy is sought and are commonly used in practice at
medium or low accuracy to form pre-conditioners for iterative solvers. 
Emerging techniques that avoid recursion -- such as recursive strong skeletonization and inverse fast multipole methods -- offer a promising path for building methods that remain performant even at higher accuracies. 
Efforts to integrate such methods with high-order quadratures and their implementation in HPC regimes is currently under way.

\paragraph{High performance algorithms}
Of course, for very large~$N$, distributed computing algorithms must be used in
order to overcome the memory constraints of a single workstation
\cite{rouet2016distributed, ma2024inherently, yamazaki2019distributed,
  guo2015mpi}.  Work in this area is relatively thin when compared to the
literature overall, but some results indicate that schemes based on slightly
more complicated data structures~\cite{li2017distributed,2015_betcke_BEMpp}
yield linearly scaling schemes. The idea here is to involve an offset data
structure which can effectively re-compress interactions whose ranks may have
grown due to the geometry of the problem. A more recent GPU-compatible
implementation~\cite{boukaram2025linear} of the strong recursive skeletonization
algorithm~\cite{2017_ho_ying_strong_RS} shows promising scaling results for
problems of size up to about~$N = 10^{6}$.

\subsection{Volume integral equations}

Discretizations of volume integral equations result in a collection of
space-filling nodes, with no geometric separation between them in any direction.
Weakly admissible schemes for compressing and inverting volume integral
equations, such as the Lippmann-Schwinger integral equation, will lead to
algorithms which scale as~$\cO(N^{3/2})$ in two dimensions and~$\cO(N^{2})$ in
three dimensions, where~$N$ is the total number of unknowns in the
discretization. As in the case of boundary integral equations along surfaces,
often times the~$\cO(N^{3/2})$ scaling schemes in two dimensions are sufficient
due to the size of~$N$ and the relatively small
constants implicit in this notation for weakly admissable schemes \cite{gopal2022accelerated}. 

However, in three dimensions, not only is~$N$ generally bigger but the scaling
of~$\cO(N^2)$ is worse, and has an even larger constant than in the two
dimensional case. It is likely necessary, in this regime, to spend the effort to
design solvers based on strong admissability conditions that will lead to
linearly scaling scheme, despite the larger pre-factor in the asympotic
complexity analysis~\cite{2017_ho_ying_strong_RS}. Direct solvers for volume
problems have recently been developed using an alternative, butterfly
compression scheme~\cite{2010_rokhlin_oneil_butterfly,
  michielssen1996multilevel}, in order to address higher frequency oscillatory
problems, such as the Helmholtz equation~\cite{sayed2021butterfly}.

It is likely that the distributed parallel algorithm of~\cite{liang2024n} can be
used to construct an efficient linearly scaling fast direct solver for volume
integral equations in two and three dimensions, but this is ongoing research and
requires a significant software development effort.

\subsection{Kernel methods in statistics and machine learning}

When the ambient dimension of the points~$\pxx_i$ responsible for generating a
kernel matrix with entries~$\mtx{C}(i,j) = k(\pxx_i,\pxx_j)$ is small, the
decision as to what factorization method to use is usually dictated by whether
or not inversion or determinant calculation is needed. Compression based on an
analytic or interpolative analysis of the kernel function, coupled with
weakly-admissable data structures tend to balance the trade-off between speed
and simplicity of algorithm.  Several codes for efficient Gaussian process
calculations exist which exploit hierarchical matrix techniques.

However, if the ambient dimension is large, competitive algorithms based on
Vecchia approximations (without rigorous error bounds) are likely the most
efficient. Algorithms for high dimensional data based on nearest-neighbor
compression have also been developed~\cite{chenhan2016inv, chavez2020scalable,
  march2015askit, march2016askit}.

\subsection{Software libraries}
\label{sec:libraries}
Here we list several available software libraries related to the topics
discussed in this review. Keep in mind that the computational infrastructure
surrounding these methods is growing and in flux, and this list is updated only
as of late 2025.

\vspace{\baselineskip}
\textbf{Fast linear algebra}:
\begin{itemize}
\item \texttt{libflame}: Fast dense linear algebra\\
  \url{https://github.com/flame/libflame/}
\end{itemize}

\vspace{\baselineskip}
\textbf{Fast multipole methods and kernel methods}:
\begin{itemize}
\item \texttt{BBFMM3D}: Black-box FMM in three dimensions\\
  \url{https://github.com/ruoxi-wang/BBFMM3D}
\item \texttt{exafmm-t}: Shared memory parallel FMM\\
  \url{https://github.com/exafmm/exafmm-t}
\item \texttt{fmm2d}: Fast multipole methods in two dimensions\\
  \url{https://github.com/flatironinstitute/fmm2d}
\item \texttt{FMM3D}: Fast multipole methods in three dimensions\\
  \url{https://github.com/flatironinstitute/fmm3d}
\item \texttt{H2Pack}: $\mathcal H^2$ kernel matrix compression\\
  \url{https://github.com/scalable-matrix/H2Pack}
\end{itemize}

\vspace{\baselineskip}
\textbf{Fast solvers for sparse matrices}:
\begin{itemize}
\item \texttt{MUMPS}: A parallel sparse direct solver\\
  \url{https://mumps-solver.org}
\item {\texttt{cuDSS}: NVIDIA GPU-accelerated sparse direct solver}\\
  \url{https://developer.nvidia.com/cudss}
\item \texttt{STRUMPACK}: Linear algebra for sparse and structured matrices\\
  \url{https://portal.nersc.gov/project/sparse/strumpack}
\item \texttt{surfacefun}: Hierarchical solver for spectral elements on surface
  PDEs~\cite{fortunato2024high}\\
  \url{https://github.com/danfortunato/surfacefun}
\end{itemize}

\vspace{\baselineskip}
\textbf{Hierarchical matrix algebra}:
\begin{itemize}
\item \texttt{AHMED}: $\mathcal{H}$-matrix methods for elliptic differential equations\\
    \url{https://www.wr.uni-bayreuth.de/en/software/ahmed/index.html}
  \item \texttt{FLAM}: Recursive skeletoniztion hierarchical linear algebra algorithms~\cite{ho2020flam}\\
  \url{https://github.com/klho/FLAM}
  \item {\texttt{HLIB}: Hierarchical matrix library, superceded by H2Lib\\
  \url{https://gitlab.mis.mpg.de/scicomp/hlib}}
\item \texttt{H2Lib}: Hierarchical matrix library\\
  \url{http://www.h2lib.org}
  \item {\texttt{HLIBpro}: H-matrix algorithms, parallelized\\
  \url{https://hlibpro.com/}}
\item {\texttt{hm-toolbox}: HODLR and HSS matrix arithmetic in MATLAB~\cite{massei2020hm}\\ 
\url{https://github.com/numpi/hm-toolbox}}
\item \texttt{hodlr\_gpu}: GPU accelerated HODLR matrix algebra\\
  \url{https://bitbucket.org/chao_chen/hodlr_gpu}
\item \texttt{HODLRlib}: OpenMP accelerated KD-tree based HODLR matrix algebra\\
  \url{https://github.com/SAFRAN-LAB/HODLR}
\item \texttt{strong-skel}: Strong skeletonization matrix factorization \\
    \url{https://github.com/victorminden/strong-skel}
\end{itemize}

\vspace{\baselineskip}
\textbf{Boundary integral equation solvers}:
\begin{itemize}
\item \texttt{bempp-cl}: Boundary element method Python package\\
  \url{https://github.com/bempp/bempp-cl}
\item \texttt{chunkie}: Fast boundary integral equations in 2D\\
  \url{https://github.com/fastalgorithms/chunkie}
\item \texttt{FMM3DBIE}: FMM-accelerated boundary integral equations in 3D\\
  \url{https://github.com/fastalgorithms/fmm3dbie}
\item \texttt{strong-skel}: Strong skeletonization-based direct solver for surfaces in
  3D\\ 
  \url{https://github.com/fastalgorithms/strong-skel}
\item \texttt{pytential}: FMM and HSS-based direct solvers in Python\\
  \url{https://github.com/inducer/pytential}
\end{itemize}

\vspace{\baselineskip}
\textbf{Statistics and kernel methods}:
\begin{itemize}
\item \texttt{george}: Fast and flexible Gaussian Process regression in Python\\
  \url{https://github.com/dfm/george}
\item \texttt{ASKIT}: HPC implementation of high-dimensional tree codes \\
  \url{https://padas.oden.utexas.edu/libaskit}
\item \texttt{GPMLE}: Fast Gaussian Process MLE using skeletonization \\
  \url{https://github.com/victorminden/GPMLE}
\item {\texttt{HiGP}: High-performance Gaussian Processes with hierarchical matrix support\\
  \url{https://github.com/huanghua1994/HiGP}}
\item {\texttt{RLCM}: Recursively Low-Rank Compressed Matrices\\
  \url{https://github.com/jiechenjiechen/RLCM}}
\item \texttt{Vecchia.jl}: Vechia approximations to the Gaussian likelihood\\
  \url{https://github.com/cgeoga/Vecchia.jl}
\item \texttt{HLIBCov}: $\mathcal{H}$-matrix arithmetic for covariance matrices\\
  \url{https://github.com/litvinen/HLIBCov}
\end{itemize}


%% file: 14-conclusion/conclusion02.tex
This survey aims to highlight a collection of techniques that have enabled the fast direct solution of many linear systems in scientific computing
that are intractable or highly challenging for iterative solvers.

The development of algorithms such as fast direct solvers has also uncovered
persistent ubiquitous analytical structures in a variety of global operators
associated with many of the fundamental PDEs of mathematical physics. The
underlying physics manifests itself as the rapid diffusion of information over
modest distances, and is present in solution operators for elliptic PDEs, time
evolution operators for parabolic PDEs, boundary-to-boundary operators such as
the Dirichlet-to-Neumann operator, and numerous other contexts. Through
numerical examples, this survey demonstrates the relatively detailed
understanding of the factors governing the rate at which information diffuses
across distances, and shows how these insights can be exploited to construct
powerful and efficient numerical algorithms.

Currently, the principles for constructing direct solvers with linear (or
near-linear) complexity scaling in both two and three dimensions are well
established. Nevertheless, developing (and implementing) solvers that are
sufficiently fast for practical application to complex three-dimensional
geometries remains ongoing work. We also believe substantial room exists for
reducing algorithmic complexity -- many methods that achieve
\textit{theoretical} linear complexity often prove too intricate to implement
efficiently in realistic settings. Recent years have witnessed significant
advances in this area through the development of randomized methods for ``black
box'' compression of various global operators.

In the regime where we fix a PDE and refine the computational discretization,
linear scaling fast direct solvers have been developed for most constant
coefficient problems and some associated variable coefficient problems
(e.g. acoustic scattering through media with variable sound-speed).  However, a
major outstanding challenge in the field involves oscillatory problems where the
number of points per wavelength is held constant, e.g. a high-frequency
time-harmonic electromagnetics problem.  Methods for these problems based on
iterative solvers and algorithms such as the High-Frequency
FMM~\cite{2006_rokhlin_wideband} often require a prohibitive number of
iterations to achieve convergence.  Recent results and evidence suggests that
many techniques described in this survey may achieve linear complexity by
replacing ``$\mathcal{H}$-matrix structures'' with ``butterfly
structures''~\cite{2016_michielssen_butterfly,2021_li_butterfly}. Research in
this direction is still at a relatively early stage, with only one or two groups
having research-grade codes~\cite{doecode_21501}.

Another area where our understanding of fast direct solvers is not yet fully
satisfactory is that of numerical stability. From a practical standpoint, the
techniques described in this survey are stable and have already see widespread
use in real-world applications. However, rigorous theoretical guarantees are
typically absent, and most existing methods possess vulnerabilities that could
theoretically trigger unacceptable error growth. Fortunately, these failure
modes simply tend not to arise in practice. Even without provably stable
techniques and rigorous error bounds, practitioners can continue to confidently
employ existing methods provided they incorporate validation procedures that
verify computed solutions meet the required accuracy standards upon
completion. Moreover, a primary application of fast direct solvers recently has
been in the construction of approximate inverses for preconditioning
purposes. In this context, the outer iteration stops only once a predetermined
accuracy threshold is satisfied, so convergence itself ensures that full
precision has been attained in the final solution, regardless of the exact
accuracy of the approximate inverse.

Given the state of research and practical use of fast direct solvers, let us
offer a prediction: While the most obvious utility of fast direct solvers will
continue to be in solving linear systems intractable to iterative methods, an
equally important use case will prove to be the modeling of multiphysics and
coupled problems whose robust formulations require non-trivial use of operator
algebra. For example, the ability to directly form boundary-to-boundary maps
(e.g.~Dirichlet-to-Neumann or Impedance-to-Impedance) or to evaluate functions
of operators enables the formulation of compact and precise equilibrium
equations defined at domain interfaces that can themselves be solved directly.

In closing, let us observe that the techniques described in this survey would
not have proven feasible without the plummeting cost of both floating-point
operations and fast storage availability over the past couple of decades. Given
their communication efficiency, fast direct solvers built on hierarchical data
structures are at this point highly competitive in a myriad of important
environments, and we predict that as most computational infrastructure continues
to become more communication-constrained the balance will keep shifting in their
favor.